\documentclass{aims} 
\usepackage{amsmath}
\usepackage{amssymb, bm, bbm}
\usepackage{esint}
\usepackage[frak]{paper_diening}
\usepackage{paralist}
\usepackage{enumitem}
\usepackage{epsfig} 
\usepackage{epstopdf} 
\usepackage[colorlinks=true]{hyperref}
\hypersetup{urlcolor=blue, citecolor=red}
\allowdisplaybreaks

\def\currentvolume{X}
 \def\currentissue{X}
  \def\currentyear{200X}
   \def\currentmonth{XX}
    \def\ppages{X--XX}
     \def\DOI{10.3934/xx.xxxxxxx}

\newtheorem{theorem}{Theorem}[section]
\newtheorem{corollary}{Corollary}

\newtheorem{lemma}[theorem]{Lemma}
\newtheorem{proposition}{Proposition}

\theoremstyle{definition}

\newtheorem{remark}{Remark}

\renewcommand{\Xi}{X}

\renewcommand{\phi}{\varphi}

\newcommand{\RN}{\setR^N}

\newcommand{\kappaA}{{\kappa_{\!\mathcal{A}}}}

\newcommand{\bg}{{\bf g}}
\newcommand{\bP}{{\bf P}}
\newcommand{\bu}{{\bf u}}

\newcommand{\meantmp}[2]{#1\langle{#2}#1\rangle}

\def\la1{\lambda_1}

\newcommand{\bT}{{\bf T}}

\newcommand{\F}{\mathcal{F}}

\newcommand{\R}{\ensuremath{\mathbb{R}}}

\newcommand{\dd}{\mathrm{d}}

\def\la1{\lambda_1}

\newcommand{\real}{\mathbb R}

\newcommand{\bQ}{{\bf Q}}

\newcommand{\kabs}[1]{\ensuremath{\vert#1\vert}}

\title[Boundary partial regularity for minimizers of discontinuous quasiconvex integrals]
{Boundary partial regularity for minimizers of \\ discontinuous quasiconvex  integrals with  general growth} 

\author[J. Ok, G. Scilla and B. Stroffolini]{}


 \email{jihoonok@sogang.ac.kr}
 \email{giovanni.scilla@uniroma1.it}
 \email{bstroffo@unina.it}

\thanks{$^*$Corresponding author: Jihoon Ok}

\begin{document}
\maketitle

\centerline{\scshape Jihoon Ok$^*$}
\medskip
{\footnotesize
 \centerline{Department of Mathematics, Sogang University}
   \centerline{35 Baekbeom-ro, Mapo-gu,  04107 Seoul}
   \centerline{Republic of Korea}
} 

\medskip

\centerline{\scshape Giovanni Scilla}
\medskip
{\footnotesize
 \centerline{Dipartimento di Scienze di Base ed Applicate per l'Ingegneria (SBAI)}
   \centerline{Sapienza Universit\`{a} di Roma, Via A. Scarpa 16, 00169 Roma}
   \centerline{Italy}
}

\medskip

\centerline{\scshape Bianca Stroffolini}
\medskip
{\footnotesize
 \centerline{Dipartimento di Ingegneria Elettrica e delle Tecnologie dell'Informazione}
   \centerline{Universit\`{a} di Napoli Federico II, Via Claudio, 80125 Napoli}
   \centerline{Italy}
}

\bigskip



\begin{abstract}
We prove the partial H\"older continuity on boundary points for minimizers of quasiconvex non-degenerate functionals 
\begin{equation*}
\mathcal{F}({\bf u}) \colon =\int_{\Omega} f(x,{\bf u},D{\bf u})\,\mathrm{d}x\,,
\end{equation*}
where $f$ satisfies a uniform VMO condition with respect to the $x$-variable, is continuous with respect to ${\bf u}$ and has a general growth with respect to the gradient variable.
\end{abstract}


\section{Introduction}

In this paper we study the boundary partial regularity of minimizers of the following non-autonomous integral functional
\begin{equation}\label{functional}
\F({\bf u}) \colon =\int_{\Omega} f(x,{\bf u},D{\bf u})\,\dd x,
\end{equation}
where $\Omega\subseteq\mathbb{R}^{n}$ ($n\ge 2$) is an open bounded set of class $C^1$ and $\bfu : \Omega\rightarrow\mathbb{R}^{N}$ with $N\ge2$  --  i.\,e., we consider vectorial minimizers. We assume that the integrand $f:\Omega\times \R^N\times R^{N\times n}$, $f(x,\bu,\bP)$, satisfies a general growth and a quasiconvexity condition with an $N$-function $\phi$. Moreover, we assume that $f$ complies with a Vanishing Mean Oscillation (VMO) with respect to the spatial variable $x$, hence it can be discontinuous. The specific assumptions on the function $f$  will be introduced in Section~\ref{sec1.1}.

\emph{A brief review of the literature on the topic.} The regularity theory for functionals of the form \eqref{functional} and related partial differential equations are one of the classical topics in the calculus of variations. The study of functionals/equations with general growth has been initiated by Marcellini and developed in a series of  seminal papers,  for example \cite{ MARCELLINI93, MARCELLINI96}.  A first $C^{1,\alpha}$ local regularity for  solutions of Euler-Lagrange elliptic systems  depending on the modulus of the gradient with general growth under suitable hypotheses was obtained by Marcellini\&Papi, \cite{MARPAPI06}.

Subsequently, for an autonomous integrand $f$ satisfying the Uhlenbeck structure~\cite{UHLENBECK77} - that is, $f(x,\bu,\bP)\equiv \varphi(|\bP|)$ and $\varphi$ is sufficiently regular and convex,  complying with \ref{ass-1phi}-\ref{ass-2phi} in Section~\ref{sec1.1} below and such that $\varphi''$ is H\"older continuous 
{off the diagonal}
 -  the minimizer of \eqref{functional} was shown to be  locally $C^{1,\alpha}$ with an excess decay estimate, see \cite{DIESTROVER09}. {Within this setting, we consider it appropriate to mention the result of global Lipschitz regularity 
for solutions to boundary problems for elliptic systems of Uhlenbeck-type obtained in \cite{CM}, where the domain is convex and its boundary is weaker than Lipschitz (see \cite[Theorem~2.1]{CM}). Going back to integral functionals,} 
if $f$ is autonomous but does not satisfy the Uhlenbeck structure, i.\,e., $f(x,\bu,\bP)\equiv g(\bP)$, then there are examples of regular and convex functions $g$ whose minimizer is locally unbounded; see, for instance, the survey paper of Mingione \cite{Mingione06Dark}. Therefore, only partial regularity is expected if we consider a general non-Uhlenbeck structure. Here, ``partial regularity" means the minimizer $\bu$ satisfies a desired regularity except a measure zero set.

Partial regularity results for quasi-convex functionals with $p$-growth, i.\,e., $\phi(t)=t^p$ in Section~\ref{sec1.1}, had been first studied in \cite{AcerbiFusco1,CarFuscoMin1,Evans1} by using  the blow-up technique, which yields decay estimates for the  so-called excess  functional. In particular, in \cite{AcerbiFusco1}, the integrand $f=f(x,\bu,\bP)$ is assumed to be H\"older continuous in $x$ and $\bu$, and the partial $C^{1,\alpha}$-regularity  is proved. Later, Foss and Mingione \cite{fossmingione1} considered the integrand $f=f(x,\bu,\bP)$ that is only continuous in $x$ and $\bu$, and proved the partial H\"older regularity, which means partial $C^{0,\alpha}$-regularity for every $\alpha\in(0,1)$, by using a different approach. They worked with an hybrid excess functional and used the $\mathcal A$-harmonic approximation and Ekeland's variational principle. The $\mathcal A$-harmonic approximation was  introduced and applied to the partial regularity in \cite{DuzaarSteffen1}. The continuity assumption for $x$ was extended to the VMO condition in \cite{BDHS}. 

We note that the partial regularity results mentioned above consider non-degenerate functionals, which means that the integrand $f$ satisfies $0<|D^2f(x,\bu,{\bf 0})|<\infty$. 
For degenerate functionals with $p$-growth, Duzaar and Mingione \cite{DUMIN04b} obtained the partial $C^{1,\alpha}$-regularity in the autonomous case by assuming an additional condition concerning the behavior of  $D^2f(\bP)$ near ${\bf 0}$ and using the $p$-harmonic approximation \cite{DUMIN04}. Then, using the same argument, B\"ogelein \cite{Bogelein} proved the partial H\"older regularity for non-autonomous problems with the VMO condition for $x$. We also refer to \cite{HAB13,Ok1,Ok2,Ok4, BeMin,  cristianarosario} for related partial regularity results for functionals or elliptic systems with variable growth conditions, e.g., $p(x)$-growth, $L^{p(\cdot)}\log L$-growth and double phase growth.

Partial regularity results for quasiconvex functionals with $p$-growth, including the papers mentioned above, deal with the superquadratic case ($p\ge 2$) and the subquadratic case ($1<p<2$) in different ways. They treat only one of the two cases, or the two cases separately, with different regularity assumptions on $D^2f$ and different approaches. However, in the general growth case with  $\phi$ satisfying \ref{ass-1phi}-\ref{ass-2phi}, the functional can be neither superquadratic nor subquadratic since it is possible that $1<\mu_1<2<\mu_2$ in   \ref{ass-2phi}, and the blow-up technique is doomed to failure. 
Therefore, a unified assumption and an approach are required in this general setting. To this aim, Diening, Stroffolini and Verde~\cite{DIESTROVER10} considered degenerate autonomous functionals (i.\,e., $f(x,\bu,\bP)\equiv g(\bP)$ in \eqref{functional}) with general growth and a unified assumption on $D^2f(\bP)$. They introduced the $\varphi$-harmonic approximation - i.e., the counterpart of the $p$-harmonic approximation in the setting of nonstandard growth - and derived excess decay estimates with a new excess functional in terms of shifted $N$-functions, which implies the partial $C^{1,\alpha}$-regularity. As for the quasiconvexity for general growth, in \cite{DIELENSTROVER12} an improved version of the $\mathcal A$-harmonic approximation  in the Orlicz setting has been introduced by using the Lipschitz truncation technique and  a duality argument.
In addition, they deduced via interpolation global Calder\'on-Zygmund estimates in Orlicz spaces for $\mathcal A$-harmonic maps in balls.
This result has been extended to non-autonomous functionals in \cite{CeladaOk,GoodSciStroff2021,Stroffo2020}. In particular, in \cite{GoodSciStroff2021},  partial H\"older-regularity  is proved for non-autonomous functionals with general growth and VMO cefficients.

Partial $C^{1,\alpha}$-regularity on boundary points for quasiconvex functionals with $p$-growth when $p\ge2$ was studied by  Grotowski \cite{Grot1} and Hamburger \cite{Ham1}, where they obtained the boundary versions of $\mathcal A$-harmonic approximation and the blow-up technique. Later, 
Beck \cite{Be09} proved partial H\"older regularity results by using the approach in \cite{fossmingione1}, and this result was extended to elliptic systems with superquadratic general growth, i.\,e., $\mu_1\ge 2$ in \ref{ass-2phi}, and VMO coefficients in \cite{Ok3}. Here, we point out that the previous results are concerned with the superquadratic case. On the other hand, Beck \cite{Be11} also considered partial $C^{1,\alpha}$-regularity on the boundary in both the subquadratic and the superquadratic cases. We recall also the contribution by Campanato \cite{Campanato}, where assuming merely the continuity of the coefficients with respect to $(x,u)$ without any further structural assumptions, he proved that the weak solution $u$ is H\"older continuous with every exponent $0<\alpha<1$ up to the boundary outside a negligible set, for the low dimensional case, $n\leq p + 2$. We also refer to \cite{BoDuMin101,BoDuMin102} for boundary partial regularity for parabolic systems.
 
As for the Hausdorff dimension of the singular set, the only known results are for integrands $f$ which are convex with respect to $D\bfu$ and H\"older continuous with respect to $(x,\bfu)$, \cite{KRIMIN06}, or for Lipschitz minimizers of quasiconvex functionals, \cite{ KRIMIN07}. Regarding the boundary regularity it is worth mentioning the paper \cite{KRIMIN10} where for regular integrands strongly convex in $D\bfu$ it is proved that $\mathcal{H}^{n-1}$-almost every boundary point is regular. In addition, using the method introduced in \cite{KRIMIN06} about fractional differentiability, there are results for the general growth case in \cite{DIEETT08, cristiana}. 
Summarizing, the existence of regular boundary points so far has been proved only for special structures \cite{JosMei83,DuzGrotKron04} or $\alpha$-H\"older continuous coefficients with $\alpha>\frac12$, see \cite{DuzKriMin07,KRIMIN10,Beck11b}. Therefore, it is still an open problem even for systems or quasiconvex functionals with H\"older continuous coefficients with small H\"older exponents and standard $p$-growth. We further refer to \cite[Section 6]{Beck11b} for existence of regular boundary points.

\par
\emph{Description of our results.} In this paper, we consider partial H\"older regularity on boundary points for quasiconvex functionals with general growth in the gradient and VMO in the $x$-variable that can be neither superquadratic nor subquadratic, and characterize the set of regular boundary points. More precisely, we prove that if the boundary and the boundary datum are of class $C^1$ then the minimizer $\bu$ of the functional \eqref{functional} satisfying the general growth condition in  Section~\ref{sec1.1} is locally H\"older continuous for every H\"older exponent $\alpha\in(0,1)$ at  any  boundary point that is Lebesgue type, in some sense, with respect to $D\bu$. Moreover, we assume that the functional is non-degenerate, and this allows us to convert the original functional with $C^1$ boundary datum to a functional with the zero boundary value.   
 
Finally we briefly comment on the strategies adopted in this paper. We consider minimizers of functionals in half balls that have the  zero values on the flat boundaries.  
Then we try to follow the approach in \cite{GoodSciStroff2021}, where, however, the main techniques consider balls as domains, and
this does not allow to apply directly these results  to the boundary case.
To overcome this difficulty, when possible, we employ the zero and the odd extensions. 
Such extensions cannot be applied to the $\mathcal A$-harmonic approximation step in Section~\ref{subsection:Aharmonic}, since the odd extension of an $\mathcal A$-harmonic map is the solution to a (in general) different homogeneous linear equation with measurable coefficients and the zero extension is more complicated. {Nonetheless, simple computations show that the method of odd reflection turns to be useful for systems with a special structure, as diagonal and Uhlenbeck-type systems.} 
{{A further issue to face is that differently to what happened in \cite{DIELENSTROVER12}, global Calder\'on-Zygmund type estimates 
 for $\mathcal A$-harmonic maps in half balls are 
not clearly known
since half balls are not domains of class $C^1$ but just Lipschitz domains {(see, e.g., \cite[Theorem~A]{JK} for some well-known negative examples about the inhomogeneous Dirichlet problem for the Laplace equation in Lipschitz domains)}. Therefore, we consider global Calder\'on-Zygmund estimates on not half-balls but  relevant $C^1$-domains  (see Theorem~\ref{thm:thm18}).
}} 
{Finally, we emphasize that in earlier papers singular sets were defined in terms of integrals involving shifted $N$-functions $\phi_a$ or the vector-valued function $\bf{V}$ (see Section~\ref{sec:basicNfunctions}). Moreover, in the subquadratic case, the singular set of boundary points had a complicated structure (see \cite[Theorem 1.1]{Be11}). On the contrary, in this paper we define it in terms of $L^1$ integrals, hence the singular set becomes more clear and smaller. This improvement is possible by deriving the reverse H\"older type estimates in \eqref{highint1} and \eqref{eq:caccioppoliIbis1}.}\\

\emph{Outline of the paper.} The paper is organized as follows. In Section~\ref{sec1.1} we specify the main assumptions on $\varphi$ and the integrand $f$, and we state the main result of the paper, Theorem~\ref{theorem-result-1}. In Section~\ref{sec:prelbas} we fix the main notation that we will use throughout the paper and recall some basic definitions and preliminary results. Specifically, Sections~\ref{sec:basicNfunctions}-\ref{sec:orliczsob} contain some facts about $N$-functions and Orlicz-Sobolev spaces, while Section~\ref{sec:somelemmas} collects few standard technical lemmas. Finally, in Section~\ref{subsection:Aharmonic}, we prove the $\mathcal{A}$-harmonic approximation lemma on upper half balls. In Section~\ref{sec:caccioppoliandhigh} we derive Caccioppoli inequalities and higher integrability results for the minimizer of our functional, and we introduce the relevant objects, as the excess functionals, that will come into play later. With Section~\ref{sec:partialreg} we start with the analysis of the partial regularity for the minimizer. In particular, in Section~\ref{sec:approxharmon} we prove the necessary result of almost $\mathcal{A}$-harmonicity in order to obtain the excess decay estimate in Section~\ref{sec:nondegenerate}. The last Section~\ref{lastsection} is entirely devoted to the proof of Theorem~\ref{theorem-result-1}.

\subsection{Assumptions and statement of the main result}
\label{sec1.1}
{We first list the main assumptions on the integrand $f$ in \eqref{functional}. 
As mentioned above, we consider a non-degenerate quasi-convexity condition
with general growth. We note that such conditions in the degenerate setting can be found, for instance, in \cite{DIELENSTROVER12} for autonomous functionals, i.e. $f(x,\bu,\bP)\equiv f(\bP)$, and in \cite{GoodSciStroff2021} for general functionals with VMO condition in $x$. Here the degenerate condition can be obtained from assumptions \ref{ass-1f}--\ref{ass-6f} by removing the ``$+1$". We also refer to \cite{CeladaOk,Stroffo2020}.
}

We start by introducing an $N$-function  $\varphi:[0,\infty)\to[0,\infty)$, and we assume that $\varphi$ satisfies 
\begin{enumerate}[font={\normalfont},label={($\varphi$\arabic*)}]
\item $\varphi\in C^1([0,\infty))\cap C^2((0,\infty))$; \label{ass-1phi}
\item $0<\mu_1-1\leq \inf_{t>0}\frac{t\varphi''(t)}{\varphi'(t)}\leq \sup_{t>0}\frac{t\varphi''(t)}{\varphi'(t)}\leq\mu_2-1$, for suitable constants $1<\mu_1\leq\mu_2$. \label{ass-2phi}
\end{enumerate}
Without loss of generality, we shall assume  that $1<\mu_1<2<\mu_2$. The properties of $\varphi$ complying with the preceding assumptions will be introduced in Section~\ref{sec:prelbas}.

Let the  integrand $f:\Omega\times\R^N\times\mathbb{R}^{N\times n}\to\R$, $f=f(x,{\bf u},{\bf P})$, in \eqref{functional} be Borel-measurable and the partial map ${\bf P}\to f(\cdot,\cdot,{\bf P})\in  C^2(\mathbb{R}^{N\times n})$. We denote by $Df$ and $D^2f$ the corresponding first and second gradients, respectively, for fixed $x$ and $\bf u$. Then, with the $N$-function $\varphi$ and some constants $0<\nu\le L$, $f$ is assumed to satisfy the following conditions:

\begin{enumerate}[font={\normalfont},label={(F\arabic*)}]
\item {\emph{($\varphi$-growth condition I)}}
\label{ass-1f} for every $x\in \Omega$, ${\bf u}\in \R^N$ and ${\bf P}\in\mathbb{R}^{N\times n}$
\begin{equation*}
\nu \phi(|{\bf P}|)+ f(x,{\bf u},{\bf 0}) \le f(x,{\bf u},{\bf P}) \le L \phi(1+|{\bf P}|) \,;
\end{equation*}
\item {\emph{($\varphi$-growth condition II)} }\label{ass-2f}  for every $x\in \Omega$, ${\bf u}\in \R^N$ and ${\bf P}\in\mathbb{R}^{N\times n}$
\begin{equation*}
 |Df(x,{\bu},{\bP})|\le L \phi'(1+|{\bP}|) \quad \text{and}\quad |D^2f(x,{\bu},{\bP})|\le L \phi''(1+|{\bP}|);
\end{equation*}

\item\emph{($\varphi$-quasiconvexity)}  \label{ass-3f} for every $x \in \Omega$, ${\bf u} \in \R^N$ and ${\bf P}\in \mathbb{R}^{N\times n}$, and every  $\bm\eta\in C^{\infty}_0(B, \mathbb{R}^N)$ with ball $B\subset \Omega$,
\begin{equation*}
\int_{B}{f (x, {\bf u}, {\bf P}+D\bm\eta(y)) -f(x, {\bf u}, {\bf P})} \, \dd y \geq \, \nu \, \int_{B}\phi'' (1+ \kabs{{\bf P}} + \kabs{D\bm\eta(y)} )\, \kabs{D\bm\eta(y)}^2 \, \dd y\, ;
\end{equation*}

\item \emph{(VMO-condition for $x$)} \label{ass-4f} {with respect to the dependence on the $x$-variable
we do not impose a continuity condition, but we assume that} the function $x\mapsto f ( x, {\bf u}, {\bf P})/\phi(1+|{\bf P}|)$ satisfies the following VMO type condition, uniformly with respect to $({\bf u}, {\bf P})\,$: for every ${\bf u} \in \R^N$ and ${\bf P}\in \mathbb{R}^{N\times n}$, and every $\Omega_r(x_0):=\Omega\cap B_r(x_0)$ with $x_0\in\overline\Omega$ and $r\in(0,1]$,
\begin{equation*}
\kabs{ f (x, {\bf u}, {\bf P})- (f( \cdot, {\bf u}, {\bf P}))_{x_0, r}} \leq  {v}_{x_0,r}(x)\,\phi(1+|{\bf P}|)\quad \mbox{ for all $x\in B_r(x_0)\,,$}
\end{equation*}
 where  ${v}_{x_0,r}:\Omega_{r}(x_0)\to[0,2L]$ are bounded functions such that
\begin{equation*}
 \lim_{\varrho\to0}{\mathcal{V}}(\varrho)=0 \quad  \mbox{ with } \quad  {\mathcal{V}}(\varrho):=\sup_{0<r\leq\varrho} \sup_{x_0\in\Omega} \dashint_{\Omega_r(x_0)}{v}_{x_0,r}(x)\,\mathrm{d}x\,,
\end{equation*}
and
\begin{equation*}
(f( \cdot, {\bf u}, {\bf P}))_{x_0, r}:=\frac{1}{|\Omega_r(x_0)|}\int_{\Omega_r(x_0)} f(x, {\bf u}, {\bf P})\,\mathrm{d}x\,;
\end{equation*}

\item \emph{(Uniform continuity for $\bf u$)} \label{ass-5f} for every $x \in \Omega$, ${\bf u}, {\bf u}_0 \in \R^N$ and ${\bf P}\in \R^{N\times n}$
\begin{equation*}
\kabs{ f (x, {\bf u}, {\bf P})- f (x, {\bf u_0}, {\bf P})}\le L\omega ( |{\bf u}-{\bf u_0}|) \phi(1+|{\bf P}|)\,, 
\end{equation*}
where $\omega:[0,\infty)\to[0,1]$ is a non-decreasing, concave modulus of continuity; i.\,e., $\lim_{t\downarrow 0}\omega(t)=\omega(0)=0$.

\item \emph{(Continuity of $D^2 f$ for $\bf P$ away from $\bf 0$)} \label{ass-6f} there exists a non-decreasing, concave function $\omega_1:[0,\infty)\to[0,1]$ with $\lim_{t\downarrow 0}\omega_1(t)=\omega_1(0)=0$ such that for every $x \in \Omega$ and  ${\bf u}\in \R^N$, and every ${\bf P}, {\bf Q}\in \R^{N\times n}$ with $0<|{\bf P}|\le \frac12 {(|{\bf Q}|+1)}$
\begin{equation*}
\kabs{ D^2 f (x, {\bf u}, {\bf Q})- D^2 f (x, {\bf u}, {\bf P}+{\bf Q}) }  \le   L \phi'' (1+\kabs{{\bf Q}})\,
{\omega_1\left(\frac{\kabs{{\bf P}}}{1+\kabs{{\bf Q}}}\right)}\,.
\end{equation*}
Without loss of generality, we assume that the map $t\mapsto \omega_1(t)/t$ is \textit{almost decreasing} with constant $L\ge 1$ (see the next section for the definition of ``almost decreasing"). 
\end{enumerate}
Note that the $\phi$-growth and the $\phi$-quasiconvexity conditions \ref{ass-2f} and \ref{ass-3f} imply the following strong Legendre–Hadamard condition:
\begin{equation}\label{LHf}
 \tilde \nu \phi''(1+|{\bP}|) |\bfxi|^2 |\bfzeta|^2 \le \langle D^2f(x,{\bu},{\bP}) \bfxi \otimes \bfzeta \, |\,  \bfxi \otimes \bfzeta  \rangle \le  \tilde L \phi''(1+|{\bP}|) |\bfxi|^2 |\bfzeta|^2
\end{equation} 
for all $x\in\Omega$, ${\bu}, {\bfxi} \in \R^n$, ${\bfzeta}\in \R^N$ and ${\bP} \in \R^{N\times n}$, and for some  $0 < \tilde \nu\le \tilde L $ depending on $n,N,\nu,L$. In addition from the first inequality in \ref{ass-2f},
\begin{equation}
|f(x,{\bf u},{\bf P})-f(x,{\bf u},{\bf Q})|\leq L|{\bf P}-{\bf Q}|\varphi'(1+|{\bf P}|+|{\bf Q}|) \leq c \varphi(1+|{\bf P}|+|{\bf Q}|)\,.
\label{(1.3celok)}
\end{equation}

Now we state our main regularity result (for the notation, we refer to Section~\ref{sec:notation} below).  

\begin{theorem}
\label{theorem-result-1}
Let $\Omega \subset \R^n$ be a bounded $C^1$ domain, ${\bf g}\in C^1(\overline\Omega,\R^N)$ and $\phi$ be  an $N$-function satisfying assumptions {\rm\ref{ass-1phi}} -- {\rm\ref{ass-2phi}}. Consider a minimizer ${\bf u} \in {\bf g}+W^{1,\phi}_0(\Omega,\R^N)$ to the functional {\rm\eqref{functional}} under assumptions {\rm\ref{ass-1f}} -- {\rm\ref{ass-6f}}. {Then the set} of regular points on the
boundary $\partial\Omega$ given by 
\begin{equation*}
\partial\Omega_{\bu}:=\bigcap_{\alpha\in(0,1)}\{x_0\in\partial\Omega:\,\, \bfu\in C^{\alpha}(U_{x_0}\cap \overline{\Omega};\R^N) \mbox{ for some }U_{x_0}\}\,,
\end{equation*}
where $U_{x_0}$ is an open neighborhood of $x_0$, satisfies $\partial\Omega\backslash\partial\Omega_{\bu}\subset \mathrm{Sing}_{\bu}(\partial\Omega)$ where
\begin{equation*}
{\begin{split}
\mathrm{Sing}_{\bu}(\partial\Omega):=& \left\{x_0\in\partial\Omega:\,\, \mathop{\lim\inf}_{\varrho\searrow 0} \dashint_{\Omega_{\varrho}(x_0)} |D{\bf u}-(D_{\bm\nu_{x_0}}{\bf u})_{\Omega_{\varrho}(x_0)}\otimes {\bm\nu_{x_0}}|\,\mathrm{d}x>0\right\}\\
&\ \cup\left\{x_0\in\partial\Omega:\,\, \mathop{\lim\sup}_{\varrho\searrow 0}\,(|D_{\bm\nu_{x_0}}{\bf u}|)_{\Omega_{\varrho}(x_0)}=+\infty\right\}\,.
\end{split}}
\end{equation*} 
where ${\bm\nu_{x_0}}$ is the inward unit normal vector at $x_0\in\partial\Omega$.
\end{theorem}

\section{Preliminaries and basic results}\label{sec:prelbas}

\subsection{Notation} \label{sec:notation}

We denote by $B_r(x_0)$ the ball of radius $r>0$ centered at $x_0\in\R^n$. We also consider the upper half ball
\begin{equation*}
B_r^+(x_0):=\{x\in\R^n:\,\, x_n>0\,,\,\, |x-x_0|<r\},
\end{equation*}
where $x_0\in\R^{n-1}\times\{0\}$, and the part of $B_r(x_0)$ in $\Omega$
\begin{equation*}
\Omega_r(x_0):=B_r(x_0)\cap \Omega.
\end{equation*}
 We write
\begin{equation*}
\Gamma_r(x_0):= \{x=(x_1,\dots,x_n)\in\R^n:\,\, |x-x_0|<r\,,\,\,x_n=0\}
\end{equation*}
for $x_0\in\R^{n-1}\times\{0\}$. In the case $x_0=0$, we will use the shorthands $B_r$, $B^+_r$, $\Omega_r$ and $\Gamma_r$ in place of $B_r(0)$, $B^+_r(0)$, $\Omega_r(0)$ and $\Gamma_r(0)$, respectively. 

For $U\subset\R^n$ and $f\in L^1(U,\R^d)$ for some $d\in \mathbb N$, we denote the average of f by
$$
(f)_{U} := \dashint_U f \, \dd x=\frac{1}{|U|} \int_U f\, \dd x 
\quad \text{and, \ in particular,}\quad
(f)_{x_0,r}:= (f)_{B_r^+(x_0)},  
$$ 
where $x_0\in \R^{n-1}\times\{0\}$.

We say that a function $\psi:I \to \R$, where $I\subset\R$ is an interval, is \textit{almost increasing}, or \textit{almost decreasing}, with constant $L\ge 1$, if  $\psi(s)\le L\psi(t)$, or $\psi(t)\le L \psi(s)$, for all $s,t\in I$ with $s\le t$. In particular, when $L=1$, we say $f$ is non-decreasing, or non-increasing. Regarding the almost decreasing condition, we introduce the following Jensen type inequality:
\begin{lemma}\label{lem:Jensen}
Suppose $\psi:[0,\infty)\to[0,\infty)$ is a non-decreasing function such that $\psi(t)/t$ is almost decreasing in $t\in(0,\infty)$ 
{with constant $L$,}
and $f\in L^1(\Omega)$. Then 
\begin{equation*}
\dashint_\Omega \psi(|f|)\, \dd x \le (L+1) \psi \Bigg( \dashint_\Omega |f|\, \dd x\Bigg).
\end{equation*}
\end{lemma} 
\proof
By the same argument in the proof of \cite[Lemma 2.2]{Ok2}, there exists a concave function such $\widetilde \psi:[0,\infty)\to[0,\infty)$ that $\psi(t)\le \widetilde \psi(t) \le (L+1)\psi(t)$, for all $t>0$. Therefore, applying Jensen's inequality to $\widetilde \psi$, we have the desired inequality.
\endproof

\subsection{Some basic facts on $N$-functions} \label{sec:basicNfunctions}

We recall basic notation and  properties  about Orlicz functions. The following definitions and results can be found, e.g., in \cite{Kras, Kufn, Bennett, Adams}. 

A function $\phi\colon [0,\infty) \to [0,\infty)$ is called  an \emph{$N$-function} if it is convex with $\phi(0)=0$ and  $\varphi$ admits  $\phi':[0,\infty)\to[0,\infty)$ such that $\phi(t)=\int_0^t\phi'(s)\,\dd s$, where $\phi'$ is right continuous, non-decreasing and satisfies $\phi'(0) = 0$, $\phi'(t)>0$ for $t>0$, and $\lim_{t\to \infty} \phi'(t)=\infty$. From now on, $\phi$ is always an $N$-function.

We say that $\phi$ satisfies the \emph{$\Delta_2$-condition} if $\phi(2t) \leq c\,
\phi(t)$ for all $t>0$ and for some $c\ge 1$. Here, we denote the smallest possible such constant $c$ by $\Delta_2(\phi)$. Note that  the $\Delta_2$-condition is equivalent to $\phi(2t) \sim \phi(t)$ and implies $\phi(t)\sim t\,\phi'(t)$ uniformly  in $t\ge 0$.

We define by the right inverse of $\phi'$
\begin{align*}
  (\phi')^{-1}(t) &:= \sup \{ s \in [0,\infty)\,:\
    \phi'(s) \leq t \},
\qquad t\ge 0,
\end{align*}
which is well-defined by the definition of $N$-function. If $\phi'$  is strictly increasing $(\phi')^{-1}$ is the usual inverse of $\phi'$.
 We further define the \textit{Young-Fenchel-Yosida conjugate function} of
$\phi$ by
\begin{align*}
  \phi^\ast(t) &:= \int_0^t (\phi')^{-1}(s)\,ds, \qquad t\ge 0.
\end{align*}
Then $\phi^*$ is again an $N$-function with $(\phi^\ast)'(t) =(\phi')^{-1}(t)$, $(\phi^\ast)^\ast = \phi$ and
\[
\phi^*(t)= \sup_{a \geq 0} \,(at - \phi(a)).
\]
For convenience, $\Delta_2(\phi,\phi^\ast)$ stands for $\Delta_2(\phi)$ and $\Delta_2(\phi^\ast)$, and $\Delta_2(\phi,\phi^\ast)<\infty$ means $\varphi$ and $\varphi^*$ satisfy the $\Delta_2$-condition. If $\Delta_2(\phi,\phi^\ast)<\infty$, we have $\phi^*(\phi'(t))\sim \phi(t)$ uniformly in $t\ge0$ and, from the last identity above, the following Young's inequality: for any $\delta>0$ there exists $c_\delta>0$ depending on $\Delta_2(\phi,\phi^\ast)$ and $\delta$
such that
\begin{equation*}
  at   \leq \delta\, \phi(t) + c_\delta\, \phi^\ast(a) \quad \text{for all }\ t,a\ge 0.
\end{equation*}
In particular, if $\delta=1$, we can take $c_\delta=1$. 

Now, we consider an $N$-function $\phi$ satisfying {\rm\ref{ass-1phi}} -- {\rm\ref{ass-2phi}}. We recall   \cite[Proposition 2.1]{GoodSciStroff2021} that collects elementary properties of such $\phi$.

\begin{proposition}\label{prop:properties}
Let $\varphi$ be an $N$-function complying with {\rm\ref{ass-1phi}} -- {\rm\ref{ass-2phi}}. Then
\begin{description}
\item[(i)] it holds that
\begin{equation*}
(\mu_1-1)\phi'(t) \le  t\,\phi''(t) \le (\mu_2-1)\phi'(t)  \ \ \Longleftrightarrow\ \ \phi'(t)\approx t \phi''(t) \,
\end{equation*}
uniformly in $t > 0$.
 The constants 
 $\mu_1$ and $\mu_2$  are called the {\em characteristics of~$\phi$};
 

\item[(ii)] it holds that
\begin{equation*}
\mu_1\leq \inf_{t>0}\frac{t\varphi'(t)}{\varphi(t)}\leq \sup_{t>0}\frac{t\varphi'(t)}{\varphi(t)}\leq\mu_2\,;
\end{equation*}

\item[(iii)] the mappings
\begin{equation*}
t\in(0,+\infty)\to \frac{\varphi'(t)}{t^{\mu_1-1}}\,,\,\, \frac{\varphi(t)}{t^{\mu_1}} \mbox{ \,\, and \,\, } t\in(0,+\infty)\to \frac{\varphi'(t)}{t^{\mu_2-1}}\,,\,\, \frac{\varphi(t)}{t^{\mu_2}}
\end{equation*}
are non-decreasing and non-increasing, respectively;
\item[(iv)] as for the functions $\varphi$ and $\varphi'$ applied to multiples of given arguments, the following inequalities hold for every $t\geq0$:
\begin{align*}
& a^{\mu_2}\varphi(t) \leq \varphi(at) \leq  a^{\mu_1}\varphi(t) \mbox{ \,\, and \,\, } a^{\mu_2-1}\varphi'(t) \leq \varphi'(at) \leq  a^{\mu_1-1}\varphi'(t) \mbox{ \, if \, } 0<a\leq1\,; \\
& a^{\mu_1}\varphi(t) \leq \varphi(at) \leq  a^{\mu_2}\varphi(t) \mbox{ \,\, and \,\,} a^{\mu_1-1}\varphi'(t) \leq \varphi'(at) \leq  a^{\mu_2-1}\varphi'(t) \mbox{\, if \,} a\geq1\,.
\end{align*}
In particular, it follows that both $\varphi$ and $\varphi^*$  satisfy the $\Delta_2$-condition with constants $\Delta_2(\phi)$ and $\Delta_2(\phi^*)$ determined by $\mu_1$ and $\mu_2$.
\end{description}
\end{proposition}

We remark from the the preceding proposition that 
\begin{equation*}
\phi(t) \sim \phi'(t)\,t\,, \qquad \phi(t) \sim \phi''(t)\,t^2\,,\qquad \phi^\ast\big( \phi'(t) \big) \sim \phi^\ast\big( \phi(t)/t \big)\sim \phi(t)\,,
\end{equation*}
and for every $a,t>0$
\begin{equation*}
\min\{a^{\frac{1}{\mu_1}},a^{\frac{1}{\mu_2}}\}\,\varphi^{-1}(t)\leq \varphi^{-1}(at)\leq \max\{a^{\frac{1}{\mu_1}},a^{\frac{1}{\mu_2}}\}\, \varphi^{-1}(t).
\end{equation*}

%

For a given $N$-function $\varphi$, 
we define \textit{shifted $N$-functions}  $\phi_a$, $a\ge 0$, 
{introduced in \cite{DIEETT08}}
by
\begin{align}
  \label{eq:phi_shifted}
  \phi_a(t):= \int _0^t \varphi_a'(s)\, \mathrm{d}s\qquad\text{with }\quad
  \phi'_a(t):=\phi'(a+t)\frac {t}{a+t}.
\end{align}
Note that $\phi_0=\phi$. If $\phi$ satisfies {\rm\ref{ass-1phi}} -- {\rm\ref{ass-2phi}}, then    families $\{\phi_a \}_{a \ge 0}$ and
$\{(\phi_a)^* \}_{a \ge 0}$ satisfy the $\Delta_2$-condition uniformly in $a \ge 0$ and
\begin{align}
&\phi_a(t) \sim \phi'_a(t)\,t =\frac{\varphi(a+t)}{(a+t)^2}t^2\sim \frac{\varphi'(a+t)}{a+t}t^2\sim \phi''(a+t)t^2\,,\label{(2.6b)}\\
& \phi(a+t)\sim [\phi_a(t)+\phi(a)]\,,\nonumber\\ 
& {\min\{a,(\mu_1-1)\} \le \frac{ t\,\phi_a''(t) }{\phi_a'(t)}\le \max\{(\mu_2-1),a\}} \,,\nonumber 
\end{align}
see \cite{DIEETT08}. Moreover, by \cite[Corollary~26]{DieKre08} we have the following lemma, which deals with the \emph{change of shift} for $N$-functions.

\begin{lemma}\label{lem:changeshift}
Let $\varphi$ be an $N$-function with $\Delta_2(\varphi,\varphi^*)<\infty$. Then for any $\eta>0$ there exists $c_\eta>0$, depending only on $\eta$ and $\Delta_2(\varphi)$, such that for all ${\bf a}, {\bf b}\in\R^d$, $d\in\mathbb N$, and $t\geq0$
\begin{equation}
\varphi_{|{\bf a}|}(t) \leq c_\eta\varphi_{|{\bf b}|}(t) + \eta \varphi_{|{\bf a}|}(|{\bf a}-{\bf b}|)\,.
\label{(5.4diekreu)}
\end{equation}
\end{lemma}

We further define $\bfV_a\,:\, \R^{m} \to \R^{m}$, where $a\ge 0$ and $m\in \mathbb N$, by
\begin{equation*}
    \bfV_a(\bfQ):=\sqrt{\phi_a'(|\bfQ|)|\bfQ|}\frac{\bfQ}{|\bfQ|}\,,
\end{equation*}
where $\phi_a'$ is defined in \eqref{eq:phi_shifted} with $\varphi$ satisfying {\rm\ref{ass-1phi}} -- {\rm\ref{ass-2phi}}. 
We write $\bfV(\bfQ)=\bfV_0(\bfQ)= \sqrt{\phi'(|\bfQ|)|\bfQ|}\frac{\bfQ}{|\bfQ|}$.
Then we have that  
\begin{equation*}
|{\bfV_a(\bfQ)}|^2 = \phi_a'(|\bfQ|)|\bfQ|\sim \phi_a(|{\bfQ}|)
\end{equation*}
and, by Young's inequality, 
\begin{equation*}
\phi'_a(|{\bf Q}|) |{\bf P}|\, \leq \, \phi_a^*(\phi'_a(|{\bf Q}|))+ \phi_a(|{\bf P}|)\, \sim\, \phi_a(|{\bf Q}|)+ \phi_a(|{\bf P}|)\, \sim\, |{\bf V}_a({\bf Q})|^2+|{\bf V}_a({\bf P})|^2\,
\end{equation*}
uniformly in $\bfP, \bfQ \in \R^{m}$ and $a\ge 0$. Moreover, we have from {\cite[Lemma 3]{DIEETT08}}, replacing $\phi$ with $\phi_1$, that 
\begin{equation}\label{eq:equivalence}
   |{ \bfV_1(\bfP) - \bfV_1(\bfQ)}|^2 \sim  \phi_{1+|{\bfP}|}(|{\bfP - \bfQ}|)\,, \quad \bfP,\bfQ\in \R^m\,,
\end{equation}
and from \cite[Lemma~A.2]{DiKaSch}, replacing $\bfV$ with $\bfV_1$, that for $\bfg\in W^{1,\phi}(B_r(x_0);\R^m)$,
\begin{equation}
\dashint_{B_r(x_0)}|{\bf V}_1(\bfg)-{\bf V}_1((\bfg)_{B_{r}(x_0)})|^2\,\mathrm{d}x \sim  \dashint_{B_r(x_0)}|{\bf V}_1(\bfg)-({\bf V}_1(\bfg))_{B_{r}(x_0)}|^2\,\mathrm{d}x\,.
\label{eq:equivalencebis}
\end{equation}

{
\begin{remark}
In recent contributions (see \cite{DFTW20,BDGP22}) a new definition of shifted $N$-function has been devised; namely, 
\begin{equation}
\phi_a(t) := \int_0^t\frac{\phi(a\vee s)}{a \vee s} s \, \mathrm{d} s, \quad t\ge 0\,,
\label{eq:newshifted}
\end{equation}
where $s_1\vee s_2:=\max\{s_1,s_2\}$ for every $s_1,s_2\in\mathbb{R}$.
The two versions of shifted $N$-function share almost all properties. The main difference is that the new definition yields the equality $(\phi^*)_a= (\phi^*)_{\phi'(a)}$, in place of the equivalence relation provided by the original definition \eqref{eq:phi_shifted}. Then this equality implies sharper constants in some estimates, as in the ``removal of shift'' lemma, \cite[Lemma 13]{BDGP22} which improves, with the new definition \eqref{eq:newshifted}, 
the result of Lemma~\ref{lem:changeshift} with ${\bf b}={\bf 0}$. 
\end{remark}
}

\subsection{Orlicz-Sobolev spaces}\label{sec:orliczsob}

$L^\phi$ and $W^{1,\phi}$ are the classical Orlicz and
Orlicz-Sobolev spaces, i.\,e.,\ $f \in L^\phi$ iff $\int
\phi(|{f}|)\,dx < \infty$ and $f \in W^{1,\phi}$ iff $f, D f
\in L^\phi$. 
Note that equipped with the Luxemburg norm $L^\phi$ is a reflexive Banach space, hence so is $W^{1,\phi}$. 
The space $W^{1,\phi}_0(\Omega;\R^N)$ denotes  the closure of $C^\infty_0(\Omega;\R^N)$ in $W^{1,\phi}(\Omega;\R^N)$, and for ${\bf g}\in W^{1,\varphi}(\Omega,\R^N)$, we denote ${\bf g}+W^{1,\varphi}_0(\Omega,\R^N):=\{{\bf f} + {\bf g}\,:\,{\bf f}\in W^{1,\phi}_0(\Omega;\R^N)\}$.
We also introduce the following notation for $W^{1,\varphi}$-functions defined on some half-ball $B_r^+(x_0)$ with $x_0\in \R^{n-1}\times\{0\}$ and which vanish (in the sense of traces) on the flat part of the boundary:
\begin{equation*}
W^{1,\varphi}_\Gamma(B_r^+(x_0);\R^N):=\{{\bf u}\in W^{1,\varphi}(B_r^+(x_0);\R^N)\,:\, {\bf u}={\bf 0} \mbox{ on }\Gamma_r(x_0)\}\,.
\end{equation*}

We next introduce Poincar\'e and Poincar\'e-Sobolev type inequalities in the Orlicz-Sobolev space on Lipschitz boundaries.  
The first lemma is a Poincar\'e type inequality  for Sobolev functions with  zero value on the flat boundary, that can be found in \cite[Lemma~2.4]{Ok3}. In this case, compared with a usual Poincar\'e type inequality,  the gradient on the right-hand side is replaced by the directional derivative $D_n{\bf u}$.  We also note that in \cite[Lemma~2.4]{Ok3} the $N$-function $\phi$ is assumed to satisfy \ref{ass-1phi} -- \ref{ass-2phi}, but the result still holds for any $N$-function $\phi$ satisfying  $\Delta_2(\phi,\phi^*)<\infty$.  

\begin{lemma}\label{thm:sob-poincare2}
Let $\varphi:[0,+\infty)\to[0,+\infty)$ be an $N$-function complying with $\Delta_2(\phi,\phi^*)<\infty$. There exists $c>0$ depending on $n, N, \Delta_2(\phi,\phi^*)$ such that if ${\bf u}\in W^{1,1}(B_r^+(x_0);\R^N)$ with ${\bf u}={\bf 0}$ on $\Gamma_r(x_0)$, where $x_0\in \R^{n-1}\times\{0\}$,
\begin{equation}
\dashint_{B_r^+(x_0)} \varphi\left(\frac{|{\bf u}|}{r}\right)\,\mathrm{d}x\leq c\dashint_{B_r^+(x_0)} \varphi\left({|D_n{\bf u}|}\right)\,\mathrm{d}x\,.
\label{eq:2.9ok}
\end{equation}
\end{lemma}

The second lemma is Poincar\'e-Sobolev type inequalities for Sobolev {functions. 
We} can deduce these results from the standard Poincar\'e-Sobolev inequality in \cite[Lemma 7]{DIEETT08}.

\begin{theorem}\label{thm:sob-poincare}
Let $\varphi:[0,+\infty)\to[0,+\infty)$ be an $N$-function complying with $\Delta_2(\phi,\phi^*)<\infty$, There exist $\alpha\in(0,1)$ and $c>0$ depending on $n, N, \Delta_2(\phi,\phi^*)$ such that 
if ${\bf u}\in W^{1,1}(B_r;\R^N)$
\begin{equation}
\dashint_{B_r} \varphi\left(\frac{|{\bf u}- ({\bf u})_{r}|}{r}\right)\,\mathrm{d}x \leq c \Bigg(\dashint_{B_r} \varphi^{\alpha}\left({|D{\bf u}|}\right)\,\mathrm{d}x\Bigg)^{\frac{1}{\alpha}};
\label{sob-poincare-ineq}
\end{equation}
moreover, if ${\bf u}\in W^{1,1}(B_r;\R^N)$ with ${\bf u}={\bf 0}$ on $\partial B_r$, or ${\bf u}={\bf 0}$ in $A\subset B_r$ with $|A|>0$,
\begin{equation}
\dashint_{B_r} \varphi\left(\frac{|{\bf u}|}{r}\right)\,\mathrm{d}x \leq c \Bigg(\dashint_{B_r} \varphi^{\alpha}\left({|D{\bf u}|}\right)\,\mathrm{d}x\Bigg)^{\frac{1}{\alpha}},
\label{eq:2.8ok}
\end{equation}
where in the second case, the constant $c$ depends also on $\frac{|B_r|}{|A|}$.
\end{theorem}

\subsection{Some useful lemmas} \label{sec:somelemmas}

The following lemma is useful in order to re-absorb certain terms (see \cite[Lemma~6.1]{GIUSTI} and also \cite[Lemma~3.1]{DIELENSTROVER12}).

\begin{lemma}\label{lem:iterationlemma}
Let $\psi$ be an $N$-function with $\Delta_2(\psi)<\infty$. For  $0<r <\varrho$, let $h\in L^\psi(B_{\varrho})$ and $g:[r,\varrho]\to\R$ be nonnegative and bounded such that for all $r\leq s<t\leq\varrho$
\begin{equation*}
g(s)\leq\theta g(t) + A \int_{B_t}\psi\left(\frac{|h(y)|}{t-s}\right)\,\mathrm{d}y+ \frac{B}{(t-s)^\beta}+C\,,
\end{equation*}
where $A,B,C\geq0$, $\beta>0$ and $\theta\in[0,1)$. Then,  {for all $0<r<\rho$,}
\begin{equation*}
g\left(r\right)\leq c(\theta,\Delta_2(\psi),\beta)\left[A\int_{B_{\varrho}}\psi\left(\frac{|h(y)|}{\varrho-r}\right)\,\mathrm{d}y + \frac{B}{(\varrho-r)^\beta}+C\right]\,.
\end{equation*}
\end{lemma}

We now state another useful iteration lemma (see, e.g., \cite[Lemma~7.3]{GIUSTI}).

\begin{lemma}\label{lem:iterationlemma2}
Let $\psi:(0,\varrho]\to\R$ be a positive and non-decreasing function satisfying
\begin{equation*}
\psi(\vartheta^{k+1}\varrho) \leq \vartheta^\lambda \psi(\vartheta^{k}\varrho) + \tilde{c}(\vartheta^k\varrho)^n \,\, \mbox{ for every }\,\,k=0,1,2,\dots,
\end{equation*}
where $\vartheta\in(0,1)$, $\lambda\in(0,n)$ and $\tilde{c}>0$. Then there exists $c=c(n,\vartheta,\lambda)>0$ such that  
\begin{equation*}
\psi(t) \leq c\left\{\left(\frac{t}{\varrho}\right)^\lambda \psi(\varrho) + \tilde{c}t^\lambda\right\} \,\, \mbox{ for every }\,\, t\in (0,\varrho]\,.
\end{equation*}
\end{lemma}

%
%
%

\subsection{$\mathcal{A}$-harmonic approximation on half balls}\label{subsection:Aharmonic}

\begin{figure}[h]
\centering
\def\svgwidth{200pt}
\begingroup%
  \makeatletter%
  \providecommand\color[2][]{%
    \errmessage{(Inkscape) Color is used for the text in Inkscape, but the package 'color.sty' is not loaded}%
    \renewcommand\color[2][]{}%
  }%
  \providecommand\transparent[1]{%
    \errmessage{(Inkscape) Transparency is used (non-zero) for the text in Inkscape, but the package 'transparent.sty' is not loaded}%
    \renewcommand\transparent[1]{}%
  }%
  \providecommand\rotatebox[2]{#2}%
  \newcommand*\fsize{\dimexpr\f@size pt\relax}%
  \newcommand*\lineheight[1]{\fontsize{\fsize}{#1\fsize}\selectfont}%
  \ifx\svgwidth\undefined%
    \setlength{\unitlength}{552.27272334bp}%
    \ifx\svgscale\undefined%
      \relax%
    \else%
      \setlength{\unitlength}{\unitlength * \real{\svgscale}}%
    \fi%
  \else%
    \setlength{\unitlength}{\svgwidth}%
  \fi%
  \global\let\svgwidth\undefined%
  \global\let\svgscale\undefined%
  \makeatother%
  \begin{picture}(1,0.57716052)%
    \lineheight{1}%
    \setlength\tabcolsep{0pt}%
    \put(0,0){\includegraphics[width=\unitlength,page=1]{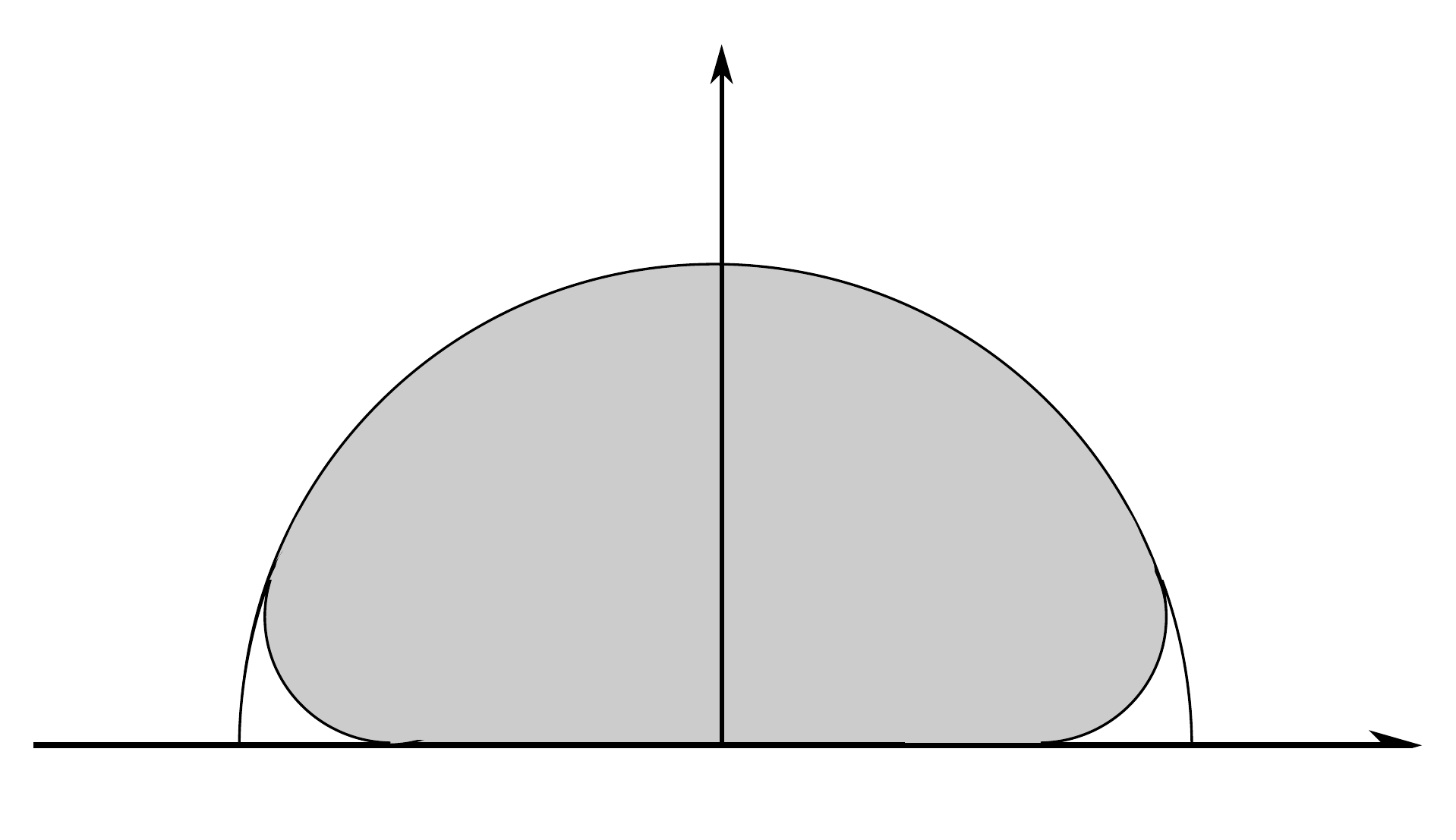}}%
    \put(0.53549382,0.51697534){\makebox(0,0)[lt]{\lineheight{1.25}\smash{\begin{tabular}[t]{l}$x_n$\end{tabular}}}}%
    \put(0.89814815,0.08796302){\makebox(0,0)[lt]{\lineheight{1.25}\smash{\begin{tabular}[t]{l}$x'$\end{tabular}}}}%
    \put(0.33950619,0.22376546){\makebox(0,0)[lt]{\lineheight{1.25}\smash{\begin{tabular}[t]{l}$U$\end{tabular}}}}%
    \put(0.80785409,0.0272314){\makebox(0,0)[lt]{\lineheight{1.25}\smash{\begin{tabular}[t]{l}$\varrho$\end{tabular}}}}%
    \put(0,0){\includegraphics[width=\unitlength,page=2]{setU.pdf}}%
  \end{picture}%
\endgroup%

\caption{The smooth neighborhood $U$.}\label{fig:setU}
\end{figure}

We introduce here a flat boundary version of the $\mathcal{A}$-harmonic approximation result in the setting of Orlicz spaces proven in \cite[Theorem~14]{DIELENSTROVER12}. 
{As remarked in the Introduction, the lack of global Calder\'on-Zygmund estimates in the Orlicz spaces $L^\phi$ for $\mathcal A$-harmonic maps in half balls leads us to adopt a trick to prove an $\mathcal A$-harmonic approximation lemma on half balls.
 The idea is to construct a suitable ``smooth neighborhood'' of the upper half ball; i.e., a domain $U$ of class $C^1$ 
containing  
the half ball in such a way that the flat part of the boundary be a proper subset of $\partial U$ (see Fig.~\ref{fig:setU}),
and then to obtain global Calder\'on-Zygmund estimates in this domain (see Theorem~\ref{thm:thm18}). 
}

Let $\mathcal{A}$ be a bilinear form on $\R^{N\times n}$. We say that $\mathcal{A}$ is {\em strongly elliptic in the sense of Legendre-Hadamard} if for all $\bm\xi\in \R^N,\bm\zeta\in\R^{n}$ it holds that 
\begin{equation}
  \kappaA \abs{\bm\xi}^2 \abs{\bm\zeta}^2\leq \langle\mathcal{A}(\bm\xi\otimes\bm\zeta)|(\bm\xi\otimes\bm\zeta)\rangle\leq L_{\mathcal{A}} \abs{\bm\xi}^2 \abs{\bm\zeta}^2
\label{(2.20)}
\end{equation}
for some $0< \kappaA \le L_{\mathcal{A}}$, and $\phi$ be an $N$-function with $\Delta_2(\phi,\phi^*)<\infty$. We say that a function $\bfw\in W^{1,\phi}(\Omega;\R^N)$  is
\emph{$\mathcal{A}$-harmonic} on $\Omega$ if it satisfies $-\divergence (\mathcal{A}D \bfw)=0$ in the sense of distributions; i.\,e.,
\begin{equation*}
\int_{\Omega} \langle\mathcal{A}D{\bf w}|D\bm\eta\rangle\,\mathrm{d}x=0\,,\quad \mbox{ for all }\ \bm\eta\in C^\infty_0(\Omega;\R^N)\,.
\end{equation*}

Let us denote by $\partial U$ the set of all $x=(x',x_n)\in \R^{n-1}\times \R=\R^n$ such that $0\le x_n\le \varrho$ and $x'$ satisfies
\begin{equation*}
\begin{cases}
|x'|^2+x_n^2 =\varrho^2  &\text{if } \ \frac{\varrho}{3}\le x_n \le \varrho,\\
(|x'|-\frac{\sqrt 2\varrho}{2})^2+(x_n-\frac{\varrho}{4})^2 =\frac{\varrho^2}{16} & \text{if } \ 0<x_n < \frac{\varrho}{3},\\
|x'| \le \frac{\sqrt2 \varrho}{2} & \text{if } \ x_n=0 
\end{cases}
\end{equation*}
(see Fig.~\ref{fig:setU}), and $U$ the interior region of $\partial U$.  Note that $\partial U$ consists of portions of spheres with radii $\varrho$ and $\varrho/4$, and is of class $C^1$ with $C^1$ semi-norm depending only on $n$ . 
Moreover, we also have $B_{\sqrt 2 \varrho/2}^+\subset U \subset B_{\varrho}^+$, hence $\partial U$ can be considered as a ``smooth neighborhood'' of $B_{\sqrt 2 \varrho/2}^+$ with the flat part of the boundary in common with that of $\partial B_{\sqrt 2 \varrho/2}^+$.

We then consider the Dirichlet problem
\begin{equation}
  \label{eq:dirich1}
  \begin{cases}
    -\divergence (\mathcal{A} D {\bf w}) = -\divergence {\bf F} &\qquad\text{in $U$},
    \\
    {\bf w}= \bfzero &\qquad\text{on $\partial U$},
  \end{cases}
\end{equation}
and have the following result of existence and uniqueness of weak solutions in Sobolev Orlicz spaces (see \cite[Theorem~18]{DIELENSTROVER12}).


\begin{theorem}
Let $\mathcal{A}$ be strongly elliptic in the sense of Legendre-Hadamard, $\phi$ be an $N$-function with $\Delta_2(\phi,\phi^*)<\infty$, and $U\subset\R^n$ as above. Then for every ${\bf F}\in L^\varphi(U;\R^{N\times n})$, the $\mathcal A$-harmonic system \eqref{eq:dirich1} admits a unique weak solution ${\bf w}\in W^{1,\varphi}_0(U;\R^N)$ such that
\begin{equation*}
\|D{\bf w}\|_{L^\varphi(U)}\leq c \|{\bf F}\|_{L^\varphi(U)}
\end{equation*}
and
\begin{equation}
\int_{U}\varphi(|D{\bf w}|)\,\mathrm{d}x \leq c \int_{U}\varphi(|{\bf F}|)\,\mathrm{d}x\,,
\label{eq:caldezyg}
\end{equation}
where $c$ depends on $n, N, \kappa_{\mathcal A}, L_{\mathcal{A}}$ and $\Delta_2(\varphi,\varphi^*)$.
\label{thm:thm18}
\end{theorem}
Note that the constant $c$ in the above theorem is independent of the structure and the size of $\partial U$, since $\partial U$ is constructed by connecting in a smooth way a flat region with only portions of spheres whose radii are $\varrho$ and $\varrho/4$.


We next introduce regularity estimates for the weak solution to 
\begin{equation}
  \label{eq:dirich2}
  \begin{cases}
    -\divergence (\mathcal{A} D {\bf h}) = 0 &\qquad\text{on $B^+_{r}$}
    \\
    {\bf h}= \bfzero &\qquad\text{on $\Gamma_{r}$.}
  \end{cases}
\end{equation}
We refer to \cite[Section 10.2]{GIUSTI}, in particular, Eq. (10.22) and Remark 10.1. 
{Note that the condition $\Delta_2(\phi,\phi^*)<\infty$ implies that $W^{1,\phi}_\Gamma(B^+_{r};\R^N)\subset W^{1,p}_\Gamma(B^+_{r};\R^N)$ for some $p>1$. Hence, in what follows we can take advantage of the classical regularity theory results.} 
\begin{theorem}\label{thm:Aharmonic}
Let $\mathcal{A}$ be strongly elliptic in the sense of Legendre-Hadamard, and $\phi$ be an $N$-function with $\Delta_2(\phi,\phi^*)<\infty$.  If  ${\bf h}\in W^{1,\phi}_\Gamma(B^+_{r};\R^N)$ is a weak solution to \eqref{eq:dirich2}, then ${\bf h} \in W^{2,\infty}(B_{r/2}^{+};\R^N)$ and we have
\begin{equation}\label{(2.21)}
\sup_{B_{r/2}^+} \left(|D {\bf h}|+r |D^2 {\bf h}|\right)\leq c \dashint_{B^+_r} |D_n{\bf h} |\,dx.
\end{equation}
where 
{$c=c(n,N,\kappa_\mathcal{A},L_{\mathcal{A}})$.} 
\end{theorem}

{The following is a flat boundary version of the $\mathcal A$-harmonic approximation in the setting of Orlicz spaces in \cite{DIELENSTROVER12}.}

\begin{theorem}
  \label{thm:Aappr_phi}
Let $\mathcal{A}$ be strongly elliptic in the sense of Legendre-Hadamard, $\phi$ be an $N$-function with $\Delta_2(\phi,\phi^*)<\infty$, and let $s>1$. Then for every
  $\epsilon>0$, there exists $\delta>0$  depending only on $n, N, \kappa_A, L_\mathcal{A},  \Delta_2(\varphi,\varphi^*)$ and
  $s$ such that the following holds.  Let $\bfv \in
  W^{1,\varphi}_\Gamma(B^+_\varrho,\R^N)$ be {\em almost $\mathcal{A}$-harmonic} in the sense that

  \begin{align}
    \label{eq:Aappr_ah}
    \biggabs{\dashint_{B^+_\varrho} \langle\mathcal{A}D \bfv | D \bm\eta\rangle\,\mathrm{d}x}
    \leq \delta \dashint_{B^+_\varrho} \abs{D \bfv}\,\mathrm{d}x
    \norm{D \bm\eta}_{L^\infty(B^+_\varrho)}
    \quad   \quad{for\,\, all }\ {\bm\eta} \in C^\infty_0(B^+_\varrho;\R^N).
  \end{align}
 Then there exists ${\bf z}
  \in W^{1, \varphi}_\Gamma(B^+_{\varrho/2};\RN)$ such that ${\bf z}$ is a weak solution to 
  \begin{equation}
  \begin{cases}
    -\divergence (\mathcal{A} D {\bf z}) =     -\divergence (\mathcal{A} D \bfv) &\qquad\text{on $B^+_{\frac{\sqrt 2 \varrho}{2}}$}
    \\
    {\bf z}= \bfzero &\qquad\text{on $\Gamma_{\frac{\sqrt 2\varrho}{2}}$,}
  \end{cases}
  \label{zsystem}
\end{equation}
and
   \begin{align*}
    \dashint_{B^+_{\varrho/2}} \varphi\bigg(\frac{\abs{{\bf z}}}{\varrho}\bigg)\,\mathrm{d}x +
    \dashint_{B^+_{\varrho/2}} \varphi(\abs{D {\bf z}})\,\mathrm{d}x \leq \epsilon
    \Bigg(\bigg(\dashint_{B^+_{\varrho}} \varphi^s(\abs{D
      \bfv}) \,\mathrm{d}x\bigg)^{\frac 1s}+\dashint_{B^+_{\varrho}} \varphi(\abs{D \bfv})\,\mathrm{d}x\Bigg)\,.
  \end{align*}
\end{theorem}

\proof
Let $U$ be the set constructed in the beginning of this subsection.  The argument of \cite[Theorem~14]{DIELENSTROVER12} holds true (with minor changes) for $U$ in place of $B$. We then prefer to omit the details, just providing a sketch of the proof.

Our aim is to find the $\mathcal{A}$-harmonic function ${\bf h}\in W^{1,\varphi}(U;\RN)$ with the same boundary values as ${\bf u}$; i.\,e., the solution of the problem
\begin{equation*}
  \begin{cases}
    -\divergence (\mathcal{A} D {\bf h}) = 0
     &\qquad\text{in $U$}
    \\
    {\bf h}= \bfu &\qquad\text{on $\partial U$}\,.
  \end{cases}
\end{equation*}
Setting ${\bf z}:={\bf h}-\bfv$, the previous problem can be rewritten as 
\begin{equation}
  \label{eq:calA2}
  \begin{cases}
    -\divergence (\mathcal{A} D {\bf z}) = -\divergence(\mathcal{A}
    D \bfv) &\qquad\text{in $U$}
    \\
    {\bf z}= \bfzero &\qquad\text{on $\partial U$}
  \end{cases}
\end{equation}
in the sense of distributions. Now, as a consequence of Theorem~\ref{thm:thm18}, there exists a unique solution ${\bf z}$ to the equation \eqref{eq:calA2} in $W^{1,\varphi}_0(U;\R^N)$.  
{Note that  ${\bf z}$ is also a weak solution to \eqref{zsystem}  and 
complies with the estimate
\begin{equation*}
\int_{B^+_{\frac{\sqrt 2}{2}\varrho}}\phi(|D{\bf z}|)\,\mathrm{d}x \leq \int_{U}\phi(|D{\bf z}|)\,\mathrm{d}x \leq c \int_{U}\phi(|D{\bf v}|)\,\mathrm{d}x  \leq c \int_{B^+_{\varrho}}\phi(|D{\bf v }|)\,\mathrm{d}x\,.
\label{eq:caldezyg}
\end{equation*}}
From this point on, the proof follows exactly (with the obvious minor changes) along the lines of \cite[Theorem~14]{DIELENSTROVER12}. Indeed, \cite[Lemma~20]{DIELENSTROVER12} can be applied to ${\bf z}$ in order to obtain the necessary variational inequalities involving test functions, then one can exploit the result of Lipschitz truncation in Orlicz spaces (\cite[Theorem~21]{DIELENSTROVER12}), which holds for open sets $\Omega$ with Lipschitz boundary whose Lipschitz semi-norms depend only on $n$, and in particular for $\Omega=U$.
\endproof

\begin{remark}\label{rem:thmmodified}
We will exploit the previous approximation result in a slightly modified version. Indeed, as already remarked in \cite[Lemma~2.7]{CeladaOk}, under the additional assumption
\begin{equation*}
\dashint_{B^+_\varrho} \varphi(|D{\bfv}|)\,\mathrm{d}x \leq \left(\dashint_{B^+_\varrho} \varphi^s(|D{\bfv}|)\,\mathrm{d}x\right)^{\frac{1}{s}}\leq \varphi(\mu)
\end{equation*}
for some exponent $s>1$ and for a constant $\mu>0$, and \eqref{eq:Aappr_ah} replaced by
  \begin{equation*}
        \biggabs{\dashint_{B^+_\varrho} \langle\mathcal{A}D \bfv | D \bm\eta\rangle\,\mathrm{d}x}
    \leq \delta \mu
    \norm{D \bm\eta}_{L^\infty(B^+_\varrho)}\,,
  \end{equation*}
it can be also seen with minor changes in the proof that the function ${\bf z}$ above satisfies
\begin{equation*}
\dashint_{B^+_{\varrho/2}} \varphi\bigg(\frac{\abs{{\bf z}}}{\varrho}\bigg)\,dx +
    \dashint_{B^+_{\varrho/2}} \varphi(\abs{D {\bf z}})\,dx \leq \epsilon \varphi(\mu)\,.
\end{equation*}
\end{remark}

%
%

\section{Caccioppoli type inequalities and higher integrability}\label{sec:caccioppoliandhigh}

\subsection{Caccioppoli type inequality I}
We derive Caccioppoli type inequalities near the boundary for the minimizer of the functional \eqref{functional},
where $f:\Omega\to \R^{N} \times \R^{N\times n}\to\R$ satisfies the following assumptions:
\begin{equation}\label{wass-f}
\nu \varphi(|{\bf P}|)   \leq f(x, \bfu, {\bf P})-f(x, \bfu, {\bf 0}) \leq L \varphi(1+|{\bf P}|) 
\ \mbox{ and }\ \ 
|f(x, \bfu, {\bf P})- f(x, \bfu_0, {\bf P})|\le L\phi(1+|{\bf P}|)
\end{equation}
for all ${\bf P}\in\R^{n\times N}$ and ${\bf u},{\bf u_0}\in\R^{N}$ and for some $0<\nu\le L$, and $\phi$ is an $N$-function with $\Delta_2(\phi,\phi^*)<\infty$. Then, using these inequalities, we obtain higher integrability results.  We note that the assumptions  \ref{ass-1f}, \ref{ass-2f} and \rm\ref{ass-5f} imply \eqref{wass-f}, see \eqref{(1.3celok)}.

 We first introduce a Caccioppoli inequality on a boundary region for the minimizer of \eqref{functional}. This is a boundary, and also {non-degenerate}, counterpart of \cite[Lemma~3.1]{GoodSciStroff2021}. 

{\begin{lemma}\label{lem:caccioppoli1}
Suppose that $\Omega$ is a  Lipschitz domain, 
$B_\varrho(x_0)\not\subset\Omega$ with $x_0\in\overline\Omega$,
$f:\Omega\times \R^{N} \times \R^{N\times n}\to\R$ satisfies \eqref{wass-f} with an $N$-function $\phi$ complying with $\Delta_2(\phi,\phi^*)<\infty$, and ${\bf g}\in W^{1,\varphi}(\Omega_\varrho(x_0);\R^N)$. Let $\bfu\in W^{1,\varphi}(\Omega_\varrho(x_0);\R^N)$  be a minimizer of \eqref{functional}, replacing $\Omega$ with $\Omega_\varrho(x_0)$, such that $\bfu - {\bf g}= {\bf 0}$ on $\partial\Omega \cap B_r(x_0)$. Then, for every $r\in[\varrho/2,\varrho)$ there holds
\begin{equation}
\int_{\Omega_r(x_0)} \varphi(1+|D{\bf u}|)\,\mathrm{d}x\leq c \int_{\Omega_\varrho(x_0)}\varphi\left(1+\frac{|\bfu-\bfg|}{\varrho-r}\right)\,\mathrm{d}x + c \int_{\Omega_\varrho(x_0)}\varphi (|D\bfg|)\,\dd x 
\label{eq:caccioppoli1}
\end{equation}
for some constant $c=c(L,\nu,\Delta_2(\phi,\phi^*))>0$.
\end{lemma}}
{\proof
We assume that $x_0=0$ for simplicity. Let $r<s<t<\rho$, $\eta\in  C^{\infty}_0(B_{t})$ with $\eta\ge0$, $\eta\equiv 1$ on $B_{s}$, and $|D\eta|\le c(n)/(t-s)$, and set $\bfxi:= \bfu - \eta (\bfu-\bfg)$. Since $\eta (\bfu-\bfg) \in W^{1,\phi}_0(\Omega_{t})$, by the minimality of $\bfu$ and \eqref{wass-f},
\[\begin{aligned}
\int_{\Omega_t} \varphi(|D\bfu|)\, \dd x  
& \le \frac{1}{\nu} \int_{\Omega_t}[ f(x,\bfu, D\bfu)-f(x,\bfu, {\bf 0})] \, \dd x
\le \frac{1}{\nu}\int_{\Omega_t} [f(x,\bfxi, D\bfxi)- f(x,\bfu, {\bf 0})] \, \dd x \\
& \le \frac{1}{\nu} \int_{\Omega_t} [c \phi(1+|D\bfxi|) + f(x,\bfxi, {\bf0})- f(x,\bfu, {\bf 0})] \, \dd x \\
& \le c_0 \int_{\Omega_t\setminus \Omega_s}  \phi (|D\bfu|) \, \dd x + c  \int_{\Omega_t}  \left[\phi \left(\frac{|\bfu-\bfg|}{t-s}\right)+\phi(|D\bfg|) + \phi(1)\right] \, \dd x .
\end{aligned}\]
This implies 
$$
\int_{\Omega_s} \varphi(1+|D\bfu|)\, \dd x \le   \frac{c_0}{1+c_0} \int_{\Omega_t}  \phi (1+ |D\bfu|) \, \dd x + c  \int_{\Omega_t}  \phi \left({\frac{|\bfu-\bfg|}{t-s}}\right)  \, \dd x + c\int_{\Omega_\varrho} [ \phi (|D\bfg|)+ \phi(1)] \, \dd x .
$$
Therefore, applying Lemma~\ref{lem:iterationlemma}, we have the desired estimate. \endproof}

Using the preceding lemma, we prove the following higher integrability result up to the boundary on a Lipschitz domain satisfying the following \textit{exterior uniform measure density condition}: there exist $c_\Omega>0$ and $r_0>0$ such that 
\begin{equation}\label{measuredensity}
|B_r| \leq c_\Omega |\Omega^{\mathrm c}\cap B_r(x_0)|\quad \text{for all $x_0\in \partial \Omega$ and all $r\in(0,r_0]$.}
\end{equation}
 
\begin{theorem}\label{thm:higLip}
Suppose that $\Omega$ is a  Lipschitz domain satisfying 
\eqref{measuredensity} with  constants $c_\Omega,r_0>0$, 
$B_\varrho(x_0)\not\subset\Omega$ with $x_0\in\overline\Omega$ and $0\le \rho\le r_0/2$,
$f:\Omega\times \R^{N} \times \R^{N\times n}\to\R$ satisfies \eqref{wass-f} with an $N$-function $\phi$ complying with $\Delta_2(\phi,\phi^*)<\infty$,
and ${\bf g}\in W^{1,\varphi}(\Omega_\varrho(x_0);\R^N)$ satisfies $\phi(|D\bfg|)\in L^{s_1}(\Omega_\varrho(x_0))$ for some $s_1>1$. 
There exist an exponent $s_0\in (1,s_1)$ and a constant $c>0$ depending only on $n,N,L,\nu,\Delta_2(\phi,\phi^*),s_1,c_\Omega$ such that  if $\bfu\in W^{1,\varphi}(\Omega_\varrho(x_0);\R^N)$ is a minimizer of \eqref{functional}, replacing $\Omega$ with $\Omega_\varrho(x_0)$, and  $\bfu - {\bf g}= {\bf 0}$ on $\partial\Omega \cap B_{\varrho}(x_0)$, then, for every $B_{8r}(y)\subset B_\varrho(x_0)$ with $y\in \overline\Omega$, 
\begin{equation}
\dashint_{\Omega_{r}(y)} \varphi^{s_0}(1+|D{\bf u}|)\,\mathrm{d}x\leq c \Bigg(\dashint_{\Omega_{8r}(y)} \varphi(1+|D{\bf u}|)\,\mathrm{d}x\Bigg)^{s_0}+ c\dashint_{\Omega_{8r}(y)} \varphi^{s_0}(|D{\bfg}|)\,\mathrm{d}x\,.
\label{highint00}
\end{equation}
\end{theorem}

{\proof
 Fix any $B_{8r}(y)\subset B_{\varrho}(x_0)$. If $B_{2r}(y)\subset \Omega$, in view of  
\cite[Theorem 2.7 and Lemma 3.1]{GoodSciStroff2021} we infer that
$$
\dashint_{B_{r}(y)} \varphi(1+|D{\bf u}|)\,\mathrm{d}x
\leq c\dashint_{B_{2r}(y)} \varphi\left(1+\frac{|{\bf u}- ({\bf u})_{B_{2r}(y)} |}{r}\right)\,\mathrm{d}x
\leq c\Bigg(\dashint_{B_{2r}(y)} \varphi^\alpha({1+|D{\bf u}|})\,\mathrm{d}x\Bigg)^\frac{1}{\alpha}
$$
for some $\alpha=\alpha(n,N,\Delta_2(\phi,\phi^*))\in(0,1)$ and $c=c(n,N,L,\nu,\Delta_2(\phi,\phi^*))>0$. On the other hand, if $B_{2r}(y)\cap \partial\Omega\neq \emptyset$, applying Lemma~\ref{lem:caccioppoli1} with $(\tilde y, 6r, 3r)$ in place of  $(x_0,\varrho,r)$ and \eqref{eq:2.8ok} with   $\bfu-\bfg$ in place of $\bfu$, we infer that
$$
\begin{aligned}
\int_{\Omega_{r}(y)} \varphi(1+|D{\bf u}|)\,\dd x  
&\le  \int_{\Omega_{3r}(\tilde y)} \varphi(1+|D{\bf u}|)\,\dd x \\
&\le c \int_{\Omega_{6r}(\tilde y)} \varphi\left(1+\frac{|\bfu-\bfg|}{r}\right)\,\dd x + c \int_{\Omega_{6r}(\tilde y)} \varphi(|D\bfg|)\,\dd x\\
&\leq  c|B_r|(1+c_\Omega)\Bigg(\dashint_{B_{6r}(\tilde y)} \varphi^\alpha(|D(\bfu-\bfg)|)\,\mathrm{d}x\Bigg)^\frac{1}{\alpha} + c\int_{\Omega_{6r}(\tilde y)} \varphi(1+|D\bfg|)\,\dd x
\end{aligned}
\label{eq:stimagehring}
$$
Here we extended $\bfu-\bfg$ to $B_r(x_0)\setminus\Omega$ by ${\bf 0}$ and used the fact that $B_{r}(y)\subset B_{3r}(\tilde y)\subset B_{6r}(\tilde y)\subset B_{8r}(y)$. 
Therefore, we obtain that for every $B_{8r}(y)\subset B_\varrho(x_0)$, 
$$
\dashint_{B_{r}(y)} \varphi(1+|D{\bf u}|)\mathbbm{1}_\Omega\,\dd x  \leq \tilde c\Bigg(\dashint_{B_{8r}(y)} [\varphi\left({1+|D{\bf u}|}\right)\mathbbm{1}_\Omega]^{\alpha}\,\mathrm{d}x\Bigg)^\frac{1}{\alpha} + \tilde c\dashint_{B_{8r}(y)} \varphi(|D\bfg|)\mathbbm{1}_\Omega\,\dd x\,
$$
for some $\tilde c=\tilde c(n,N,L,\nu,\Delta_2(\phi,\phi^*),c_\Omega)>0$, where $\mathbbm{1}_\Omega$ is the characteristic function of  $\Omega$.
Then, by virtue of Gehring's Lemma (see e.g. \cite[Theorem~6.6]{GIUSTI}) we prove \eqref{highint00}.
\endproof

From the preceding theorem with a covering argument, we derive the following two corollaries that will be used later. 
The first corollary is a higher integrability result on the flat boundary. Its interior version in the degenerate case can be found in  \cite[Lemma 3.2]{GoodSciStroff2021}.

{\begin{corollary}\label{cor:higflat}
Suppose $x_0\in \R^{n-1}\times \{0\}$, and $f: B_{\varrho}^+(x_0)\times \R^{N} \times \R^{N\times n}\to\R$ satisfies \eqref{wass-f} with an $N$-function $\phi$ complying with $\Delta_2(\phi,\phi^*)<\infty$. 
Let $\bfu\in W^{1,\varphi}_\Gamma(B_{\varrho}^+(x_0);\R^N)$ be a minimizer of \eqref{functional} with $\Omega=B_{\varrho}^+(x_0)$. 
There exist an exponent $s_0>1$ and a constant $c>0$ depending only on $n,N,L,\nu,\Delta_2(\phi,\phi^*)$ such that
 for every $s\in(1,s_0]$ and every radius $r\in[\varrho/2,\varrho)$, one has 
\begin{equation}
\dashint_{B_r^+(x_0)} \varphi^{s}(1+|D{\bf u}|)\,\mathrm{d}x \leq {c}\left(\frac{\varrho}{\varrho-r}\right)^{n(s-1)}\Bigg(\dashint_{B_\varrho^+(x_0)}\varphi(1+|D{\bf u}|)\,\mathrm{d}x\Bigg)^{s}\,.
\label{highint0}
\end{equation}
{Moreover, we also have for any $t\in(0,1]$
\begin{equation}
\dashint_{B_{\varrho/2}^+(x_0)} \varphi(1+|D_n{\bf u}|)\,\mathrm{d}x \leq c_t\Bigg(\dashint_{B_\varrho^+(x_0)}\varphi(1+|D_n{\bf u}|)^t\,\mathrm{d}x\Bigg)^{\frac{1}{t}}\,,
\label{highint1}
\end{equation}
where $c_t>0$ depends on $n,N,L,\nu,\Delta_2(\phi,\phi^*)$ and $t$.}
\end{corollary}}

\proof
Set $\Omega=\R^n_+$. Note that  $\R^n_+$ satisfies the uniform measure density condition \eqref{measuredensity} with $c_\Omega=2$ and any $r_0>0$. By Theorem~\ref{thm:higLip} with $\bfg\equiv \bf0$,  we infer that 
$$
\dashint_{B^+_{\varrho-r}(y)} \varphi^{s}(1+|D{\bf u}|)\,\mathrm{d}x\leq c \Bigg(\dashint_{B^+_{8(\varrho-r)}(y)} \varphi(1+|D{\bf u}|)\,\mathrm{d}x\Bigg)^{s}\,
$$
for every $r\in[\varrho/2,\varrho)$  and every $y\in \overline{B_r^+(x_0)}$. Therefore, by a standard covering argument, we prove 
{\eqref{highint0}. 

Moreover, from \eqref{highint0}, \eqref{eq:caccioppoli1} and \eqref{eq:2.9ok} we deduce that for every $\frac12\le \tau_1 < \tau_2\le 1$

$$
\Bigg(\dashint_{B_{\tau_1\varrho}^+(x_0)} \varphi(|D_n{\bf u}|)^{s_0}\,\mathrm{d}x \Bigg)^{\frac{1}{s_0}} \leq \frac{c}{(\tau_2-\tau_1)^{\beta_0}} \dashint_{B_\varrho^+(x_0)}\varphi(1+|D_n{\bf u}|)\,\mathrm{d}x\,.
$$
for some $\beta_0=\beta_0(n,s_0)>0$. For $t\in(0,1)$, set $\theta\in(0,1)$ such that $t=\frac{s_0(1-\theta)}{s_0-\theta}$. Then an application of H\"older's inequality and Young's inequality yields
$$\begin{aligned}
&\Bigg(\dashint_{B_{\tau_1\varrho}^+(x_0)} \varphi(|D_n{\bf u}|)^{s_0}\,\mathrm{d}x \Bigg)^{\frac{1}{s_0}}  \leq \frac{c}{(\tau_2-\tau_1)^{\beta_0}}\dashint_{B_\varrho^+(x_0)}\varphi(1+|D_n{\bf u}|)\,\mathrm{d}x\\
&\leq \frac{c}{(\tau_2-\tau_1)^{\beta_0}}\Bigg(\dashint_{B_{\tau_2\varrho}^+(x_0)}\varphi(1+|D_n{\bf u}|)^{s_0}\,\mathrm{d}x\Bigg)^{\frac{\theta}{s_0}}\Bigg(\dashint_{B_{\tau_2\varrho}^+(x_0)}\varphi(1+|D_n{\bf u}|)^t\,\mathrm{d}x\Bigg)^{\frac{1-\theta}{t}}\\
&\leq \frac{1}{2}\Bigg(\dashint_{B_{\tau_2\varrho}^+(x_0)}\varphi(1+|D_n{\bf u}|)^{s_0}\,\mathrm{d}x\Bigg)^{\frac{1}{s_0}}+ \frac{c}{(\tau_2-\tau_1)^{\frac{\beta_0}{1-\theta}}}\Bigg(\dashint_{B_{2\varrho}^+(x_0)}\varphi(1+|D_n{\bf u}|)^t\,\mathrm{d}x\Bigg)^{\frac{1}{t}}.
\end{aligned}$$ 
Finally, applying Lemma~\ref{lem:iterationlemma}, we have \eqref{highint0}.}
\endproof

\endproof}

The second corollary is a global higher integrability result in an upper-half ball.
 
\begin{corollary}\label{cor:highalfball}
Suppose that $x_0\in \R^{n-1}\times \{0\}$, $f: B_{\varrho}^+(x_0)\times \R^{N} \times \R^{N\times n}\to\R$ satisfies \eqref{wass-f} with an $N$-function $\phi$ complying with $\Delta_2(\phi,\phi^*)<\infty$,
 and ${\bfg}\in W^{1,\varphi}(B^+_\varrho(x_0),\R^N)$ satisfies $\varphi(|D{\bfg}|)\in L^{s_1}(B^+_\varrho(x_0))$ for some $s_1>1$. 
 There exists an exponent $s_0\in(1,s_1)$ and a constant $c>0$ depending only on $n,N,L,\nu,\Delta_2(\phi,\phi^*),s_1$ such that if $\bfu\in W^{1,\varphi}(B^+_\varrho(x_0);\R^N)$ is the minimizer of \eqref{functional} with $\Omega=B^+_\varrho(x_0)$, and  $\bfu =  \bfg $ on $\partial B^+_\varrho(x_0)$, 
 then $\phi(|D\bfu|)\in L^{s_0}(B_{\varrho}^+(x_0))$ with the following estimate: for every $s\in (1,s_0]$
\begin{equation}\label{globalhigestimate}
\int_{B^+_{\varrho}(x_0)} \varphi^{s}(1+|D{\bfu}|)\,\mathrm{d}x\leq  c \int_{B^+_\varrho(x_0)} \varphi^{s}(1+|D{\bf g}|)\,\mathrm{d}x\,.
\end{equation}
 
\end{corollary}

\proof
First, we observe that $B_\varrho^+(x_0)$ satisfies the exterior measure density condition in \eqref{measuredensity} with $c_\Omega=2$ and $r_0=\varrho$; hence, we obtain \eqref{highint00} for all $B_{8r}(y)\subset B_\varrho(x_0)$ with $r\le \varrho/16$ and $y\in \overline{B_\varrho(x_0)}$. Therefore, 
{applying a standard covering argument and the fact that, by \eqref{wass-f} and the minimality of $\bu$,
$$\begin{aligned}
\int_{B_\varrho^+ (x_0)} \varphi(|D\bu|)\, \dd x  &\le \frac{1}{\nu}\int_{B_\varrho^+ (x_0)} f(x,\bu , D\bu) - f(x,\bu , {\bf 0})  \, \dd x \\
&\le \frac{1}{\nu}\int_{B_\varrho^+ (x_0)} f(x,\bg , D\bg) - f(x,\bg , {\bf 0})+ f(x,\bg , {\bf 0}) - f(x,\bu , {\bf 0})  \, \dd x  \\
&\le c \int_{B_\varrho^+ (x_0)} \phi(1+|D\bg|)  \, \dd x  \,,
\end{aligned}$$
 we obtain \eqref{globalhigestimate}.}
\endproof

\subsection{Caccioppoli type inequality II}

Next, we  derive the following Caccioppoli inequality of second type for  the minimizers of \eqref{functional}, where $f:\Omega\times \R^N\times \R^{N\times n}$ satisfies \ref{ass-1f}  --  \ref{ass-5f} with an $N$-function complying with \ref{ass-1phi}  --  \ref{ass-2phi}. We note that $f$ satisfies \eqref{wass-f}, and $\Delta_2(\phi)$ and $\Delta_2(\phi^*)$ are finite and depend on $\mu_1$ and $\mu_2$.

{\begin{lemma}\label{lem:3.5}
Suppose that $x_0\in \R^{n-1}\times \{x_0\}$,  $\varphi$ satisfies {\rm\ref{ass-1phi} -- \ref{ass-2phi}}, and $f:B_\varrho^+(x_0)\times \R^N\times \R^{N\times n}$ does  {\rm\ref{ass-1f} -- \ref{ass-5f}}. Let $\bfu\in W^{1,\varphi}_\Gamma(B_\varrho^+(x_0))$ be a minimizer of \eqref{functional}. Then, for every $\bm \xi\in \R^N$ and every $s\in(1,s_0]$,
\begin{equation}
\begin{split}
\dashint_{B_{\varrho/2}^+(x_0)} \varphi_{1+|\bm\xi|}(|D\bfu - \bm\xi\otimes {\bf e}_n|)\,\mathrm{d}x  &
\leq c \dashint_{B_{\varrho}^+(x_0)} \varphi_{1+|\bm\xi|}\left(\frac{|\bfu - x_n\bm\xi|}{\varrho}\right)\,\mathrm{d}x \\
&\hspace{-2cm} + c \varphi(1+|\bm\xi|)\Bigg[\omega\Bigg(\dashint_{B_{\varrho}^+(x_0)}|{\bf u}|+|{\bf u}-\bm\xi x_n|\,\mathrm{d}x\Bigg)^{1-\frac{1}{s}}+ {\mathcal{V}}(\varrho)^{1-\frac{1}{s}}\Bigg]
\end{split}
\label{eq:3.2ok}
\end{equation}
for some $c=c(n,N,\mu_1,\mu_2,\nu,L)>0$, where $s_0>1$ is the constant determined in Lemma~\ref{cor:higflat} with $\phi$ and $f$ given in the present lemma; hence it depends only on $n,N,\mu_1,\mu_2,\nu,L$, and $\mathcal V$ and $\omega$ are given in {\rm\ref{ass-4f}} and {\rm\ref{ass-5f}}, respectively.
\end{lemma}}

\proof
We may adapt the interior argument of \cite[Lemma~3.4]{GoodSciStroff2021} to the boundary setting, with minor modifications due to the different assumptions \ref{ass-1f} -- \ref{ass-2f} therein. Assuming, without loss of generality, that $x_0=0$, we fix any $\frac{\varrho}{2}\leq r<\tau<t\leq \frac{3\varrho}{4}$ with $\tau:=\frac{r+t}{2}$, and a cut-off function $\eta\in C_0^\infty(B_\tau;[0,1])$ such that $\eta\equiv1$ on $B_r$ and $|D\eta|\leq \frac{4}{t-r}$ on $B_\tau$. Correspondingly, we define the test functions $\bm\zeta:=\eta({\bf u} - x_n\bm\xi)\in W^{1,\varphi}_\Gamma(B_\tau^+;\R^N)$ and $\bm\psi:=(1-\eta)({\bf u} - x_n\bm\xi)\in W^{1,\varphi}_\Gamma(B_\tau^+;\R^N)$. From this point, the proof goes on as in \cite[Lemma~3.4]{GoodSciStroff2021}. We then omit the details. 
\endproof

For a minimizer $\bfu\in W^{1,\phi}_\Gamma(B_\varrho^+(x_0);\R^N)$ of \eqref{functional} with $\Omega=B_\varrho^+(x_0)$ and $x_0\in \R^{n-1}\times\{0\}$, we set
\begin{equation}
\Theta(x_0,\varrho):=\varrho\varphi^{-1}\Bigg(\dashint_{B_{\varrho}^+(x_0)}\varphi(1+|D_n{\bf u}|)\,\mathrm{d}x\Bigg)\,.
\label{eq:Theta}
\end{equation}
Then, we obtain some upper bounds for the remainder term inside $\omega$ in \eqref{eq:3.2ok}.

{\begin{lemma}\label{lem:remainder}
Let $\phi$ be an $N$-function complying with \ref{ass-1phi} -- \ref{ass-2phi}, and $\bfu\in W^{1,\phi}_\Gamma(B_\varrho^+(x_0);\R^N)$ with $x_0\in \R^{n-1}\times\{0\}$. Then, for every $\bm\xi\in\R^N$, if the smallness assumption
\begin{equation}
\dashint_{B_{\varrho}^+(x_0)} \varphi_{1+|{\bm \xi}|}(|D{\bf u}-\bm\xi\otimes {\bf e}_n|)\,\mathrm{d}x \leq \Lambda\varphi(1+|\bm\xi|)
\label{smallnesscon}
\end{equation}
holds for some 
{$\Lambda>0$,}
 then
\begin{equation}
\Bigg( \dashint_{B^+_{\varrho}(x_0)}|{\bf u}|^{\mu_1}\,\mathrm{d}x \Bigg)^{\frac{1}{\mu_1}}
\leq c\Theta(x_0,\varrho)
\leq c {(1+\Lambda)^{\frac{1}{\mu_2}}}\varrho (1+ |\bm\xi|)\,
\label{restostima2} 
\end{equation}
for some  
$c=c(n,N,\mu_1,\mu_2)>0$.

\end{lemma}}

\proof
First, we note from \ref{ass-2phi} that $\phi(t^{1/\mu_1})$ is convex for $t\ge 0$. Applying Jensen's inequality, the Poincar\'e type estimate in Lemma~\ref{thm:sob-poincare2} and the change-shift formula \eqref{(5.4diekreu)} with ${\bf a}={\bf 0}$, and using  assumption \eqref{smallnesscon}, we obtain
\begin{equation*}
\begin{split}
\varphi\Bigg(\Bigg(\dashint_{B_\varrho^+(x_0)} \left[\frac{|{\bf u}|}{\varrho}\right]^{\mu_1}\,\mathrm{d}x\Bigg)^{\frac{1}{\mu_1}}\Bigg)
&\leq \dashint_{B_\varrho^+(x_0)} \varphi\left(\frac{|{\bf u}|}{\varrho}\right)\,\mathrm{d}x  \leq c\dashint_{B_\varrho^+(x_0)} \varphi(|D_n{\bf u}|)\,\mathrm{d}x \\
& \leq c \dashint_{B_{\varrho}^+(x_0)} \varphi(|D{\bf u}-\bm\xi\otimes {\bf e}_n|)\,\mathrm{d}x + c \varphi(|\bm\xi|) \\
& \leq c\dashint_{B_{\varrho}^+(x_0)} \varphi_{1+|{\bm\xi}|}(|D{\bf u}-\bm\xi\otimes {\bf e}_n|)\,\mathrm{d}x + c\varphi(1+|{\bm\xi}|) \\
& \leq c(1+\Lambda)\varphi(1+|{\bm\xi}|)\,,
\end{split}
\end{equation*}
which yields \eqref{restostima2} up to applying $\varphi^{-1}$ to both the sides and Proposition~\ref{prop:properties}(iv).
\endproof

Next, we establish an higher integrability result for $\varphi_{1+|\bm\xi|}(|D{\bf u}-\bm\xi\otimes {\bf e}_n|)$ under assumption \eqref{smallnesscon}. 

\begin{lemma}\label{corollary3.2}
Suppose that $x_0\in \R^{n-1}\times \{x_0\}$,  $\varphi$ satisfies {\rm\ref{ass-1phi} -- \ref{ass-2phi}}, and $f:B_\varrho^+(x_0)\times \R^N\times \R^{N\times n}$ does  {\rm\ref{ass-1f} -- \ref{ass-5f}}. There exist an exponent $\sigma>1$ and a  constant $c>0$ depending on $n,N,\mu_1,\mu_2,\nu,L$ such that 
if ${\bf u}\in W^{1,\varphi}_\Gamma(B_r^+(x_0);\R^N)$ is a minimizer of \eqref{functional} and the inequality in \eqref{smallnesscon} holds for some for some ${\bm\xi}\in\R^{N}$ and 
{$\Lambda>0$,} 
then,  for every $s\in(1,s_0]$,
\begin{equation}
\begin{split}
&\Bigg(\dashint_{B_{\varrho/2}^+(x_0)}  \varphi_{1+|\bm\xi|}^\sigma(|D{\bf u}-\bm\xi\otimes {\bf e}_n|)\,\mathrm{d}x\Bigg)^{\frac{1}{\sigma}} \\
&\leq c \dashint_{B_{\varrho}^+(x_0)}  \varphi_{1+|\bm\xi|}(|D{\bf u}-\bm\xi\otimes {\bf e}_n|)\,\mathrm{d}x+ c\varphi(1+|{\bm \xi}|)\left[{(1+\Lambda)^{\frac{1}{\mu_2}}}\omega\left(\varrho(1+|{\bm \xi}|)\right)^{1-\frac{1}{s}}+ {\mathcal{V}}(\varrho)^{1-\frac{1}{s}}\right]\,,
\end{split}
\label{eq:caccioppoliIbis}
\end{equation} 
where $s_0>1$ is the exponent determined in Lemma~\ref{cor:higflat} with $\phi$ and $f$ given in the present lemma; hence it depends only on $n,N,\mu_1,\mu_2,\nu,L$, and $\mathcal V$ and $\omega$ are given in {\rm\ref{ass-4f}} and {\rm\ref{ass-5f}}, respectively.

{
Moreover, we also have that for every $t\in(0,1]$
\begin{equation}
\begin{split}
&\Bigg(\dashint_{B_{\varrho/2}^+(x_0)}  \varphi_{1+|\bm\xi|}^\sigma(|D{\bf u}-\bm\xi\otimes {\bf e}_n|)\,\mathrm{d}x\Bigg)^{\frac{1}{\sigma}} \leq c_t \Bigg(\dashint_{B_{\varrho}^+(x_0)}  \varphi_{1+|\bm\xi|}(|D{\bf u}-\bm\xi\otimes {\bf e}_n|)^t\,\mathrm{d}x\Bigg)^{\frac{1}{t}}
\\
&\quad + c\varphi(1+|{\bm \xi}|)\left[(1+\Lambda)^{\frac{1}{\mu_2}}\omega\left(\varrho(1+|{\bm \xi}|)\right)^{1-\frac{1}{s}}+ {\mathcal{V}}(\varrho)^{1-\frac{1}{s}}\right]\,,
\end{split}
\label{eq:caccioppoliIbis1}
\end{equation}
for some $c_t>0$ depending on $n,N,\mu_1,\mu_2,\nu,L$ and $t$.}
\end{lemma}
{\proof
We follow the argument in the proof of Theorem~\ref{thm:higLip}. We extend $\bfu$ and $\bfu- \bm\xi x_n$ to $B_{\varrho}(x_0)\setminus B^+_{\varrho}(x_0)$ by ${\bf 0}$, and  consider any ball $B_{8r}(y)\subset B_\varrho(x_0)$. If $B_{2r}(y)\subset B_{\varrho}^+(x_0)$,
in view of \cite[Lemma 3.4]{GoodSciStroff2021} with $\varrho=2r$, $x_0=y$, $\bfu_0=(\bfu)_{y,2r}$, and ${\bf Q}= {\bm\xi}\otimes {\bf e}_n$, we obtain
\begin{equation}\label{hig2pf1}
\begin{aligned}
&\dashint_{B_{r}(y)} \varphi_{1+|\bm\xi|}(|D{\bf u}-\bm\xi\otimes {\bf e}_n|)\,\mathrm{d}x 
\leq c \dashint_{B_{2r}(y)} \varphi_{1+|\bm\xi|}\left(\frac{|\bfu-(\bfu)_{y,2r}-(\bm\xi\otimes {\bf e}_n)(x-y)|}{2r}\right)\,\mathrm{d}x\\
&\qquad + c\varphi(1+|{\bm \xi}|)\Bigg[\omega\Bigg(\dashint_{B_{2r}(y)}|\bfu-(\bfu)_{y,2r}|+ |{\bm\xi}||x-y|\, \dd x\Bigg)^{1-\frac{1}{s}}+ {\mathcal{V}}(2r)^{1-\frac{1}{s}}\Bigg]\,.
\end{aligned}
\end{equation}
Here, we observe that
$$
\dashint_{B_{2r}(y)}|\bfu-(\bfu)_{y,2r}|+ |{\bm\xi}||x-y|\, \dd x \le c \dashint_{B_{2r}(y)}|\bfu|\, \dd x + c |{\bm\xi}|r \le c \Bigg(\frac{1}{\varrho(1+|{\bm\xi}|)}\dashint_{B_{2r}(y)}|\bfu|\, \dd x+1 \Bigg) \varrho (1+|\bm\xi|)
$$
which with the concavity of $\omega$ and $\omega(0)=0$, precisely, $\omega(ct)\le  c\omega(t)$ when $c\ge 1$, yields
$$
\omega\Bigg(\dashint_{B_{2r}(y)}|\bfu-(\bfu)_{y,2r}|+ |{\bm\xi}||x-y|\, \dd x\Bigg)^{1-\frac{1}{s}} 
\le c \Bigg(\frac{1}{\varrho(1+|{\bm\xi}|)}\dashint_{B_{2r}(y)}|\bfu|\, \dd x+1 \Bigg) \omega(\varrho(1+|\bm\xi|))^{1-\frac{1}{s}}.
$$
Moreover, as $\dashint_{B_{2r}(y)}\bfu-(\bfu)_{y,2r}-(\bm\xi\otimes {\bf e}_n)(x-y)\,\dd x={\bf 0}$, by the Sobolev--Poincar\'e type inequality in \eqref{sob-poincare-ineq},
$$
 \dashint_{B_{2r}(y)} \varphi_{1+|\bm\xi|}\left(\frac{|\bfu-(\bfu)_{y,2r}-(\bm\xi\otimes {\bf e}_n)(x-y)|}{2r}\right)\,\mathrm{d}x \leq c  \Bigg(\dashint_{B_{2r}(y)} \varphi_{1+|\bm\xi|}^{\alpha}(|D\bfu-\bm\xi\otimes {\bf e}_n|)\,\mathrm{d}x\Bigg)^{\frac{1}{\alpha}}  
$$
for some $\alpha\in(0,1)$. Therefore, inserting the preceding two estimates into \eqref{hig2pf1}, we obtain
$$\begin{aligned}
&\dashint_{B_{r}(y)} \varphi_{1+|\bm\xi|}(|D{\bf u}-\bm\xi\otimes {\bf e}_n|)\,\mathrm{d}x 
\leq c   \Bigg(\dashint_{B_{2r}(y)} \varphi_{1+|\bm\xi|}^{\alpha}(|D\bfu-\bm\xi\otimes {\bf e}_n|)\,\mathrm{d}x\Bigg)^{\frac{1}{\alpha}}\\
&\qquad + c\frac{ \varphi(1+|{\bm \xi}|)\omega(\varrho(1+|\bm\xi|))^{1-\frac{1}{s}}}{\varrho(1+|{\bm\xi}|)}\dashint_{B_{2r}(y)}|\bfu|\, \dd x+c\varphi(1+|{\bm \xi}|)\left[\omega\left(\varrho(1+|{\bm \xi}|)\right)^{1-\frac{1}{s}}+ {\mathcal{V}}(\varrho)^{1-\frac{1}{s}}\right]\,.
\end{aligned}$$

On the other hand, if $B_{2r}(y)\cap (\R^{n-1}\times\{0\})\neq \emptyset$, by Lemma~\ref{lem:3.5}, the Sobolev--Poincar\'e type inequality in \eqref{eq:2.8ok} and the same argument in the case that $B_{2r}(y)\subset B_{\varrho}^+(x_0)$ in above, we obtain
$$
\begin{aligned}
&\dashint_{B_{r}(y)}  \varphi_{1+|\bm\xi|}(|D\bfu - \bm\xi\otimes {\bf e}_n|)\,\mathrm{d}x 
\le c \dashint_{B_{3r}^+(\tilde y)} \varphi_{1+|\bm\xi|}(|D\bfu - \bm\xi\otimes {\bf e}_n|)\,\mathrm{d}x \\
& \leq c \dashint_{B_{6r}^+(\tilde y)} \varphi_{1+|\bm\xi|}\left(\frac{|\bfu - x_n\bm\xi|}{6r}\right)\,\mathrm{d}x  + c \varphi(1+|\bm\xi|)\Bigg[\omega\Bigg(\dashint_{B_{6r}^+(\tilde y)}|{\bf u}|+|\bm\xi x_n|\,\mathrm{d}x\Bigg)^{1-\frac{1}{s}}+ {\mathcal{V}}(6r)^{1-\frac{1}{s}}\Bigg]\\
&\leq c   \Bigg(\dashint_{B_{8r}(y)} \varphi_{1+|\bm\xi|}^{\alpha}(|D\bfu-\bm\xi\otimes {\bf e}_n|)\,\mathrm{d}x\Bigg)^{\frac{1}{\alpha}}\\
&\qquad + c\frac{ \varphi(1+|{\bm \xi}|)\omega(\varrho(1+|\bm\xi|))^{1-\frac{1}{s}}}{\varrho(1+|{\bm\xi}|)}\dashint_{B_{8r}(y)}|\bfu|\, \dd x+c\varphi(1+|{\bm \xi}|)\left[\omega\left(\varrho(1+|{\bm \xi}|)\right)^{1-\frac{1}{s}}+ {\mathcal{V}}(\varrho)^{1-\frac{1}{s}}\right]\,.
\end{aligned}
$$
Here, $\tilde y\in \R^{n-1}\times \{0\}$ is chosen to satisfy $B_{r}(y)\subset B_{3r}(\tilde y)\subset B_{6r}(\tilde y)\subset B_{8r}(y)$. 

Consequently, since $\bfu \in L^{\mu_1}(B_{\varrho}(x_0))$, by virtue of Gehring's lemma there exists $\sigma=\sigma(n,N,\mu_1,\mu_2,\nu,L)\in (1,\mu_1)$ such that for every  $\frac{\varrho}{2}\le r_1 < r_2\le \varrho$,
$$
\begin{aligned}
&\Bigg( \dashint_{B^+_{r_1}(x_0)}  \varphi_{1+|\bm\xi|}^{\sigma}(|D\bfu - \bm\xi\otimes {\bf e}_n|)\,\mathrm{d}x \Bigg)^{\frac{1}{\sigma}} \\
&\leq c  \left(\frac{r_2}{r_2-r_1}\right)^{\frac{n(\sigma-1)}{\sigma}} \dashint_{B^+_{r_2}(x_0)} \varphi_{1+|\bm\xi|}(|D\bfu-\bm\xi\otimes {\bf e}_n|)\,\mathrm{d}x \\
&\qquad + c\frac{ \varphi(1+|{\bm \xi}|)\omega(\varrho(1+|\bm\xi|))^{1-\frac{1}{s}}}{\varrho(1+|{\bm\xi}|)}\Bigg(\dashint_{B_{r_2}^+(x_0)}|\bfu|^{\sigma}\, \dd x\Bigg)^{\frac{1}{\sigma}}\\
&\qquad +c\left(\frac{r_2}{r_2-r_1}\right)^{\frac{n(\sigma-1)}{\sigma}}\varphi(1+|{\bm \xi}|)\left[\omega\left(\varrho(1+|{\bm \xi}|)\right)^{1-\frac{1}{s}}+ {\mathcal{V}}(\varrho)^{1-\frac{1}{s}}\right]\\
&\leq c  \left(\frac{r_2}{r_2-r_1}\right)^{\frac{n(\sigma-1)}{\sigma}} \dashint_{B^+_{r_2}(x_0)} \varphi_{1+|\bm\xi|}(|D\bfu-\bm\xi\otimes {\bf e}_n|)\,\mathrm{d}x  + c (1+\Lambda)^{\frac{1}{\mu_2}} \varphi(1+|{\bm \xi}|)\omega(\varrho(1+|\bm\xi|))^{1-\frac{1}{s}}\\
&\qquad +c\left(\frac{r_2}{r_2-r_1}\right)^{\frac{n(\sigma-1)}{\sigma}}\varphi(1+|{\bm \xi}|)\left[\omega\left(\varrho(1+|{\bm \xi}|)\right)^{1-\frac{1}{s}}+ {\mathcal{V}}(\varrho)^{1-\frac{1}{s}}\right]\,,
\end{aligned}
$$
where we apply Lemma~\ref{lem:remainder} to the last inequality.
Therefore, when $r_2=\varrho$ and $r_1=\varrho$,  we obtain \eqref{eq:caccioppoliIbis}. On the other hand, the estimate in  \eqref{eq:caccioppoliIbis} can be obtained by using the same argument as in the proof of \eqref{highint1}.
\endproof}

We conclude this section by introducing the excess functional and other tools useful in the sequel. 
For a minimizer $\bfu \in W^{1,\phi}_{\Gamma}(B_\varrho^+(x_0))$ of \eqref{functional} with $\Omega=B_\varrho^+(x_0)$ and $x_0\in \R^{n-1}\times\{0\}$, we define the \emph{excess functional} as
\begin{equation}
\Phi(x_0,\varrho):=\dashint_{B_\varrho^+(x_0)}\varphi_{1+|(D_n{\bf u})_{x_0,\varrho}|}(|D{\bf u}-(D_n{\bf u})_{x_0,\varrho}\otimes {\bf e}_n|)\,\mathrm{d}x
\label{(3.5)}
\end{equation}
and
\begin{equation}
\Psi(x_0,\varrho):=\dashint_{B_\varrho^+(x_0)}\varphi\left(1+\frac{|{\bf u}|}{\varrho}\right)\,\mathrm{d}x\,.
\label{(3.6)}
\end{equation}
Moreover, we define also 
\begin{equation}
H(x_0,\varrho):=\frac{1}{(1+\Lambda)^{\frac{1}{\mu_2}}+(2L)^{1-\frac{1}{s}}}\left((1+\Lambda)^{\frac{1}{\mu_2}}\omega(\varrho(1+|(D_n{\bf u})_{x_0,\varrho}|))^{1-\frac{1}{s}} + {\mathcal{V}}(\varrho)^{1-\frac{1}{s}}\right)\,,
\label{(3.8)}
\end{equation}
and
\begin{equation}
\widetilde{H}(x_0,\varrho):=\frac{1}{(1+\Lambda)^{\frac{1}{\mu_2}}+(2L)^{1-\frac{1}{s}}}\left((1+\Lambda)^{\frac{1}{\mu_2}}\omega(\Theta(x_0,\varrho))^{1-\frac{1}{s}} + {\mathcal{V}}(\varrho)^{1-\frac{1}{s}}\right)\,,
\label{(3.8bis)}
\end{equation}
{where $s=s(n,N,\mu_1,\mu_2,\nu,L)\in(1,s_0)$ is determined in the proof of Lemma~\ref{lem:lemma3.7} below (see the line above from \eqref{(3.18)})} and $\Theta(x_0,\varrho)$ is the ``Morrey-type'' excess defined in \eqref{eq:Theta}. 
Since $\omega\leq 1$ and ${\mathcal{V}}(\varrho)\leq 2L$, we have that $H(x_0,\varrho)\leq1$, $\widetilde{H}(x_0,\varrho)\leq1$. 

Under the 
assumption $\Phi(x_0,\varrho)\leq \Lambda \varphi(1+|(D_n{\bf u})_{x_0,\varrho}|)$ in \eqref{smallnesscon} 
{with $\Lambda\in(0,1]$}, by virtue of Lemma~\ref{lem:remainder}, 
\begin{equation*}
\widetilde{H}(x_0,\varrho)\leq c H(x_0,\varrho)\,.
\end{equation*}
We can rewrite the Caccioppoli inequality \eqref{eq:caccioppoliIbis} as
\begin{equation}
\Phi(x_0,\varrho/2) \leq c \Phi(x_0,\varrho) + c\varphi(1+|(D_n{\bf u})_{x_0,\varrho}|)H(x_0,\varrho)\,.
\label{eq:caccioppoliI}
\end{equation} 

In the case $x_0=0$, we will use the shorthands $\Phi(\varrho)$, $\Psi(\varrho)$, $\Theta(\varrho)$, $H(\varrho)$ and $\widetilde{H}(\varrho)$ in place of $\Phi(0,\varrho)$, $\Psi(0,\varrho)$, $\Theta(0,\varrho)$, $H(0,\varrho)$ and $\widetilde{H}(0,\varrho)$, respectively.

\section{Boundary partial regularity}\label{sec:partialreg}

In this section, we provide a linearization strategy for the minimization problem along the lines of \cite[Sections~3.2 and 3.3]{GoodSciStroff2021}, where an analogous analysis has been performed in the interior and degenerate setting. 

Throughout this section, we always assume that $f:B_{\varrho}^+(x_0)\times \R^N\times \R^{N\times n}\to \R$ satisfies  \ref{ass-1f} -- \ref{ass-6f} with an $N$-function $\phi$ complying with \ref{ass-1phi} -- \ref{ass-2phi}, and that ${\bf u}\in W^{1,\phi}_\Gamma(B_\varrho^+(x_0);\R^N)$ is a minimizer of \eqref{functional} with $\Omega=B_\varrho^+(x_0)$ and $x_0\in\R^{n-1}\times \{0\}$.

\subsection{Approximate $\mathcal{A}$-harmonicity}\label{sec:approxharmon}

 In this subsection we show that $\bfu - (D_n\bfu)_{x_0,\rho}x_n$ is  an almost $\mathcal{A}$-harmonic function for a suitable elliptic bilinear form $\mathcal{A}$, see Lemma~\ref{lem:lemma3.9} below.
To prove this result, we employ  a suitable comparison function that is obtained by a freezing argument in the variables $(x,{\bf u})$ based on Ekeland's variational principle. We recall below a version of this classical tool, whose proof can be found, e.g., in \cite[Theorem~5.6]{GIUSTI}. 

\begin{lemma}[Ekeland's principle]\label{lem:ekeland}
Let $(X,d)$ be a complete metric space, and assume that $F : X\to[0,\infty]$ be not identically $\infty$ and lower semicontinuous with respect to the metric topology on $X$. If for some $u\in X$ and some $\kappa>0$, there
holds
\begin{equation*}
F(u)\leq \inf_XF + \kappa\,,
\end{equation*}
then there exists $v\in X$ with the properties 
\begin{equation*}
d(u,v)\leq 1 \mbox{\,\, and \,\,} F(v)\leq F(w)+\kappa d(v,w) \quad \forall w\in X\,.
\end{equation*}
\end{lemma}

We start with setting 
\begin{equation*}
g(\bm\xi) \equiv g_{{x}_0,\varrho}(\bm\xi) := (f(\cdot, {\bf 0},\bm\xi\otimes {\bf e}_n))_{{x}_0,\varrho} \quad \mbox{ for all $\ \bm\xi\in\R^{N}$}
\end{equation*}
and
\begin{equation*}
K(x_0,\varrho):=\widetilde{H}(x_0,\varrho)\Psi(x_0,\varrho)\,,
\end{equation*}
where $\widetilde{H}(x_0,\varrho)$ and $\Psi(x_0,\varrho)$ are defined as in \eqref{(3.8bis)} and \eqref{(3.6)}, respectively. Then, we introduce the following metric $d$ in $X:={\bf u}+W_0^{1,\phi}(B_{\varrho/2}^+(x_0);\R^N)$:
\begin{equation*}
d({\bf w}_1,{\bf w}_2):= \frac{1}{c_*\varphi^{-1}(K(x_0,\varrho))}\dashint_{B_{\varrho/2}^+(x_0)}|D{\bf w}_1-D{\bf w}_2|\,\mathrm{d}x\,, \quad {\bf w}_1,{\bf w}_2\in{\bf u}+W_0^{1,\phi}(B_{\varrho/2}^+(x_0);\R^N)\,,
\end{equation*}
where $c_*>0$ is determined in Lemma~\ref{lem:lemma3.7} below. 
Then $(X,d)$ is a complete metric space and the following functional is lower semicontinuous in this metric topology:
\begin{equation}
\mathcal{G}[{\bf w}]:= \dashint_{B_{\varrho/2}^+(x_0)}g(D{\bf w})\,\mathrm{d}x, \qquad  \bfw\in {\bf u}+W^{1,\phi}_0(B_{\varrho/2}^+(x_0);\R^N)\,,
\label{eq:gfunctional}
\end{equation}
 
In the next lemma, we find a suitable comparison map ${\bf v}\in X= {\bf u}+W^{1,\phi}_0(B_{\varrho/2}^+(x_0);\R^N)$ by applying Ekeland's variational principle.

\begin{lemma}\label{lem:lemma3.7}
Under the setting above, there exists a minimizer ${\bf v}\in {\bf u}+W^{1,\phi}_0(B_{\varrho/2}^+(x_0);\R^N)$ of the functional 

\begin{equation*}
\widetilde{\mathcal{G}}[{\bf w}]:= \dashint_{B_{\varrho/2}^+(x_0)} g(D{\bf w})\,\mathrm{d}x + \frac{K(x_0,\varrho)}{\varphi^{-1}(K(x_0,\varrho))}\dashint_{B_{\varrho/2}^+(x_0)}|D{\bf v}-D{\bf w}|\,\mathrm{d}x\,, \quad \bfw\in {\bf u}+W^{1,\phi}_0(B_{\varrho/2}^+(x_0);\R^N)\,,
\end{equation*}
that satisfies 
\begin{equation*}
\dashint_{B_{\varrho/2}^+(x_0)}|D{\bf v}-D{\bf u}|\,\mathrm{d}x \leq c_*\varphi^{-1}(K(x_0,\varrho))
\end{equation*}
for some constant $c_*=c_*(n,N,\mu_1,\mu_2,\nu,L)>0$. Moreover, ${\bf v}$ fulfills the following Euler-Lagrange variational inequality:
\begin{equation*}
\left|\dashint_{B_{\varrho/2}^+(x_0)}\langle Dg(D{\bf v})|D{\bm\eta}\rangle\,\mathrm{d}x\right|\leq \frac{K(x_0,\varrho)}{\varphi^{-1}(K(x_0,\varrho))}\dashint_{B_{\varrho/2}^+(x_0)}|D\bm\eta|\,\mathrm{d}x
\end{equation*}
for every $\bm\eta\in C^{\infty}_0(B_{\varrho/2}^+(x_0);\R^N)$.
\end{lemma}

\proof
The proof of this result can be obtained as in \cite[Lemma~3.8]{GoodSciStroff2021}, where its interior counterpart has been proved, by also exploiting the higher integrability results in Corollary~\ref{cor:higflat} and Corollary~\ref{cor:highalfball} with $(\bu,\bg,\varrho)$ in place of $(\tilde{\bf v}, \bu, \varrho/2)$ to infer that for some  $s\in(1,s_0)$ depending on $n,N,\mu_1,\mu_2,\nu,L$,
\begin{equation}
\Bigg(\dashint_{B_{\varrho/2}^+} \varphi^{s}(1+|D\tilde{\bf v}|)\,\mathrm{d}x\Bigg)^{\frac{1}{s}} 
\leq c \Bigg(\dashint_{B_{\varrho/2}^+} \varphi^{s}(1+|D\bfu|)\,\mathrm{d}x\Bigg)^{\frac{1}{s}}
\leq c \dashint_{B_{3\varrho/4}^+}\varphi(1+|D{\bf u}|)\,\mathrm{d}x\leq c\Psi(\varrho)\,,
\label{(3.18)}
\end{equation}
where $\tilde{\bf v}\in X$ is a minimizer of the functional \eqref{eq:gfunctional}. We then omit the details.
\endproof

\begin{corollary}\label{cor3.4celadaok}
Let  ${\bf v}\in {\bf u}+W^{1,\phi}_0(B_{\varrho/2}^+(x_0);\R^N)$ be as in Lemma~\ref{lem:lemma3.7}. Then there exists $\tau=\tau(n,N,\mu_1,\mu_2,\nu,L)\in(0,1)$ such that
\begin{equation}
\dashint_{B_{\varrho/4}^+(x_0)}\varphi(|D{\bf u}-D{\bf v}|)\,\mathrm{d}x \leq c [\widetilde{H}(x_0,\varrho)]^{1-\tau}\Psi(x_0,\varrho)
\label{(3.13celok)}
\end{equation}
for some constant $c=c(n,N,\mu_1,\mu_2,\nu,L)>0$.
\end{corollary}

\proof
We follow the proof of \cite[Corollary~3.4]{CeladaOk}. First, we derive the following higher integrability result for $\tilde{\mathcal G}$:
\begin{equation}
\Bigg(\dashint_{B^+_{\varrho/4}(x_0)}\varphi^s(1+|D{\bf v}|)\,\mathrm{d}x\Bigg)^{\frac{1}{s}} \leq c \dashint_{B^+_{\varrho/2}(x_0)}\varphi(1+|D{\bf v}|)\,\mathrm{d}x + c K(x_0,\varrho)\,,
\label{cor:omparison-pf1}
\end{equation}
where $s>1$ and $c>0$ depends on $n,N,\mu_1,\mu_2,\nu,L$.
The proof is exactly the same as the one of Theorem~\ref{thm:higLip} by obtaining the Caccioppoli type estimates for the functional $\tilde{\mathcal G}$ in interior and boundary regions. For the Caccioppoli type estimate in interior regions,  we refer to \cite[Step 1 in the proof of Corollary~3.4]{CeladaOk}, and, thus, obtain that for every $B_{2r}(y)\subset B^+_{\varrho}(x_0)$,
\begin{equation*}
\begin{aligned}
\dashint_{B_{r}(y)}\varphi(1+|D{\bf v}|)\,\mathrm{d}x 
& \leq c \dashint_{B_{2r}(y)}\varphi\left(1+\left|\frac{\bfv-(\bfv)_{y,2r}}{r}\right|\right)\,\mathrm{d}x + c K(x_0,\varrho)\\
& \leq c \Bigg(\dashint_{B_{2r}(y)}\varphi\left(1+|D\bfv |\right)^{\alpha}\,\mathrm{d}x\Bigg)^{\frac{1}{\alpha}} + c K(x_0,\varrho)\,.
\end{aligned}
\end{equation*}
Here, we used the Sobolev--Poincar\'e type inequality in \eqref{sob-poincare-ineq}.
In a similar way (see, e.g., Lemma~\ref{lem:caccioppoli1}), when $B_{2r}(y)\subset B_{\varrho}(x_0)$ with $y\in \Gamma_r(x_0)$, we also obtain the above estimate, extending $\bfv$ to $B_{2r}(y)\setminus B^+_{2r}(y)$ by ${\bf0}$, replacing $(\bfv)_{y,2r}$ with ${\bf 0}$ and  using the Sobolev--Poincar\'e type inequality in \eqref{eq:2.8ok}. Therefore, by Gehring's lemma, we obtain \eqref{cor:omparison-pf1}.

Then, by following \cite[Step 1 in the proof of Corollary~3.4]{CeladaOk}, we obtain \eqref{(3.13celok)}, where $\tau\in(0,1)$ is explicitly defined such that
\begin{equation*}
\frac{\tau}{s}+(1-\tau)\mu_2=1\,,
\end{equation*}
where $\mu_2$ denotes the bigger of the characteristics of $\varphi$.
\endproof

We are now in the position to prove the approximate $\mathcal{A}$-harmonicity of $\bfu$, where
\begin{equation}
\mathcal{A}:=\frac{D^2g((D_n{\bf u})_{x_0,\varrho}\otimes {\bf e}_n)}{\varphi''(1+|(D_n{\bf u})_{x_0,\varrho}|)}\equiv \frac{\left(D^2f(\cdot, {\bf 0} ,(D_n{\bf u})_{x_0,\varrho}\otimes {\bf e}_n)\right)_{x_0,\varrho}}{\varphi''(1+|(D_n{\bf u})_{x_0,\varrho}|)}\,.
\label{(4.1celadaok)}
\end{equation}
We notice that, by  \eqref{LHf}, $\mathcal{A}$ defined above satisfies the Legendre -- Hadamard  condition in \eqref{(2.20)}. Before stating the next lemma, we define
\begin{equation}\label{tildeomega1}
\tilde\omega_1(t):= \omega_1(t^{1/2})t^{1/2}\,,
\end{equation}
where $\omega_1(\cdot)$ is defined in \ref{ass-6f}. Then one can see that $\omega_1(t)/t$ and $\tilde \omega_1(t)/t$ are almost decreasing in $t\in(0,\infty)$ with the same constant $L\ge1$.


{\begin{lemma}\label{lem:lemma3.9}
Suppose that
\begin{equation}
\Phi(x_0,\varrho)\leq \varphi(1+|(D_n{\bf u})_{x_0,\varrho}|).
\label{(3.19)}
\end{equation}
Then, ${\bf u}$ is \emph{approximately $\mathcal{A}$-harmonic} on the upper half-ball $B_{\varrho/4}^+(x_0)$, in the sense that there exists $\beta_1=\beta_1(n,N,\mu_1,\mu_2,\nu,L)\in(0,\frac{1}{2})$ such that 
\begin{equation*}
\begin{split}
&\Biggl|\dashint_{B_{\varrho/4}^+(x_0)} \langle\mathcal{A}(D{\bf u}-(D_n{\bf u})_{x_0,\varrho}\otimes {\bf e}_n)|D\bm\eta\rangle\,\mathrm{d}x\Biggr| 
\le c (1+|(D_n{\bf u})_{x_0,\varrho}|)\|D\bm\eta\|_\infty \\
& \times \left\{H(x_0,\varrho)^{\beta_1}+\frac{\Phi(x_0,\varrho)}{\varphi(1+|(D_n{\bf u})_{x_0,\varrho}|)}+\tilde\omega_1\left(H(x_0,\varrho)^{\beta_1}+\frac{\Phi(x_0,\varrho)}{\varphi(1+|(D_n{\bf u})_{x_0,\varrho}|)}\right)\right\}
\end{split}
\end{equation*}
holds for every $\bm\eta\in C^\infty_0(B_{\varrho/4}^+(x_0);\R^N)$ for some constant $c = c(n,N,\mu_1, \mu_2,\nu, L)>0$, where
\begin{equation}
\beta_1:=\min\left\{\frac{1}{\mu_2},1-\frac{1}{\mu_1},1-\tau\right\}<\frac{1}{2}\,.
\label{beta1}
\end{equation}
In \eqref{beta1}, $\tau$ is that of Corollary~\ref{cor3.4celadaok}, and $\mu_1$ and $\mu_2$ are the characteristics of $\varphi$.
\end{lemma}}

\proof
The proof is exactly as in the interior setting \cite[Lemma~3.9]{GoodSciStroff2021} (see also \cite[Lemma~4.3]{CeladaOk}), considering the comparison map ${\bf v}\in {\bf u} + W^{1,\varphi}_0(B_{\varrho/2}^+,\R^N)$ provided by Lemma~\ref{lem:lemma3.7}, with minor changes due to the presence of function $\tilde\omega_1$ defined in \eqref{tildeomega1}, {that was $t^{\frac{\beta_0+1}{2}}$ for some $\beta_0>0$ in \cite{CeladaOk,GoodSciStroff2021}. For instance, we use Lemma~\ref{lem:Jensen} in the proof instead of H\"older's inequality.} We omit the details. 
\endproof

We introduce the \emph{hybrid excess functional}
\begin{equation}
\Phi_*(x_0,\varrho):=\Phi(x_0,\varrho)+\varphi(1+|(D_n{\bf u})_{x_0,\varrho}|) H(x_0,\varrho)^{\beta_1}\,.
\label{(4.11)}
\end{equation}
where $\beta_1$ is the exponent defined in \eqref{beta1}. Since $H(x_0,\varrho)\leq 1$, we deduce, in particular, that
$H(x_0,\varrho)\leq [H(x_0,\varrho)]^{\beta_1}$. Thus, the Caccioppoli inequality \eqref{eq:caccioppoliI} can be re-read as
\begin{equation*}
\Phi(x_0,\varrho/2)\leq c\Phi_*(x_0,\varrho)\,,
\end{equation*}
where $c=c(n,N,\mu_1,\mu_2,\nu,L)$.

\subsection{The excess-decay estimate}\label{sec:nondegenerate}

We start by establishing excess improvement estimates. The strategy of the proof is to approximate the given minimizer with $\mathcal{A}$-harmonic functions, for which suitable decay estimates are available from Theorem~\ref{thm:Aappr_phi}.

\begin{lemma}\label{lem:lemma3.12}
For every $\varepsilon\in(0,1)$ there exist $\delta_1,\delta_2\in(0,1]$, where $\delta_i=\delta_i(n,N,\mu_1,\mu_2, \nu,L,\omega_1(\cdot),\varepsilon)$, $i=1,2$, with the following property: if
\begin{align}
\frac{\Phi(x_0,\varrho)}{\varphi(1+|(D_n{\bf u})_{x_0,\varrho}|)}\leq \delta_1 \label{(4.5a)}\\
H(x_0,\varrho)^{\beta_1}\leq \delta_2 \label{(4.5b)}
\end{align}
then the excess improvement estimate
\begin{equation}
\Phi(x_0,\vartheta\varrho)\leq c_{\rm dec}\vartheta^2\left[1+\frac{\varepsilon}{\vartheta^{n+2}}\right]\Phi_*(x_0,\varrho)
\label{(4.6)}
\end{equation}
holds for every $\vartheta\in(0,1)$ for some constant $c_{\rm dec}=c_{\rm dec}(n,N,\mu_1,\mu_2,\nu,L)>0$, where $\Phi_*$ is defined in \eqref{(4.11)}.
\end{lemma}

\proof
The proof is an adaptation of \cite[Lemma~3.10]{GoodSciStroff2021} to the flat boundary setting. Without loss of generality, we assume that $x_0=0$, and we set $\bm\xi_\varrho:= (D_n{\bf u})_{0,\varrho}$ and ${\bf Q}_\varrho:=\bm\xi_\varrho \otimes {\bf e}_n$. 
We fix any 
{$\vartheta\in(0,\frac{1}{16}]$.}
Note that if 
{$\vartheta\in(\frac{1}{16},1)$,}
the inequality in \eqref{(4.6)} is trivial. 
 

\noindent
\textbf{Step 1: $\mathcal A$-harmonic approximation.}
For $\delta_1>0$ small enough, the assumption \eqref{(3.19)} is satisfied. Thus, if $\mathcal{A}$ is defined as in \eqref{(4.1celadaok)}, in view of Lemma~\ref{lem:lemma3.9}, for every $\bm\eta\in C^\infty_0(B_{\varrho/4};\R^N)$ with $\|D\bm\eta\|_\infty\leq1$ we have
\begin{equation}
\begin{split}
&\left|\dashint_{B_{\varrho/4}^+}\langle\mathcal{A}\left(\frac{D{\bf u}-{\bf Q}_{\varrho}}{1+|\bm\xi_\varrho|}\right) |D\bm\eta\rangle\,\mathrm{d}x\right|\\
&\leq {c}_1 \left\{\left(H(\varrho)^{\beta_1}+\frac{\Phi(\varrho)}{\varphi(1+|\bm\xi_\varrho|)}\right)^{\frac{1}{2}}+\omega_1\left(\left(H(\varrho)^{\beta_1}+\frac{\Phi(\varrho)}{\varphi(1+|\bm\xi_\varrho|)}\right)^{\frac{1}{2}}\right)\right\}\left(\frac{\Phi_*(\varrho)}{\varphi(1+|\bm\xi_\varrho|)}\right)^{\frac{1}{2}}\,,
\end{split}
\label{(4.10)}
\end{equation}
for some constant ${c}_1>0$. 
Note that, by assumptions \eqref{(4.5a)}-\eqref{(4.5b)}, it holds that
\begin{equation}
\frac{\Phi_*(\varrho)}{\varphi(1+|\bm\xi_\varrho|)}\leq \delta_1+\delta_2\,.
\label{(4.10bis)}
\end{equation}
Now, we define the $N$-function
\begin{equation*}
\psi(t):=\frac{\varphi_{1+|\bm\xi_\varrho|}((1+|\bm\xi_\varrho|)t)}{\varphi(1+|\bm\xi_\varrho|)}\,,\quad t\geq0\,,
\end{equation*}
and taking into account the properties of $\varphi$ and the fact that $\psi$ is, actually, a shifted function (see \cite[eq. (4.2)]{CeladaOk} for details),
we can show that, for a suitable constant ${c}_2>0$,
\begin{equation}
t^2\leq {c}_2 \psi(t)\,,\quad t\in[0,1]\,.
\label{(4.10zero)}
\end{equation}
From Lemma~\ref{corollary3.2}, we have also
\begin{equation}
\begin{split}
\Bigg(\dashint_{B_{\varrho/4}^+}\psi\left(\frac{|D{\bf u}-{\bf Q}_{\varrho}|}{1+|\bm\xi_\varrho|}\right)^{s_0}\,\mathrm{d}x\Bigg)^{\frac{1}{s_0}} & = \Bigg(\dashint_{B_{\varrho/4}^+}\left[\frac{\varphi_{1+|\bm\xi_\varrho|}(|D{\bf u}-{\bf Q}_{\varrho}|)}{\varphi(1+|\bm\xi_\varrho|)}\right]^{s_0}\,\mathrm{d}x\Bigg)^{\frac{1}{s_0}} \\
 & \leq \frac{c}{\varphi(1+|\bm\xi_\varrho|)}\dashint_{B_{\varrho}^+}\varphi_{1+|\bm\xi_\varrho|}(|D{\bf u}-{\bf Q}_{\varrho}|)\,\mathrm{d}x + c H(\varrho)^{\beta_1}\\
& \leq {c}_3 \frac{\Phi_*(\varrho)}{\varphi(1+|\bm\xi_\varrho|)}
\end{split}
\label{(4.10tris)}
\end{equation}
for some constant ${c}_3>0$. Then, with the constants ${c}_1, {c}_2, {c}_3$ determined above and \eqref{(4.10bis)}, by choosing $\delta_1$ and $\delta_2$ small enough we obtain
\begin{equation}
\kappa:=\max\{{c}_1, \sqrt{{c}_2{c}_3}\}\left[\frac{\Phi_*(\varrho)}{\varphi(1+|\bm\xi_\varrho|)}\right]^{\frac{1}{2}}\leq \max\{{c}_1, \sqrt{{c}_2{c}_3}\}\sqrt{\delta_1+\delta_2}<1\,.
\label{(4.13)}
\end{equation}
Now, combining \eqref{(4.10tris)}, \eqref{(4.13)} and \eqref{(4.10zero)}, we get
\begin{equation}
\Bigg(\dashint_{B_{\varrho/4}^+}\psi\left(\frac{|D{\bf u}-{\bf Q}_{\varrho}|}{1+|\bm\xi_\varrho|}\right)^{s_0}\,\mathrm{d}x\Bigg)^{\frac{1}{s_0}} \leq {c}_3\frac{\kappa^2}{\max\{{c}_1, \sqrt{{c}_2{c}_3}\}^2}\leq \frac{\kappa^2}{{c}_2}\leq \psi(\kappa)
\label{(4.14)}
\end{equation} 
and, plugging \eqref{(4.13)} and \eqref{(4.5a)} into \eqref{(4.10)}, we infer
\begin{equation*}
\dashint_{B_{\varrho/4}^+}\langle\mathcal{A}\left(\frac{D{\bf u}-{\bf Q}_{\varrho}}{1+|\bm\xi_\varrho|}\right)|D\bm\eta\rangle\,\mathrm{d}x \leq \frac{{c}_1\left(\sqrt{\delta_1+\delta_2}+\omega_1(\sqrt{\delta_1+\delta_2})\right)}{\max\{{c}_1, \sqrt{{c}_2{c}_3}\}}\kappa\,.
\end{equation*}
Therefore, up to choosing $\delta_1,\delta_2$ small enough, we can apply Theorem~\ref{thm:Aappr_phi} and the subsequent Remark~\ref{rem:thmmodified} to the function $\displaystyle{\bf v}:=\frac{{\bf u}-\bm\xi_\varrho x_n}{1+|\bm\xi_\varrho|}$ with $s=s_0$ and $\mu=\kappa$, so that there exists an $\mathcal{A}$-harmonic function ${\bf h}$ in $B_{\varrho/2}^+$, with ${\bf h}={\bf v}={\bf 0}$ on $\Gamma_{\varrho/2}$, such that
\begin{equation}
\frac{1}{\varphi(1+|\bm\xi_\varrho|)}\dashint_{B_{\varrho/4}^+}\varphi_{1+|\bm\xi_\varrho|}(|D{\bf u}-{\bf Q}_{\varrho}-(1+|\bm\xi_\varrho|)D{\bf h}|)\,\mathrm{d}x = \dashint_{B_{\varrho/4}^+} \psi(|D{\bf v}-D{\bf h}|)\,\mathrm{d}x\leq \varepsilon \psi(\kappa)\,.
\label{(4.15)}
\end{equation} 
In addition, since
\begin{equation*}
\psi(\kappa)\leq  c \frac{\varphi((1+|\bm\xi_\varrho|)(1+\kappa))}{\varphi(1+|\bm\xi_\varrho|)(1+\kappa)^2}\kappa^2\leq c\kappa^2\, \quad (\kappa\in(0,1)),
\end{equation*}
applying \eqref{(4.13)}, we conclude that
\begin{equation}
\dashint_{B_{\varrho/4}^+}\varphi_{1+|\bm\xi_\varrho|}(|D{\bf u}-{\bf Q}_{\varrho}-(1+|\bm\xi_\varrho|)D{\bf h}|)\,\mathrm{d}x\leq {c}\varepsilon\Phi_*(\varrho)\,.
\label{(4.16)}
\end{equation}

\noindent
\textbf{Step 2: Estimates for $\mathcal{A}$-harmonic maps.}
From \eqref{(4.14)} and \eqref{(4.15)} we infer that
\begin{equation*}
\dashint_{B_{\varrho/4}^+}\psi(|D{\bf h}|)\,\mathrm{d}x \leq c\dashint_{B_{\varrho/4}^+}\psi(|D{\bf v}|)\,\mathrm{d}x +  c\dashint_{B_{\varrho/4}^+}\psi(|D{\bf v}-{D\bf h}|)\,\mathrm{d}x \leq c(1+\varepsilon)\psi(\kappa) \leq c\psi(\kappa)\,,
\end{equation*}
which combined with the property \eqref{(2.21)}, Jensen's inequality and \eqref{(4.13)}, gives
\begin{equation}
\sup_{B_{\varrho/8}^+}|D{\bf h}|\leq c \psi^{-1}\Bigg(\dashint_{B_{\varrho/4}^+}\psi(|D{\bf h}|)\,\mathrm{d}x\Bigg)\leq c\kappa \leq  c \sqrt{\delta_1+\delta_2} \leq \frac{1}{2}
\label{eq:hineq}
\end{equation}
up to choose $\delta_1$ and $\delta_2$ small enough.   
Here, we recall $\vartheta\in(0,\frac{1}{16}]$.  Taking into account \eqref{eq:hineq},  we observe that 
\begin{equation}
\frac{1+|\bm\xi_\varrho|}{2}\leq 1+|{\bm\xi}_{\varrho}+(1+|\bm\xi_\varrho|)(D_n{\bf h})_{\vartheta\varrho}| \leq \frac{3(1+|\bm\xi_\varrho|)}{2}\,,
\label{(4.17)}
\end{equation}
and 
\begin{equation*}
(1+|\bm\xi_\varrho|)(1+(|D_n{\bf h}|)_{\varrho/4})\geq (1+|\bm\xi_\varrho|)(1+\vartheta(|D_n{\bf h}|)_{\varrho/4}|)\geq \frac{2}{3}(1+|\bm\xi_\varrho|)(1+(|D_n{\bf h}|)_{\varrho/4})
\end{equation*}
which together with \eqref{(2.6b)} and Jensen's inequality implies that
\begin{equation}\begin{aligned}
\varphi_{1+|\bm\xi_\varrho|}(\vartheta(1+|\bm\xi_\varrho|)(|D_n{\bf h}|)_{\varrho/4}) 
& \le c \vartheta^2 \varphi_{1+|\bm\xi_\varrho|}((1+|\bm\xi_\varrho|)(|D_n{\bf h}|)_{\varrho/4}) \\
& \leq \vartheta^2 \dashint_{B_{\varrho/4}^+} \varphi_{1+|\bm\xi_\varrho|}((1+|\bm\xi_\varrho|)|D_n{\bf h}|)\,\mathrm{d}x\,.
\end{aligned}
\label{(4.18bis)}
\end{equation}

In addition, applying the Calder\'on-Zygmund estimate \eqref{eq:caldezyg} to $\bfh$, we have
\begin{equation}
\begin{split}
\dashint_{B_{\varrho/4}^+}\varphi_{1+|\bm\xi_\varrho|}((1+|\bm\xi_\varrho|)|D{\bf h}|)\,\mathrm{d}x 
&= {\varphi(1+|\bm\xi_\varrho|)\dashint_{B_{\varrho/4}^+} \psi(|D{\bf h}|)\,\mathrm{d}x} \\
&{\leq c \varphi(1+|\bm\xi_\varrho|) \dashint_{B_{\varrho/4}^+} \psi(|D{\bf v}|)\,\mathrm{d}x}  \leq c \Phi(\varrho)\,.
\end{split}
\label{cald-zygm}
\end{equation}

Finally, since $\bfh-(D_n{\bf h})_{\vartheta\varrho}x_n$ is also an $\mathcal A$-harmonic map with $\bfh-(D_n{\bf h})_{\vartheta\varrho}x_n={\bf 0}$ on $\Gamma_{\varrho/4}$, 
in view of Theorem~\ref{thm:Aharmonic} together with the Poincar\'e inequality, we have

\begin{equation}\label{oscillationdecay}
\begin{aligned}
\sup_{B_{\vartheta\varrho}^+}|D{\bf h}-(D_n{\bf h})_{\vartheta\varrho}\otimes {\bf e}_n|
& \leq c\dashint_{B_{2\vartheta\varrho}^+}|D_n{\bf h}-(D_n{\bf h})_{\vartheta\varrho}|\,\mathrm{d}x \\
& \leq c\vartheta\varrho \dashint_{B_{2\vartheta\varrho}^+}|D^2{\bf h}|\,\mathrm{d}x  \leq  c\vartheta\varrho \sup_{B_{\varrho/8}^+}|D^2{\bf h}| \leq  c\vartheta \dashint_{B_{\varrho/4}^+}|D_n{\bf h}|\,\mathrm{d}x.
\end{aligned}
\end{equation}

\noindent
\textbf{Step 3: Estimation of $\Phi(\vartheta\varrho)$.}
First, we denote by $\overline\bfu$ and $\overline\bfh$ the odd extensions of  $\bfu$ and $\bfh$ in $B_{\varrho/4}$, respectively. Then we note that, for every $r\in(0,\varrho/4]$,
$$
(D\overline\bfu)_{B_r}= \Bigg(\dashint_{B_r^+}D_n\bu \,\dd x\Bigg) \otimes {\bf e}_n   = {\bm\xi}_r \otimes {\bf e}_n = {\bf Q}_r,
\quad \text{hence } \  |(D\overline\bfu)_{B_r}|=|(D_n\bu)_r| = |{\bm\xi}_r|,
$$
and also
$$
(D\overline\bfh)_{B_r}= ( D_n\bfh)_{r} \otimes {\bf e}_n\,.
$$
Now, we estimate $\Phi(\vartheta\varrho)$. By \eqref{eq:equivalence} and \eqref{eq:equivalencebis} with the preceding identities, we obtain 
\begin{equation*}
\begin{aligned}
\Phi(\vartheta\varrho)
&= \dashint_{B_{\vartheta\varrho}^+}\varphi_{1+|{\bm\xi}_{\vartheta\varrho}|}(|D\bfu-{\bf Q}_{\vartheta\varrho}|)\,\mathrm{d}x= \dashint_{B_{\vartheta\varrho}}\varphi_{1+|(D\overline{\bfu})_{B_{\vartheta\varrho}}|}(|D{\overline\bfu}-(D\overline{\bfu})_{B_{\vartheta\varrho}}|)\,\mathrm{d}x\\
&\leq c  \dashint_{B_{\vartheta\varrho}}|{\bf V}_1(D{\overline{\bf u}})-{\bf V}_1((D\overline{\bfu})_{B_{\vartheta\varrho}})|^2\,\mathrm{d}x
\leq  c  \dashint_{B_{\vartheta\varrho}}|{\bf V}_1(D{\overline{\bf u}})-({\bf V}_1(D\overline{\bfu}))_{B_{\vartheta\varrho}}|^2\,\mathrm{d}x\\
&\leq  \dashint_{B_{\vartheta\varrho}}|{\bf V}_1(D{\overline{\bfu}})-{\bf V}_1\left((D{\overline\bfu})_{B_\varrho}+(1+|(D{\overline\bfu})_{B_\varrho}|)(D{\overline\bfh})_{B_{\vartheta\varrho}}\right)|^2\,\mathrm{d}x\\
&\leq c \dashint_{B_{\vartheta\varrho}}\varphi_{1+|(D{\overline\bfu})_{B_\varrho}+(1+|(D{\overline\bfu})_{B_\varrho}|)(D{\overline\bfh})_{B_{\vartheta\varrho}}|}\left(\left|D\overline\bfu - [(D{\overline\bfu})_{B_\varrho}+(1+|(D{\overline\bfu})_{B_\varrho}|)(D{\overline\bfh})_{B_{\vartheta\varrho}}]\right|\right)\,\mathrm{d}x\\
&\leq c \dashint_{B_{\vartheta\varrho}^+}\varphi_{1+|{\bm \xi}_{\varrho}+(1+|\bm\xi_\varrho|)(D_n{\bf h})_{\vartheta\varrho}|}\left(\left|D\bfu - {\bf Q}_{\varrho}-(1+|\bm\xi_\varrho|)(D_n{\bf h})_{\vartheta\varrho}\otimes {\bf e}_n\right|\right)\,\mathrm{d}x\,.
\end{aligned}
\end{equation*}
Then, from \eqref{(4.17)}, we immediately infer that
\begin{equation*}
\Phi(\vartheta\varrho) \leq c \dashint_{B_{\vartheta\varrho}^+}\varphi_{1+|{\bm \xi}_{\varrho}|}\left(\left|D\bfu - {\bf Q}_{\varrho}-(1+|\bm\xi_\varrho|)(D_n{\bf h})_{\vartheta\varrho}\otimes {\bf e}_n\right|\right)\,\mathrm{d}x\,.
\end{equation*}
Finally, by \eqref{oscillationdecay},   \eqref{(4.16)}, \eqref{(4.18bis)} and \eqref{cald-zygm},
\begin{equation*}
\begin{split}
\Phi(\vartheta\varrho) & \leq c\dashint_{B_{\vartheta\varrho}^+}\varphi_{1+|\bm\xi_\varrho|}(|D{\bf u}-{\bf Q}_{\varrho}-(1+|\bm\xi_\varrho|)D{\bf h}|)\,\mathrm{d}x \\
&\qquad + \dashint_{B_{\vartheta\varrho}^+}\varphi_{1+|\bm\xi_\varrho|}((1+|\bm\xi_\varrho|)|D{\bf h}-(D_n{\bf h})_{\vartheta\varrho}\otimes {\bf e}_n|)\,\mathrm{d}x \\
& \leq c \vartheta^{-n} \dashint_{B_{\varrho}^+}\varphi_{1+|\bm\xi_\varrho|}(|D{\bf u}-{\bf Q}_{\varrho}-(1+|\bm\xi_\varrho|)D{\bf h}|)\,\mathrm{d}x + \varphi_{1+|\bm\xi_\varrho|}(\vartheta(1+|\bm\xi_\varrho|)(|D_n{\bf h}|)_{\varrho/4}) \\
& \leq c \vartheta^{-n} \varepsilon\Phi_*(\varrho)  + c \vartheta^2 \dashint_{B_{\varrho/4}^+} \varphi_{1+|\bm\xi_\varrho|}((1+|\bm\xi_\varrho|)|D{\bf h}|)\,\mathrm{d}x \\
& \leq c \vartheta^{-n} \varepsilon\Phi_*(\varrho)  + c \vartheta^2 \Phi(\varrho)
 \leq c \vartheta^2\left(1+\frac{\varepsilon}{\vartheta^{n+2}}\right)\Phi_*(\varrho)\,,
\end{split}
\end{equation*}
and the proof is complete.
\endproof

\subsection{Proof of Theorem~\ref{theorem-result-1}}\label{lastsection}

We first prove the partial regularity for the minimizers $\bfu\in W^{1,\phi}_\Gamma(B_R^+)$ of \eqref{functional} with $\Omega=B_R^+$,

\begin{theorem}
\label{theorem-result-2}
Let  $\phi$ be an $N$-function complying with {\rm\ref{ass-1phi}} -- {\rm\ref{ass-2phi}}, and $f:B_{R}^+\times \R^N\times \R^{N\times n}\to \R$ satisfies  {\rm\ref{ass-1f}} -- {\rm\ref{ass-6f}} with $\Omega=B_R^+$. If  $\bfu\in W^{1,\phi}_\Gamma(B_R^+)$ is a minimizer of \eqref{functional} with $\Omega=B_R^+$, then the set of regular points on ${\Gamma_R}$ given by 
\begin{equation*}
\Gamma_{\bfu}:=\bigcap_{\alpha\in(0,1)}\left\{x_0\in{\Gamma_R}:\,\, \bfu\in C^{\alpha}(U_{x_0}\cap \overline{B_R^+};\R^N) \mbox{ for some }U_{x_0}\subset B_R\right\}\,,
\end{equation*}
where $U_{x_0}$ is an open neighborhood of $x_0$, satisfies
${\Gamma_R}\backslash\Gamma_{\bfu}\subset\mathrm{Sing}_{\bu}(\Gamma_R)$ where
\begin{equation*}
{\begin{split}
\mathrm{Sing}_{\bu}(\Gamma_R):=& \Bigg\{x_0\in\Gamma_R:\,\, \mathop{\lim\inf}_{\varrho\searrow 0} \dashint_{B^+_\varrho(x_0)} |D{\bf u}-(D_{n}{\bf u})_{x_0,\varrho}\otimes {\bf e}_n|\,\mathrm{d}x>0\bigg\}\\
&\ \cup\left\{x_0\in\Gamma_R:\,\, \mathop{\lim\sup}_{\varrho\searrow 0}\,(|D_{n}{\bf u}|)_{x_0,\varrho}=+\infty\right\}\,.
\end{split}}
\end{equation*} 
\end{theorem}

\proof
We fix any $\alpha\in(0,1)$, and set
\begin{equation}\label{lambda}
\lambda:=n-\mu_1(1-\alpha)\in(n-\mu_1,n)\,.
\end{equation}

\noindent
{\bf Step 1: Boundary decay estimates and Morrey type estimates.} 
We start by fixing various parameters. First, we choose $\vartheta>0$ such that 
\begin{equation}
\vartheta\leq \vartheta_1:=\min\left\{(6\cdot 2^{\mu_2}c_{\rm dec})^{-\frac{1}{2}},\frac{1}{2^{\frac{\mu_2}{\mu_1}}}, \frac{1}{2^\frac{\mu_2}{n-\lambda}}\right\}<1\,,
\label{eq:choosetheta}
\end{equation}
where $c_{\rm dec}$ is the constant in Lemma~\ref{lem:lemma3.12}, hence $\vartheta_1$ depends on $n,N,\mu_1,\mu_2,\nu,L,\alpha$. Correspondingly, let $\delta_i=\delta_i(n,N,\mu_1,\mu_2,\nu,L,\omega_1(\cdot),\vartheta)$, $i=1,2$, be the constants determined in Lemma~\ref{lem:lemma3.12}, applied with the choice $\varepsilon=\vartheta^{n+2}$. Then we choose $\varepsilon_*>0$ such that
\begin{equation}
\varepsilon_*\leq\epsilon_1:=\min\left\{\delta_1, \frac{\delta_2}{2}, \frac{\vartheta^n}{2c_{\eta_1}},\frac{\vartheta^n}{2^{\mu_2+1}c_{\eta_2}}\right\}\,,
\label{(3.46verena)}
\end{equation}
where $c_{\eta_1}$ and $c_{\eta_2}$ are the constants in the change-shift formula \eqref{(5.4diekreu)} with $\eta_1:=\frac{1}{2}$ and $\eta_2:=2^{-\mu_2-1}$, respectively. Finally, with a fixed $\epsilon_*$, we choose constants $\kappa_1,\varrho_1>0$ so small that
\begin{equation}
\left(\frac{\omega(\kappa_1)^{1-\frac{1}{s}}}{1+(2L)^{1-\frac{1}{s}}}\right)^{\beta_1}=\varepsilon_* 
\quad\text{and}\quad
\left(\frac{\mathcal{V}(\varrho_1)^{1-\frac{1}{s}}}{1+(2L)^{1-\frac{1}{s}}}\right)^{\beta_1}=\varepsilon_*\,,
\label{(3.47verena)}
\end{equation}
where $s>1$ is the exponent appearing in the definitions of $H$ and $\widetilde{H}$ in \eqref{(3.8)} and \eqref{(3.8bis)}, and $\beta_1$ is given in \eqref{beta1}.
As a consequence, $\epsilon_1$ depends on the choice of  $\vartheta$, $\kappa_1$ does on the choice of $\varepsilon_*$ and $\omega(\cdot)$, and $\varrho_1$ on that of $\varepsilon_*$ and $\mathcal{V}(\cdot)$.

Now, we prove that the excess-decay estimate in \eqref{(4.6)} can be iterated, as the  conditions in  \eqref{(4.5a)} and \eqref{(4.5b)} are also satisfied on any smaller upper half-ball $B_{\vartheta^m\varrho}^+(x_0)\subset B^+_R$ with $m\in\mathbb{N}$, $x_0\in \R^{n-1}\times \{0\}$ and  $\rho\in(0, \rho_1]$. Precisely, we prove the following:
if the conditions 
\begin{equation}
\frac{\Phi(x_0,\varrho)}{\varphi(1+|(D_n{\bf u})_{x_0,\varrho}|)}\leq \varepsilon_* 
\quad \mbox{ and }\quad 
\varrho\, \varphi^{-1}\Bigg(\dashint_{B_{\varrho}^+(x_0)}\varphi(1+|D{\bf u}|)\,\mathrm{d}x\Bigg)\leq\kappa_* 
\label{eq:0step}
\end{equation}
hold for some $\epsilon_*\in(0,\epsilon_1]$,  $\kappa_*\in(0,\kappa_1]$ and $\varrho\in(0,\varrho_1]$, where $\Phi(x_0,\varrho)$ is defined in \eqref{(3.5)}, then 
\begin{equation}
\frac{\Phi(x_0,\vartheta^m\varrho)}{{\varphi(1+|(D_n{\bf u})_{x_0,\vartheta^m\varrho}|)}}\leq \varepsilon_*
 \quad \mbox{ and }\quad 
\vartheta^m\varrho\,\varphi^{-1}\Bigg(\dashint_{B_{\vartheta^m\varrho}^+(x_0)}\varphi(1+|D{\bf u}|)\,\mathrm{d}x\Bigg) \leq\kappa_*
\label{eq:kstep}
\end{equation}
for every $m=0,1,2,\dots.$.

We argue by induction on $m$. As usual, we omit the explicit dependence on $x_0$. Since \eqref{eq:kstep} are trivially true for $m=0$ by assumption \eqref{eq:0step}, our aim is to show that if \eqref{eq:kstep} holds for some $m\geq1$, then the corresponding inequalities hold with $m+1$ in place of $m$.
We first prove the second inequalities in \eqref{eq:kstep}.
Applying the shift-change formula \eqref{(5.4diekreu)} with $\eta=\eta_1=\frac{1}{2}$, the first inequality in \eqref{eq:kstep} at step $m$, and the facts that $c_{\eta_1}\vartheta^{-n}\varepsilon_*\le \frac{1}{2}$ from  \eqref{(3.46verena)} and $ 2^{\mu_2}\vartheta^{\mu_1} \leq1$ from \eqref{eq:choosetheta}, and Jensen's inequality, we obtain
\begin{equation}\label{estim1}
\begin{aligned}
\dashint_{B_{\vartheta^{m+1}\varrho}^+}\varphi(1+|D\bfu|) \,\mathrm{d}x 
& \leq 2^{\mu_2-1}\Bigg(\dashint_{B_{\vartheta^{m+1}\varrho}^+}\varphi(|D\bfu-(D_n{\bf u})_{\vartheta^m\varrho}\otimes{\bf e}_n|)\,\mathrm{d}x +\phi(1+|(D_n{\bf u})_{\vartheta^m\varrho}|) \Bigg) \\
& \leq 2^{\mu_2-1}\left(c_{\eta_1} \vartheta^{-n} \Phi(\vartheta^m\varrho)+\frac{3}{2}\varphi(1+|(D_n{\bf u})_{\vartheta^m\varrho}|)\right) \\
& \leq 2^{\mu_2-1}\left(c_{\eta_1}\vartheta^{-n} \varepsilon_* + \frac 32 \right) \varphi(1+|(D_n{\bf u})_{\vartheta^m\varrho}|) \\
& \leq 2^{\mu_2-1}\left(c_{\eta_1}\vartheta^{-n} \varepsilon_* + \frac 32 \right)\dashint_{B_{\vartheta^{m}\varrho}^+}\varphi(1+|D\bfu|)\,\dd x \\
& \leq 2^{\mu_2}\varphi\left(\frac{\kappa_*}{\vartheta^{m}\varrho}\right) 
\leq 2^{\mu_2}\vartheta^{\mu_1} \varphi\left(\frac{\kappa_*}{\vartheta^{m+1}\varrho}\right) 
\leq  \phi \left(\frac{\kappa_*}{\vartheta^{m+1}\varrho}\right) 
\end{aligned}
\end{equation}
which corresponds to \eqref{eq:kstep} with $m+1$ in place of $m$. 

Next, we prove  the first inequality in \eqref{eq:kstep} with $m+1$ in place of $m$.
From \eqref{eq:kstep} at step $m$ and the choices of $\varepsilon_*$, $\kappa_1$ and $\varrho_1$ as in \eqref{(3.46verena)}-\eqref{(3.47verena)}, we have
\begin{align*}
&\frac{\Phi(\vartheta^m\varrho)}{{\varphi(1+|(D_n{\bf u})_{\vartheta^m\varrho}|)}}\leq \varepsilon_*\leq\delta_1\,,\\
&H(\vartheta^m\varrho)^{\beta_1}\leq \left[\frac{\omega(\kappa_*)^{1-\frac{1}{s}} + {\mathcal{V}}(\varrho)^{1-\frac{1}{s}}}{1+(2L)^{1-\frac{1}{s}}}\right]^{\beta_1} 
\leq \left[\frac{\omega(\kappa_1)^{1-\frac{1}{s}}}{1+(2L)^{1-\frac{1}{s}}}\right]^{\beta_1}+\left[\frac{{\mathcal{V}}(\varrho_1)^{1-\frac{1}{s}}}{1+(2L)^{1-\frac{1}{s}}}\right]^{\beta_1}<2\varepsilon_*\leq \delta_2\,,
\end{align*}
and
\begin{equation*}
\begin{split}
\Phi_*(\vartheta^m\varrho)=\Phi(\vartheta^m\varrho) + \varphi(1+|(D_n{\bf u})_{\vartheta^m\varrho}|) [H(\vartheta^m\varrho)]^{\beta_1} & \leq 3\epsilon_*\varphi(1+|(D_n{\bf u})_{\vartheta^m\varrho}|)\,.
\end{split}
\end{equation*}
Moreover, using the shift-change formula \eqref{(5.4diekreu)} with $\eta=\eta_2:=2^{-\mu_2-1}$, the first inequality in \eqref{eq:kstep} at step $m$  and  the fact that $2^{\mu_2-1} c_{\eta_2} \vartheta^{-n}\epsilon_* \leq  \frac{1}{4}$ from \eqref{(3.46verena)},
$$\begin{aligned}
\varphi(1+|(D_n{\bf u})_{\vartheta^m\varrho}|) &  \leq 2^{\mu_2-1}\varphi(1+|(D_n{\bf u})_{\vartheta^{m+1}\varrho}|) +  2^{\mu_2-1}\dashint_{B_{\vartheta^{m+1}\varrho}^+}\varphi(|D_n{\bf u}-(D_n{\bf u})_{\vartheta^{m}\varrho}|)\, \dd x\\
& \leq 2^{\mu_2-1} \varphi(1+|(D_n{\bf u})_{\vartheta^{m+1}\varrho}|) +  2^{\mu_2-1} c_{\eta _2}\vartheta^{-n}\Phi(\vartheta^{m}\varrho) +\tfrac{1}{4} \varphi(1+|(D_n{\bf u})_{\vartheta^{m}\varrho}|)\\
& \leq 2^{\mu_2-1} \varphi(1+|(D_n{\bf u})_{\vartheta^{m+1}\varrho}|) +  (2^{\mu_2-1} c_{\eta_2} \vartheta^{-n}\epsilon_* +\tfrac{1}{4})  \varphi(1+|(D_n{\bf u})_{\vartheta^{m}\varrho}|)\\
& \leq 2^{\mu_2-1} \varphi(1+|(D_n{\bf u})_{\vartheta^{m+1}\varrho}|) +  \tfrac{1}{2}  \varphi(1+|(D_n{\bf u})_{\vartheta^{m}\varrho}|).
\end{aligned}$$
Therefore, we obtain
$$
\varphi(1+|(D_n{\bf u})_{\vartheta^m\varrho}|)  \leq 2^{\mu_2} \varphi(1+|(D_n{\bf u})_{\vartheta^{m+1}\varrho}|),
$$ 
hence
$$
\Phi_*(\vartheta^m\varrho) \leq 3\cdot 2^{\mu_2}\epsilon_*\varphi(1+|(D_n{\bf u})_{\vartheta^{m+1}\varrho}|)\,.
$$
Finally, by virtue of Lemma~\ref{lem:lemma3.12}, applied with radius $\vartheta^m\varrho$ in place of $\varrho$, and the fact that  $6\cdot 2^{\mu_2}c_{\rm dec}\vartheta^2\leq 1$ from \eqref{eq:choosetheta}, we get
\begin{equation*}
\begin{split}
\Phi(\vartheta^{m+1}\varrho)\leq 2c_{\rm dec}\vartheta^2\Phi_*(\vartheta^m\varrho)&\leq 6\cdot 2^{\mu_2}c_{\rm dec}\epsilon_*\vartheta^2\varphi(1+|(D_n{\bf u})_{\vartheta^{m+1}\varrho}|) \\
& \leq \epsilon_* \varphi(1+|(D_n{\bf u})_{\vartheta^{m+1}\varrho}|)\,
\end{split}
\end{equation*}
which yields  the first inequality in \eqref{eq:kstep} with  $m+1$ in place of $m$.


Consequently, from \eqref{eq:kstep}, applying the shift-change formula \eqref{(5.4diekreu)} with $\eta=\eta_1=\frac{1}{2}$, \eqref{eq:kstep} and \eqref{(3.46verena)}, we obtain
\begin{equation*}
\begin{split}
\dashint_{B_{\vartheta^{m+1}\varrho}^+}\varphi(1+|D{\bf u}|)\,\mathrm{d}x & \leq 2^{\mu_2-1}\left(c_{\eta_1}\vartheta^{-n} \Phi(\vartheta^m\varrho)+\tfrac{3}{2}\varphi(1+|(D_n\bfu)_{\vartheta^m\varrho}|)\right) \\
& \leq \left(2^{\mu_2-1}c_{\eta_1}\vartheta^{-n} \varepsilon_* +3\cdot2^{\mu_2-2}\right)\dashint_{B_{\vartheta^{m}\varrho}^+}\varphi(1+|D{\bf u}|)\,\mathrm{d}x \\
& \leq 2^{\mu_2} \Bigg(\dashint_{B_{\vartheta^{m}\varrho}^+}\varphi(1+|D{\bf u}|)\,\mathrm{d}x  \Bigg)\,,\\
\end{split}
\end{equation*}
whence, by \eqref{eq:choosetheta},
$$
\int_{B_{\vartheta^{m+1}\varrho}^+}\varphi(1+|D{\bf u}|)\,\mathrm{d}x  \leq 2^{\mu_2} \vartheta^n \int_{B_{\vartheta^{m}\varrho}^+}\varphi(1+|D{\bf u}|)\,\mathrm{d}x \leq \vartheta^{\lambda}\int_{B_{\vartheta^{m}\varrho}^+}\varphi(1+|D{\bf u}|)\,\mathrm{d}x \,. 
$$

Consequently, by iterating the above inequality for $m=0,1,2,\dots$, we have that for every $r\in(0,\varrho]$,
\begin{equation}
\int_{B_{r}^+}\varphi(1+|D{\bf u}|)\,\mathrm{d}x \leq \vartheta^{-\lambda} \left(\frac{r}{\varrho}\right)^{\lambda}\int_{B_{\varrho}^+}\varphi(1+|D{\bf u}|)\,\mathrm{d}x\,,
\label{morrey-boundary}
\end{equation}
if the inequalities in \eqref{eq:0step} hold.

\noindent
{\bf Step 2: Interior decay estimates and Morrey type estimates.} 

In Step 1, we derived Morrey-type estimates \eqref{morrey-boundary} on half balls by proving excess-decay estimates in \eqref{eq:kstep}. 
Moreover, in the same way with minor modifications, we can also get the interior counterpart of  \eqref{morrey-boundary} on balls in $B_r^+$ by proving interior versions of  Lemma~\ref{lem:lemma3.9} and the results in Step 1. In fact,  these were obtained in \cite{GoodSciStroff2021} for degenerate problems. Therefore, we only state the interior counterpart of the results in Step 1 without proofs.

There exists a small $\vartheta_2=\vartheta_2(n,N,\mu_1,\mu_2,\nu,L,\alpha)\in(0,1)$ such that the following holds:  for any $\vartheta\in(0,\vartheta_2]$ there exists $\epsilon_2=(n,N,\mu_1,\mu_2,\nu,L,\omega_1(\cdot),\vartheta)$ such that if  $\epsilon_*\leq \epsilon_2$ and if $\kappa_*\in(0,\kappa_2)$ and $B_{\varrho}(y)\subset B_R^+$ with $\varrho\in(0,\varrho_2]$, where $\kappa_2$ and $\varrho_2$ are the constants satisfying the interior counterpart of \eqref{(3.47verena)}, then
\begin{equation}
\frac{\Phi_{\text{int}}(y,\vartheta^m\varrho)}{{\varphi(1+|(D{\bf u})_{B_{\vartheta^m\varrho}(y)}|)}}\leq \varepsilon_*
 \quad \mbox{ and }\quad 
\vartheta^m\varrho\,\varphi^{-1}\Bigg(\dashint_{B_{\vartheta^m\varrho}(y)}\varphi(1+|D{\bf u}|)\,\mathrm{d}x\Bigg) \leq\kappa_*
\label{eq:kstepint}
\end{equation}
for every $m=0,1,2,\dots$, where 
$$
\Phi_{\text{int}}(y,r) :=\dashint_{B_{r}(y)}\varphi_{1+|(D{\bf u})_{B_r(y)}|}(|D{\bf u}-(D{\bf u})_{B_r(y)}|)\, \dd x\,,
$$
if \eqref{eq:kstepint} hold when $m=0$. Moreover, \eqref{eq:kstepint} implies
\begin{equation}
\int_{B_{r}(y)}\varphi(1+|D{\bf u}|)\,\mathrm{d}x \leq \vartheta^{-\lambda} \left(\frac{r}{\varrho}\right)^{\lambda}\int_{B_{\varrho}(y)}\varphi(1+|D{\bf u}|)\,\mathrm{d}x.
\label{morrey-interior}
\end{equation}

\noindent
{\bf Step 3: $C^{\alpha}$-continuity.} 
Now, we fix $\vartheta:=\min\{\vartheta_1,\vartheta_2\}$. Then $\epsilon_1$ and $\epsilon_2$ are determined, and set $\epsilon_*:=\min\{\epsilon_1,\epsilon_2\}$, so that also $\kappa_1$, $\kappa_2$, $\varrho_1$, and $\varrho_2$ are determined. Finally we set $\kappa_*:=\min\{\kappa_1,\kappa_2\}$.

Suppose that for $B_{3\varrho_0}^+(x_0)\subset B_R^+$ with $0< \varrho_0\leq \frac{\min\{\varrho_1,\varrho_2\}}{3}$ and $x_0\in\Gamma_R$,  the following holds:
\begin{equation}
\frac{\Phi(\tilde y,2\varrho_0)}{\varphi(1+|(D_n{\bf u})_{\tilde y,2\varrho_0}|)}\leq \frac{\vartheta^n\varepsilon_*}{c_1}\leq \varepsilon_* 
\quad \mbox{ and }\quad 
2\varrho_0\, \varphi^{-1}\Bigg(\dashint_{B_{2\varrho_0}^+(\tilde y)}\varphi(1+|D{\bf u}|)\,\mathrm{d}x\Bigg)\leq\kappa_*
\label{eq:0step1}
\end{equation}
for every $\tilde y \in \Gamma_{\varrho_0}(x_0)$, where $c_1\geq 1$ is determined in \eqref{c1} below. Then, we prove that $\bfu \in C^{\alpha}(B_{\varrho_0}^+(x_0))$.

Let  $\overline \bfu$ be the odd extension of $\bfu$ in $B_{3\varrho_0}(x_0)$, and $y=(y',y_n)\in B^+_{\varrho_0}(x_0)\cup \Gamma_{\varrho_0}(x_0)$ be arbitrary, and set $\tilde y=(y',0)$. We first assume $r \in [y_n, \rho_0]$. Then, since $B_{r}(y)\subset B_{2r}(\tilde y) \subset B_{2\rho_0}(\tilde y)\subset B_{3\rho_0}(x_0)$, by \eqref{eq:0step1}, we have \eqref{morrey-boundary} with $B_{2\varrho_0}^+(\tilde y)$ in place of $B_\varrho^+$, hence
\begin{equation}
\begin{aligned}
\int_{B_{r}(y)}\varphi(1+|D{\overline\bfu}|)\,\mathrm{d}x
&\leq 2\int_{B_{2r}^+(\tilde y)}\varphi(1+|D{\bf u}|)\,\mathrm{d}x
\leq 2 \vartheta^{-\lambda} \left(\frac{r}{\rho_0}\right)^{\lambda}\int_{B_{2\varrho_0}^+(\tilde y)}\varphi(1+|D{\bf u}|)\,\mathrm{d}x\\
& \leq c \left(\frac{r}{\rho_0}\right)^{\lambda}\int_{B^+_{3\varrho_0}(x_0)}\varphi(1+|D{\bfu}|)\,\mathrm{d}x\, 
\quad  \text{for all }\ r\in [y_n, \rho_0].
\end{aligned}\label{morrey-boundary1}
\end{equation}

We next consider  $r\in(0,y_n)$, when $y_n>0$. Note that $\vartheta^{m+1}\varrho\leq y_n < \vartheta^m\varrho$ for some $m\in\mathbb N\cup\{0\}$. Then, by \eqref{eq:equivalence} and \eqref{eq:equivalencebis},
$$\begin{aligned}
&\Phi_{\text{int}}(y,y_n) \leq c \dashint_{B_{y_n}(y)} |\bfV_1(D\bfu)-(\bfV_1(D\bfu))_{B_{y_n}(y)}|^2\,\dd x\\
&\leq c\vartheta^{-n} \dashint_{B_{\vartheta^m\varrho}(\tilde y)} |\bfV_1(D\overline\bfu)-(\bfV_1(D\overline\bfu))_{B_{\vartheta^m\varrho}(\tilde y)}|^2\,\dd x \\
&\leq c\vartheta^{-n} \dashint_{B_{\vartheta^m\varrho}(\tilde y)} |\bfV_1(D\overline\bfu)-\bfV_1((D\overline\bfu)_{B_{\vartheta^m\varrho}(\tilde y)})|^2\,\dd x  \\
&\leq c \vartheta^{-n} \dashint_{B_{\vartheta^m\varrho}(\tilde y)} \phi_{1+|(D\overline\bfu)_{B_{\vartheta^m\varrho}(\tilde y)}|}(|D\overline\bfu-(D\overline\bfu)_{B_{\vartheta^m\varrho}(\tilde y)}|)\,\dd x\\
&\leq  c \vartheta^{-n} \dashint_{B^+_{\vartheta^m\varrho}(\tilde y)} \phi_{1+|(D_n\bfu)_{\tilde y,\vartheta^m\varrho}|}(|D\bfu-(D_n\bfu)_{\tilde y, \vartheta^m\varrho}\otimes {\bf e}_n|)\,\dd x  \\
&= c \vartheta^{-n} \Phi(\tilde y, \vartheta^m\varrho)\,,\\
\end{aligned}$$
where we used the fact that $(D\overline\bfu)_{B_{\vartheta^m\varrho}(\tilde y)}=(D_n\bfu)_{\tilde y,\vartheta^m\varrho}\otimes {\bf e}_n$. Moreover,
$$\begin{aligned}
\phi(1+|(D_n\bfu)_{\tilde y, \vartheta^m\varrho}|) &  \leq 2^{\mu_2-1}\varphi(1+|(D_n{\bf u})_{B_{y_n}(y)}|) +  2^{\mu_2-1}\dashint_{B_{y_n}(y)}\varphi(|D_n{\bf u}-(D_n\bfu)_{\tilde y, \vartheta^m\varrho}|)\, \dd x\\
& \leq 2^{\mu_2-1} \varphi(1+|(D_n{\bf u})_{B_{y_n}(y)}|) +  2^{\mu_2-1} c_{\eta _2}\vartheta^{-n}\Phi(\tilde y,\vartheta^{m}\varrho) +\tfrac{1}{4} \varphi(1+|(D_n{\bf u})_{\tilde y, \vartheta^{m}\varrho}|)\\
& \leq 2^{\mu_2-1} \varphi(1+|(D_n{\bf u})_{B_{y_n}(y)}|) +  (2^{\mu_2-1} c_{\eta_2} \vartheta^{-n}\epsilon_* +\tfrac{1}{4})  \varphi(1+|(D_n{\bf u})_{\tilde y, \vartheta^{m}\varrho}|)\\
& \leq 2^{\mu_2-1} \varphi(1+|(D_n{\bf u})_{B_{y_n}(y)}|) +  \tfrac{1}{2}  \varphi(1+|(D_n{\bf u})_{\tilde y, \vartheta^{m}\varrho}|).
\end{aligned}$$
From the preceding two inequalities, we obtain
\begin{equation}\label{c1}
\frac{\Phi_{\text{int}}(y,y_n)}{\Phi(1+|(D_n\bfu)_{B_{y_n}(y)}|)} \leq   c_1 \vartheta^{-n} \frac{ \Phi(\tilde y, \vartheta^m\varrho)}{\phi(1+|(D_n\bfu)_{\tilde y, \vartheta^m\varrho}|)}  
\end{equation}
for some $c_1=c_1(n,N,\mu_1,\mu_2,\nu,L)\geq 1$. Therefore, by \eqref{eq:0step1}, we prove  that
\begin{equation}
\frac{\Phi_{\text{int}}(y,y_n)}{\Phi(1+|(D_n\bfu)_{B_{y_n}(y)}|)} \leq \varepsilon_* 
\label{eq:0step2}
\end{equation}
In addition, by arguing as for \eqref{estim1}, we get
$$
\begin{aligned}
\dashint_{B_{y_n}(y)}\varphi(1+|D\bfu|) \,\mathrm{d}x 
& \leq 2^{\mu_2-1}\dashint_{B_{y_n}(y)}\varphi(|D\bfu-(D_n{\bf u})_{\tilde y,\vartheta^m\varrho}\otimes{\bf e}_n|)\,\mathrm{d}x+ 2^{\mu_2-1}\phi(1+|(D_n{\bf u})_{\tilde y,\vartheta^m\varrho}|)  \\
& \leq 2^{\mu_2-1}\left(c_{\eta} \vartheta^{-n} \Phi(\tilde y, \vartheta^m\varrho)+\frac{3}{2}\varphi(1+|(D_n{\bf u})_{\tilde y, \vartheta^m\varrho}|)\right) \\
& \leq 2^{\mu_2-1}\left(c_{\eta}\vartheta^{-n} \varepsilon_* + \frac 32 \right) \varphi(1+|(D_n{\bf u})_{\tilde y, \vartheta^m\varrho}|) \\
& \leq 2^{\mu_2-1}\left(c_{\eta}\vartheta^{-n} \varepsilon_* + \frac 32 \right)\dashint_{B_{\vartheta^{m}\varrho}^+(\tilde y)}\varphi(1+|D\bfu|)\,\dd x \\
& \leq 2^{\mu_2}\varphi\left(\frac{\kappa_*}{\vartheta^{m}\varrho}\right) 
\leq 2^{\mu_2}\vartheta^{\frac{1}{\mu_1}} \varphi\left(\frac{\kappa_*}{\vartheta^{m+1}\varrho}\right) 
\leq  \phi \left(\frac{\kappa_*}{y_n}\right)\,,
\end{aligned}
$$
which implies
\begin{equation}
y_n\, \phi^{-1}\Bigg( \dashint_{B_{y_n}(y)} \phi(1+|D\bfu|)\,\dd x\Bigg) \le \kappa_*\,.
\label{eq:0step3}
\end{equation}
Therefore, in view of Step 2, in particular, \eqref{morrey-interior},  with \eqref{eq:0step2} and \eqref{eq:0step3}, we obtain
\begin{equation}
\int_{B_{r}(y)}\varphi(1+|D{\bf u}|)\,\mathrm{d}x \leq \vartheta^{-\lambda} \left(\frac{r}{y_n}\right)^{\lambda}\int_{B_{y_n}(y)}\varphi(1+|D{\overline\bfu}|)\,\mathrm{d}x\,, 
\quad r< y_n\,.
\label{morrey-interior1}
\end{equation}
Now, combining the previous two estimates \eqref{morrey-boundary1} and \eqref{morrey-interior1}, we get
\[
\dashint_{B_{r}(y)}\varphi(1+|D{\overline\bfu}|)\,\mathrm{d}x
\leq c  \left(\frac{r}{\rho_0}\right)^{\lambda-n}\dashint_{B_{3\varrho_0}^+(x_0)}\varphi(1+|D{\bf u}|)\,\mathrm{d}x
\quad
\text{for every $y\in B_{\varrho_0}^+$ and $r\in (0,\varrho_0]$},
\] 
where the constant $c$ depends on $n,N,\mu_1,\mu_2,\nu,L,\alpha$, which together with the choice of $\lambda$ in \eqref{lambda} implies the $C^{\alpha}$ continuity of $\bfu$ in $B_{\varrho}^+(x_0)$. Indeed, for a.e. $x,y\in B_{\varrho}^+$ with $\varrho_1=|x-y|$, applying Poincar\'e inequality and the last inequality above, we infer
$$\begin{aligned}
|\bfu(x)-\bfu(y)|&= \sum_{i=0}^\infty |(\overline\bfu)_{B_{2^{-i}\varrho_1}(x)}-(\overline\bfu)_{B_{2^{-i+1}\varrho_1}(x)}|+\sum_{i=1}^\infty |(\overline\bfu)_{B_{2^{-i}\varrho_1}(y)}-(\overline\bfu)_{B_{2^{-i+1}\varrho_1}(y)}|\\
&\qquad + |(\overline\bfu)_{B_{2\varrho_1}(x)}-(\overline\bfu)_{B_{\varrho_1}(y)}| \\
&\leq  c \sum_{i=0}^\infty \Bigg(\dashint_{B_{2^{-i+1}\varrho_1}(x)}|\overline\bfu-(\overline\bfu)_{B_{2^{-i+1}\varrho_1}(x)}|\,\dd z + \dashint_{B_{2^{-i+1}\varrho_1}(y)} |\overline\bfu-(\overline\bfu)_{B_{2^{-i+1}\varrho_1}(y)}| \, \dd z\Bigg)\\
&\leq  c \sum_{i=0}^\infty 2^{-i}\varrho_1 \Bigg(\dashint_{B_{2^{-i+1}\varrho_1}(x)}|D\overline\bfu|\,\dd z + \dashint_{B_{2^{-i+1}\varrho_1}(y)} |D\overline\bfu|\,, \dd z\Bigg)\\
&\leq c \sum_{i=0}^\infty 2^{-i}\varrho_1 \varphi^{-1}\Bigg(\left(\frac{2^{-i+1}\varrho_1}{\rho_0}\right)^{\lambda-n}\dashint_{B_{3\varrho_0}^+(x_0)}\varphi(1+|D{\bf u}|)\,\mathrm{d}x\Bigg)\\
&\leq c \left[\sum_{i=0}^\infty \left(\frac{2^{-i}\varrho_1}{\rho_0}\right)^{\frac{\lambda-n}{\mu_1}+1}\right] \varrho_0\,\varphi^{-1}\Bigg(\dashint_{B_{3\varrho_0}^+(x_0)}\varphi(1+|D{\bf u}|)\,\mathrm{d}x\Bigg)\\
&\leq c  |x-y|^\alpha \varrho_0^{1-\alpha}\,\varphi^{-1}\Bigg(\dashint_{B_{3\varrho_0}^+(x_0)}\varphi(1+|D{\bf u}|)\,\mathrm{d}x\Bigg).
\end{aligned}$$


\noindent 
{\bf Step 4: Choice of regular points on the flat boundary.} 
Finally, we prove that if $x_0\in\Gamma_R\setminus \mathrm{Sing}_{\bu}(\Gamma_R)$ 
then, for some $\varrho_0\in(0,\frac{1}{2})$ satisfying $\varrho_0\leq \frac{\min\{\varrho_1,\varrho_2\}}{3}$ and $B^+_{3\varrho}(x_0)\subset B_R^+$, the inequalities in \eqref{eq:0step1} hold for every $\tilde y\in\Gamma_{\varrho_0}(x_0)$, hence $x_0\in \Gamma_{\bfu}$. For this, we start with assuming that   $x_0\in\Gamma$ satisfies
\begin{equation}
\mathop{\lim\inf}_{\varrho\searrow 0} \Phi(x_0,\varrho)=\mathop{\lim\inf}_{\varrho\searrow 0} \dashint_{B_\varrho^+(x_0)} {\phi}_{1+|(D_{n}{\bf u})_{x_0,\varrho}|} (|D{\bf u}-(D_{n}{\bf u})_{x_0,\varrho}\otimes {\bf e}_n|)\,\mathrm{d}x=0
\label{regularcond1}
\end{equation}
and
\begin{equation}
m_{x_0}:=\mathop{\lim\sup}_{\varrho\searrow 0}\, (|D_{n}{\bf u}|)_{x_0,\varrho}<+\infty\,.
\label{regularcond2}
\end{equation}

We first recall the constants $\epsilon_*$, $\kappa_*$, $\varrho_1$ and $\varrho_2$ determined in Step 3. We then set 
\begin{equation}
\sigma:=\left(\frac{\epsilon_* \phi(1)}{2c_2}\right)^{\mu_1}\label{eq:sigmacond},
\end{equation}
where $c_2$ is  determined in \eqref{epsilon*} below, respectively.
Then, in view of \eqref{regularcond1} and \eqref{regularcond2}, we can find $\varrho_0>0$ such that
\begin{equation}
c_2\varphi(2+m_{x_0})\left[\omega\left(\varrho_0(2+M_{x_0})\right)^{1-\frac{1}{s}}+ {\mathcal{V}}(\varrho_0)^{1-\frac{1}{s}}\right] \le \frac{\epsilon_*\phi(1)}{2},
\label{eq:radius0}
\end{equation}
\begin{equation}
\varrho_0\leq \min\left\{ \frac{\kappa_*}{2\varphi^{-1}\big(c_3\varphi(m_{x_0}+1)\big)},\frac{\varrho_1}{3},\frac{\varrho_2}{3}\right\}\,,
\label{eq:radius}
\end{equation}
where the constant $c_3>0$ is determined in \eqref{kappa*}  below,
\begin{equation}
\dashint_{B_{14\varrho_0}^+(x_0)}  |D{\bf u}-(D_n\bu)_{x_0,14\varrho_0}\otimes {\bf e}_n|\,\mathrm{d}x \leq\sigma
\quad \mbox{ and }\quad 
(|D_n\bfu|)_{x_0,14\varrho_0} \leq m_{x_0}+1\,.
\label{regcond2}
\end{equation}

Now, we prove \eqref{eq:0step1} for $\tilde y\in \Gamma_{\varrho_0}(x_0)$. By \eqref{eq:equivalence} and \eqref{eq:equivalencebis}, we have 
\begin{equation*}
\begin{aligned}
\Phi(\tilde y, 2\varrho_0) & \leq c \dashint_{B_{2\varrho_0}(\tilde y)} |\bfV_1(D\overline\bfu)-(\bfV_1(D\overline\bfu))_{B_{2\varrho_0}(\tilde y)}|^2\,\dd x\\
& \leq c \dashint_{B_{3\varrho_0}(x_0)} |\bfV_1(D\overline\bfu)-(\bfV_1(D\overline\bfu))_{B_{3\varrho_0}(x_0)}|^2\,\dd x\\
& \leq c \dashint_{B_{3\varrho_0}(x_0)} |\bfV_1(D\overline\bfu)-\bfV_1({\bm\xi}_0\otimes{\bf e}_n)|^2\,\dd x\\
& \leq c \dashint_{B_{3\varrho_0}^+(x_0)} {\phi}_{1+|\bm\xi_0|} (|D{\bf u}-{\bm\xi}_0\otimes{\bf e}_n|)\,\mathrm{d}x
\end{aligned}
\end{equation*}
for some $c_2=c_2(n,N,\mu_1,\mu_2)>0$, where $\overline \bu$ is the odd extension of $\bu$ and $\bm\xi_0=(D_{n}{\bf u})_{x_0,14\varrho_0}$. Note that by Lemma~\ref{lem:changeshift}, \eqref{eq:caccioppoli1}, \eqref{eq:2.9ok}, Lemma~\ref{lem:Jensen}
we have 
$$\begin{aligned}
\dashint_{B_{6\varrho_0}^+(x_0)} {\phi}_{1+|{\bm\xi}_0|} (|D{\bf u}-{\bm\xi}_0\otimes {\bf e}_n|)\,\mathrm{d}x & \leq c \dashint_{B_{6\varrho_0}^+(x_0)} {\phi}(|D{\bf u}|)\, \mathrm{d}x + c \phi(1+ |{\bm\xi}_0|) \\
& \leq c \dashint_{B_{7\varrho_0}^+(x_0)} {\phi}(|D_n{\bf u}|)\, \mathrm{d}x + c \phi(1+ |{\bm\xi}_0|)\\
& \leq c \Bigg(\dashint_{B_{14\varrho_0}^+(x_0)} {\phi}(|D_n{\bf u}|)^{\frac{1}{\mu_2}}\, \mathrm{d}x\Bigg)^{\mu_2} + c \phi(1+ |{\bm\xi}_0|)\\
& \leq  \Lambda \phi(1+ |{\bm\xi}_0|)
\end{aligned}$$
for some $\Lambda>0$ depending on $n,N,\mu_1,\mu_2,\nu,L$.
Therefore, from the preceding estimate, we can exploit the reverse H\"older estimate \eqref{eq:caccioppoliIbis1} in Lemma~\ref{corollary3.2}  with $t=\frac{1}{\mu_2}$. Using this, Lemma~\ref{lem:Jensen}, and \eqref{regcond2}, 
we obtain
\begin{equation}\begin{aligned}
&\Phi(\tilde y, 2\varrho_0)\leq c \dashint_{B_{3\varrho_0}^+(x_0)}  \varphi_{1+|\bm\xi_0|}(|D{\bf u}-\bm\xi_0\otimes {\bf e}_n|)\,\mathrm{d}x
 \\
&\leq c \varphi_{1+|\bm\xi_0|}\Bigg(\dashint_{B_{6\varrho_0}^+(x_0)}  |D{\bf u}-\bm\xi_0\otimes {\bf e}_n|\,\mathrm{d}x\Bigg)+ c\varphi(1+|{\bm \xi}_0|)\left[\omega\left(\varrho_0(1+|{\bm \xi_0}|)\right)^{1-\frac{1}{s}}+ {\mathcal{V}}(\varrho_0)^{1-\frac{1}{s}}\right]\,\\
&\leq c_2 \sigma^{\mu_1}+ c_2\varphi(2+m_{x_0})\left[\omega\left(\varrho_0(2+M_{x_0})\right)^{1-\frac{1}{s}}+ {\mathcal{V}}(\varrho_0)^{1-\frac{1}{s}}\right],
\label{epsilon*}
\end{aligned}\end{equation}
where $c_2>0$ depends only on $n,N,\mu_1,\mu_2,\nu,N$.
Then, taking into account \eqref{eq:sigmacond} and \eqref{eq:radius0}, \eqref{epsilon*} implies
$$
\frac{\Phi(\tilde y, 2\varrho_0)}{\phi(1+|(D_n\bfu)_{\tilde y, 2\varrho_0}|)} \le  \varepsilon_{*} \,.
$$
This is the first inequality in \eqref{eq:0step1}.


Moreover, by \eqref{eq:caccioppoli1}, \eqref{eq:2.9ok}, \eqref{highint1} with $t=\frac{1}{\mu_2}$, Lemma~\ref{lem:Jensen},  \eqref{regcond2}, we have 
\begin{equation}\label{kappa*}
\begin{split}
\dashint_{B_{2\varrho_0}^+(\tilde y)}\varphi(1+|D\bfu|)\,\mathrm{d}x & \leq c \dashint_{B_{3\varrho_0}^+(x_0)}\varphi(1+|D\bfu|)\,\mathrm{d}x  \leq c \dashint_{B_{4\varrho_0}^+(x_0)}\varphi(1+|D_n\bfu|)\,\mathrm{d}x\\
& \leq c+c \varphi\Bigg(\dashint_{B_{8\varrho_0}^+(x_0)}|D_n\bfu|\,\mathrm{d}x\Bigg) \le c_3 \phi(1+m_{x_0})
\end{split}
\end{equation}
for some $c_3>0$ depending only on $n,N,\mu_1,\mu_2,\nu,\Lambda$. Therefore, by \eqref{eq:radius} we obtain
$$
\dashint_{B_{2\varrho_0}^+(\tilde y)}\varphi(1+|D\bfu|)\,\mathrm{d}x \leq \varphi\left(\frac{\kappa_*}{2\varrho_0}\right).
$$
This implies the second inequality in \eqref{eq:0step1}, and the proof is concluded.
\endproof

We are now in position to prove the main result of the paper.

\proof[Proof of Theorem~\ref{theorem-result-1}] \

We convert the minimizing problem in Theorem~\ref{theorem-result-1} on a boundary region near each boundary point  into the one in Theorem~\ref{theorem-result-2} by a standard flattening argument.  Then Theorem~\ref{theorem-result-1} follows from Theorem~\ref{theorem-result-2}.

Since $\partial \Omega\in C^1$, for each $\tilde x\in \partial \Omega$   there exist a coordinate system $x=(x_1,\dots,x_n)$ with the origin at $\tilde x$ and $\bm\nu_{\tilde x}=(0,\dots,0,1)$ and, in this coordinate system, a $C^1$ function $\gamma:\R^{n-1}\to\R$ such that $B_{r_0}(0)\cap\Omega=\{x=(x',x_n)\in B_r(0): x_n>\gamma(x') \}$ for some sufficiently small  $r_0>0$, and $\gamma_{x_i}(0)=0$, $i=1,2,\dots,n-1$.  Note that  $\sup\{|D'\gamma(x')|\,:\, |x'|\le r\}  \le M_r$ for every $\tilde x\in \Omega$ with $M_r>0$ satisfying that $\lim_{r\to 0}M_r=0$, where $D'\gamma=(\gamma_{x_1},\dots,\gamma_{x_{n-1}})$. From now on, we fix $\tilde x\in\partial\Omega$, the relevant coordinate system and the $C^1$-function $\gamma$. Note that since the functional \eqref{functional} and the assumptions in Theorem~\ref{theorem-result-1} are invariant under rotations of the coordinate system, without loss of generality we write the relevant coordinate system as $x=(x_1,\dots,x_n)$.

Then we define $\bT:\R^n\to\R^n$ by $\bT(x',x_n)=(x',x_n-\gamma(x'))$ and its inverse by $\bT^{-1}(y)=(y',y_n+\gamma(y'))$.  We choose sufficiently small $R\in (0,r_0/4)$ such that 
{for every $x\in B_{2R}$, $y\in \bT(B_{2R})$ and $r\in (0,2R]$,
\begin{align}
&\frac{1}{4}|B_r| \le |\Omega_r(x)|\,; \label{T0}\\
&2^{-1} \le |D\bT(x)|\le  2\,, \quad \text{hence } \ 2^{-1}\le |D\bT^{-1}(y)| \le  2\,; \label{DT}\\
&B_{r/2}(\bT(x))  \subset \bT(B_{r}(x))\subset   B_{2r}(\bT(x))\,,\quad
B_{r/2}(\bT^{-1}(y))  \subset \bT^{-1}(B_{r}(y))\subset   B_{r}(\bT^{-1}(y))\,. \label{T1}
\end{align}
Note that $B_{R}^+ \subset \bT(B_{2R}\cap \Omega)$, $|U|\sim |\bT(U)|$ and $|V|\sim |\bT^{-1}(V)|$.}
In addition,
since $D\bT$, $D\bT^{-1}$ and $D\bg$ are uniformly continuous, we denote
$$
\omega_2(r):=\sup_{|x_1-x_2|\le r} |D\bT(x_1)-D\bT(x_2)|, \quad \omega_3(r):=\sup_{|y_1-y_2|\le r} |D\bT^{-1}(x_1)-D\bT^{-1}(x_2)|,
$$
$$
\omega_4(r):=\sup_{|x_1-x_2|\le r} |D\bg(x_1)-D\bg(x_2)|.
$$
Note that $\omega_2,\omega_3,\omega_4$ are non-decreasing, vanishing at $0$, and, without loss of generality, assumed to be concave. Hence, we have $\omega_i(cr)\le c\omega_i(r)$, $i=2,3,4$, for all $c\ge 1$ and $r>0$.

Now we set
\[
\tilde \bfu(y):= \bfu(\bT^{-1}(y))-\bg(\bT^{-1}(y))
\]
and
\[
\tilde f(y, {\bf w}, \bQ):= |D{\bT}^{-1}(y)| \, f(\bT^{-1}(y), {\bf w}+\bg(\bT^{-1}(y)),   \bQ\, D\bT(\bT^{-1}(y))+ D\bg(\bT^{-1}(y)) ) .
\] 
Then we see that $\tilde \bfu$ is a minimizer of 
\[
{\bf v}\in W^{1,\phi}_\Gamma(B^+_{R})\ \ \mapsto \ \ \int_{B_{R}^+} \tilde f(y,{\bf v},D{\bf v})\,\mathrm{d}y\,.
\]

Moreover, $\tilde f$ satisfies the analog of \ref{ass-1f} -- \ref{ass-6f} in the setting of the $y$-coordinate system with $B_R^+$ in place of $\Omega$,
{where relevant constants $\nu$ and $L$ depend on the ones for $f$, $n$, $N$, $\mu_1$, $\mu_2$, and $\|D\bg\|_{\infty}$. }
 We only show that $\tilde f$ satisfies the VMO assumption \ref{ass-4f}. The other assumptions are relatively easy to check, since they are pointwise conditions, and we refer to \cite{Be09,Ok3} for the proof.

Now, we start proving \ref{ass-4f} for $\tilde f$. Let $r\in (0,R]$ and  $y_0,y\in \overline{B_R^+}$, and  set $x_0:=\bT^{-1}(y_0)$, $x:=\bT^{-1}(y)$,
$$
{\bf w}_y:={\bf w}+\bg(\bT^{-1}(y)), \quad \bQ_y:=\bQ\, D\bT(\bT^{-1}(y))+ D\bg(\bT^{-1}(y)).
$$
$$
{\bf u}_x:={\bf w}_{\bT(x)}={\bf w}+\bg(x), \quad \bP_x:=\bQ_{\bT(x)}=\bQ\, D\bT(x)+ D\bg(x).
$$
Then we first observe that 
\begin{equation*}
\begin{aligned}
&| \tilde f (y, {\bf w}, {\bf Q})- (\tilde f( \cdot, {\bf w}, {\bf Q}))_{B_{r}(y_0) \cap B_R^+}|\\
&= \bigg||D{\bT}^{-1}(y)| \, f(x, {\bf u}_x,   \bP_{x})-\dashint_{B_{r}(y_0) \cap B_R^+} f (\bT^{-1}(\tilde y), {\bf w}_{\tilde y}, {\bf Q}_{\tilde y})|D{\bT}^{-1}(\tilde y)| \,\dd \tilde y\bigg|\\
&\leq  ||D{\bT}^{-1}(y)|-(|D{\bT}^{-1}|)_{B_{r}(y_0) \cap B_R^+}| | f(x, {\bf u}_x,   \bP_{x})|\\
&\qquad +\bigg| \dashint_{B_{r}(y_0) \cap B_R^+} (f(x, {\bf u}_x,   \bP_{x})-f (\bT^{-1}(\tilde y), {\bf u}_{x}, {\bf P}_{x}))|D{\bT}^{-1}(\tilde y)| \,\dd \tilde y\bigg|\\
&\qquad +\bigg| \dashint_{B_{r}(y_0) \cap B_R^+} (f (\bT^{-1}(\tilde y), {\bf u}_{x}, {\bf P}_{x})-f (\bT^{-1}(\tilde y), {\bf w}_{\tilde y}, {\bf Q}_{\tilde y}))|D{\bT}^{-1}(\tilde y)| \,\dd \tilde y\bigg|\\
&=: I+II+III.
\end{aligned}
\end{equation*}
We then estimate the three terms  $I$, $II$, $III$, separately. By \ref{ass-1f}, \eqref{DT}, and the definition of $\omega_3$, we have
$$
I\leq  c \,\omega_3(2r) \phi(1+|\bQ|+\|D\bg\|_\infty) \le c (1+\|D\bg\|_\infty)^{\mu_2}\omega_3(2r) \phi(1+|\bQ|)\,.
$$
As for $II$, by the change of variable $\tilde x=\bT^{-1}(\tilde y)$, \eqref{T0} -- \eqref{T1} and  \ref{ass-4f}, we obtain
$$\begin{aligned}
II&=\bigg| \frac{1}{|B_{r}(y_0) \cap B_R^+|}\int_{\bT^{-1}(B_{r}(y_0) \cap B_R^+)} f(x, {\bf u}_x,   \bP_{x})-f (\tilde x, {\bf u}_{x}, {\bf P}_{x}) \,\dd \tilde x\bigg|\\
&\le c\dashint_{\Omega_{2r}(x_0)} |f(x, {\bf u}_x,   \bP_{x})- (f(\cdot, {\bf u}_x,   \bP_{x}))_{\Omega_{2r}(x_0)}|+|f (\tilde x, {\bf u}_{x}, {\bf P}_{x})- (f(\cdot, {\bf u}_x,   \bP_{x}))_{\Omega_{2r}(x_0)}| \,\dd \tilde x \\
&\le c\Bigg( v_{x_0,2r}(x)+\dashint_{\Omega_{2r}(x_0)} v_{x_0,2r}(\tilde x) \,\dd \tilde x\Bigg)  \phi(1+|\bQ|+\|D\bg\|_\infty) \\
& \le c (1+\|D\bg\|_\infty)^{\mu_2} ( v_{\bT^{-1}(y_0),2r}\left(\bT^{-1}(y)) |D\bT^{-1}(y)| +\mathcal{V}(2r) \right)  \phi(1+|\bQ|)\,.
\end{aligned}$$
Finally, with the change of variable $\tilde x=\bT^{-1}(\tilde y)$ and a similar argument as in the preceding inequalities for $II$, we get the estimate
$$
III\le c \dashint_{\Omega_{2r}(x_0) } |f (\tilde x, {\bf u}_{x}, {\bf P}_{x})-f (\tilde x, {\bf u}_{\tilde x}, {\bf P}_{x})|+|f (\tilde x, {\bf u}_{\tilde x}, {\bf P}_{x})-f (\tilde x, {\bf u}_{\tilde x}, {\bf P}_{\tilde x})| \,\dd \tilde x\,.
$$
Note that, by \ref{ass-5f}, \eqref{(1.3celok)}, \eqref{DT} and the definitions of $\omega_2$ and $\omega_4$,
$$\begin{aligned}
|f (\tilde x, {\bf u}_{x}, {\bf P}_{x})-f (\tilde x, {\bf u}_{\tilde x}, {\bf P}_{x})| 
&\leq c \omega( |\bg(x)-\bg(\tilde x)|)\phi(1+|\bQ|+\|D\bg\|_\infty)\\
&\leq c(1+\|D\bg\|_\infty)^{\mu_2+1} \omega(r)\phi(1+|\bQ|)\,,
\end{aligned}$$
and
$$\begin{aligned}
|f (\tilde x, {\bf u}_{\tilde x}, {\bf P}_{x})-f (\tilde x, {\bf u}_{\tilde x}, {\bf P}_{\tilde x})| 
& \leq c\phi'(|{\bf P}_{x}|+|{\bf P}_{\tilde x}|) |{\bf P}_{x}-{\bf P}_{\tilde x}|\\
& \leq c \phi'(1+|\bQ|+\|D\bg\|_\infty)\, |\bQ\omega_2(r)+
\omega_4(r)|\\
& \leq c (1+\|D\bg\|_\infty)^{\mu_2-1}(\omega_2(r)+
\omega_4(r)) \phi(1+|\bQ|)\,.
\end{aligned}$$
Combining the preceding three estimates, we then have
$$
III\le c (1+\|D\bg\|_\infty)^{\mu_2+1} (\omega(r)+\omega_2(r)+\omega_4(r)) \phi(1+|\bQ|)\,.
$$
Consequently, collecting the estimates for $I$, $II$ and $III$, we obtain
$$\begin{aligned}
| \tilde f (y, {\bf w}, {\bf Q})- (\tilde f( \cdot, {\bf w}, {\bf Q}))_{B_{r}(y_0) \cap B_R^+}| \le  \tilde v_{y_0,r}(y) \phi(1+|\bQ|),
\end{aligned}$$ 
where
$$\begin{aligned}
\tilde v_{y_0,r}(y) := &c (1+\|D\bg\|_\infty)^{\mu_2+1} \\
&\times \left\{v_{\bT^{-1}(y_0),2r}(\bT^{-1}(y)) |D\bT^{-1}(y)| +\mathcal{V}(2r)+\omega(r)+\omega_2(r)+\omega_3(r)+\omega_4(r)\right\}.
\end{aligned}$$
Moreover, since 
$$\begin{aligned}
\dashint_{B_r(y_0)\cap B^+_R}v_{\bT^{-1}(y_0),2r}(\bT^{-1}(y)) |D\bT^{-1}(y)|\,\mathrm{d}y  
&= \frac{1}{|B_r(y_0)\cap B^+_R|} \int_{\bT^{-1}(B_r(y_0)\cap B^+_R)}v_{x_0,2r}(x) \,\mathrm{d}x\\
& \le c \dashint_{\Omega_{2r}(x_0)}v_{x_0,2r}(x) \,\mathrm{d}x \le c \mathcal V(2r)\,,
\end{aligned}$$
for every $B_r(y_0)$ with $y_0\in \overline {B^+_R}$, 
$$
\lim_{r\to 0} \underbrace{\Bigg[\sup_{0<\varrho\le r} \sup_{y_0\in \overline{B^+_{R}}} \dashint_{B_\varrho(y_0)\cap B^+_R} v_{y_0,r} (y)\,\dd y\Bigg]}_{=:\tilde{\mathcal V}(r)} \le \lim_{r\to0} c\left(\mathcal{V}(2r)+\omega(r)+\omega_2(r)+\omega_3(r) \right)=0.
$$
This implies the validity of assumption \ref{ass-4f} for $\tilde f$. 
{Here, we stress that $\|D\bg \|_\infty$ can be easily handled due to the term ``1+" which implies the non-degenerate condition. Within the degenerate setting, instead, 
the flattening argument above seems quite unclear.
}

Finally, we will prove that
\begin{equation}
0 \not\in \mathrm{Sing}_{\tilde{\bu}}(\partial\Omega)
 \quad  \Longleftrightarrow \quad  
\bT(0) = 0 \not\in \mathrm{Sing}_{\tilde{\bu}}(\Gamma_r)\,.
\label{eq:implication}
\end{equation}

We only show the implication ``$\Rightarrow$'', as the other one can 
be obtained in a similar way. Note that ${\bm\nu}_0={\bf e}_n$ and $\bT({\bm\nu_0})=\bT({\bf e}_n)={\bf e}_n$. 

{
We first show that
$$
\mathop{\lim\sup}_{\varrho\searrow 0}\,(|D_{n}\tilde {\bf u}|)_{\varrho}<\infty\,. 
$$
For this, we preliminarily note that
$$
\begin{aligned}
(|D_{n}{\tilde{\bf u}}|)_{\varrho}
\le \dashint_{B^+_{\varrho}}|D_n\bu(\bT^{-1}(y))D_n\bT^{-1}(y)|\,\dd y  + \dashint_{B^+_{\varrho}}|D_n\bg(\bT^{-1}(y))D_n\bT^{-1}(y)|\,\dd y \,. 
\end{aligned}
$$
Now, the second integral on the right hand side converges as $\varrho\to 0$ since $D\bg$, $\bT^{-1}$ and $D_n\bT^{-1}$ are continuous at $0$. As for the first integral, by the change of variable $x=\bT^{-1}(y)$, we have
$$
\dashint_{B^+_{\varrho}}\left|D_n\bu(\bT^{-1}(y))D_n\bT^{-1}(y)\right|\,\dd y \leq c \dashint_{\Omega_{2\varrho}}|D_n\bu(x)|\,|D_n\bT^{-1}(\bT(x))|\,|D\bT(x)|\,\dd x\leq c \dashint_{\Omega_{2\varrho}}|D_n\bu(x)|\,\dd x.
$$
Therefore, since  the limsup of the right hand side above is finite  as $\varrho\to 0$, so is the limsup of $(|D_{n}{\tilde{\bf u}}|)_{\varrho}$.

We now show that 
$$
\mathop{\lim\inf}_{\varrho\searrow 0} \dashint_{B^+_\varrho} |D\tilde{\bf u}-(D_{n}\tilde{\bf u})_{\varrho}\otimes {\bf e}_n|\,\mathrm{d}y=0\,.
$$
Since $\tilde \bu={\bf 0}$ on $\Gamma_R$, we extend $\tilde \bu$ to $B_R$ such that $\tilde\bu(y',y_n)=-\tilde\bu(y',-y_n)$. Then, with $\bP_0={\bf Q}_{{\bf T}(0)}:=[(D\bu)_{\Omega_{2\varrho}}-(D\bg)_{\Omega_{2\varrho}}]D\bT^{-1}(\bT(0))$, we have
$$\begin{aligned}
&\dashint_{B^+_\varrho} |D\tilde{\bf u}-(D_{n}\tilde{\bf u})_{\varrho}\otimes {\bf e}_n|\,\mathrm{d}y 
= \dashint_{B_\varrho} |D\tilde{\bf u}-(D\tilde{\bf u})_{B_\varrho}|\,\mathrm{d}y
\le  c  \dashint_{B_\varrho} |D\tilde{\bf u}-\bP_0|\,\mathrm{d}y  =  c  \dashint_{B^+_\varrho} |D\tilde{\bf u}-\bP_0|\,\mathrm{d}y  \\
&\le c \dashint_{\Omega_{2\varrho}} |[D{\bf u}(x)-D\bg (x)]D\bT^{-1}(\bT(x))-\bP_0|\,\mathrm{d}x \\
&\le  c |D\bT^{-1}(\bT(0))|\dashint_{\Omega_{2\varrho}} |D{\bf u}(x)-(D\bu)_{\Omega_{2\varrho}}|+|D\bg (x)-(D\bg)_{\Omega_{2\varrho}}| \,\mathrm{d}x\\
&\quad +c \omega_3(\varrho)\dashint_{\Omega_{2\varrho}} |D\bg (x)| \,\mathrm{d}x +  c \omega_3(\varrho)\dashint_{\Omega_{2\varrho}} |D\bu-(D_n\bu)_{\Omega_{2\varrho}}\otimes {\bf e}_n| \,\mathrm{d}x+ c \omega_3(\varrho)(|D_n\bu|)_{\Omega_{2\varrho}}\,,
\end{aligned}$$
where $\omega_3$ is the modulus of continuity of $D\bT^{-1}$. 
Since $0\not\in \mathrm{Sing}_{\bu}(\partial\Omega)$ and $D\bg$ is continuous, the liminf of the right hand side as $\varrho\to0$ is zero, and so is the liminf of the left hand side.
}

\smallskip

\noindent {\bf Conclusion.} Suppose $x_0\not\in \mathrm{Sing}_{\tilde{\bu}}(\partial\Omega)$. We set $x_0=0$. Then by \eqref{eq:implication}, $\bT(x_0)=0\not\in \mathrm{Sing}_{\tilde{\bu}}(\Gamma_R)$. Therefore, in view of Theorem~\ref{theorem-result-2}, we infer that $\tilde \bu$ is locally $C^{0,\alpha}$ at $0$ for every $\alpha\in(0,1)$ in the $y$-coordinate system, which together with the $C^1$ regularity of $\bT$ and $\bg$ implies that 
$\bu$ is locally $C^{0,\alpha}$ at $x_0=0$ for every $\alpha\in(0,1)$ in the $x$-coordinate system, hence $x_0\in \Omega_{\bu}$. 
\endproof

\section*{Acknowledgments} 

{J.\ Ok was supported by the National Research Foundation of Korea funded by the
Korean Government (NRF-2022R1C1C1004523).}
G. Scilla has been supported by the Italian Ministry of Education, University and Research through the MIUR – PRIN project 2017BTM7SN “Variational methods for stationary and evolution problems with singularities and interfaces”. The research of B. Stroffolini was supported by PRIN Project 2017TEXA3H.










\end{document}